\documentclass[envcountchap]{svmono}

\usepackage{makeidx}     
\usepackage{graphicx}    
\usepackage{amsfonts,amsmath}
\renewcommand{\ge}{\varepsilon}
\DeclareMathOperator{\rge}{rge}
\newcommand{\R}{\mathbb{R}}
\newcommand{\Rn}{{\mathbb{R}^n}}
\newcommand{\de}{\partial}
\newcommand{\weakto}{\rightharpoonup}
\DeclareMathOperator{\essinf}{ess\, inf}
\DeclareMathOperator{\cat}{cat}
\let\eps=\ge
\newcommand{\N}{\mathbb{N}}
\newcommand{\C}{\mathbb{C}}
\newcommand{\re}{\operatorname{Re}}

\newcommand{\irn}{\int\limits_{\Rn}}

\newcommand{\Wn}{W^{1,2}(\Rn)}

\newcommand{\dw}{\dot{w}}
\newcommand{\dz}{\dot{z}}
 
\newcommand{\wz}{\widetilde z}

\renewcommand{\a }{\alpha }
\renewcommand{\b }{\beta }
\renewcommand{\d }{\delta }
\let\e=\ge

\newcommand{\n }{\nabla }

\makeindex             

\begin{document}

\author{Simone~Secchi}
\title{Nonlinear differential equations on noncompact domains}
\subtitle{Ph.D. Thesis}
\maketitle

\frontmatter

%
%
%

\thispagestyle{empty}
\vspace*{3.5cm}
\begin{flushright}

Vorrei esprimere la mia gratitudine a tutti le persone che mi sono state vicine
in questi anni: la mia famiglia, i miei amici, e tutti quelli che hanno discusso
con me di matematica. Un ringraziamento particolare va al mio supervisore,
il prof. Antonio Ambrosetti. Da lui ho imparato tutto quello che so sulla ricerca
scientifica, e non solo le tecniche.

\end{flushright}
%
%
%

\preface

The aim of this Ph.d.~thesis is to present some very recent results
 concerning specific problems that share a common feature: a certain
 kind of \emph{non--compactness}.

As explained in the first chapter, while \emph{compact} problems have a
 well-established existence theory, much less is known for non--compact
 problems which have attracted a great attention in recent years.
 Compactness usually breaks down for topological reasons, for example
 because an equation should be solved on an unbounded domain. Hence the
 elementary tools of variational analysis cannot be applied, and one is
 compelled to find new ways out.

We present three possible approaches to non--compact problems:
\begin{description}
\item[(1)] Concentration--Compactness;
\item[(2)] Perturbation methods;
\item[(3)] Bifurcation theory.
\end{description}

The first one was introduced by P.~L.~Lions in 1984
 (\cite{Lions85-1,Lions85-2}), and immediately adapted to a wide range
 of differential equations and, more generally, variational problems. As
 an application, we study here the existence of best constants for a
 class of Hardy--like inequalities in $\Rn$. See section 1.7. The main
 result is Theorem 1.11. This result will appear in our forthcoming
 paper \cite{SeSmWi}.

\vspace{20pt}

The second approach was introduced in recent years by Ambrosetti and
 Badiale (see \cite{AmbBad1,AmbBad2}) as a generalization of previous
results by Ambrosetti, Coti Zelati and Ekeland (\cite{ACE}). It is still 
a ``young'' tool,  though several applications have already appeared in
the literature. We present new results concerning the existence of
closed geodesics on infinite cylinders, following our paper \cite{Se},
and to some Schr\"{o}dinger--like elliptic partial differential
equations, following \cite{ams}. Moreover, applications to a different
problem where a vector potential function appears as a perturbation term
are described, following \cite{CS}. The main results are Theorems 21 and
23 in chapter 3, Theorem 25 in chapter 4 and Theorems 30 and 31 in
chapter 5. More comments on the history and background material on these
topics are  contained in the next chapters.

\vspace{20pt}

Finally, Bifurcation Theory has also been used to prove existence of
 nontrivial solutions for some classes of nonlinear elliptic equaitons
 on $\Rn$. More recently, Fitzpatrick and Pejsachowicz have defined a
 degree for Fredholm maps to find more general bifurcation results for
 such equations on $\Rn$. Following \cite{SeStu}, we discuss in chapter
 6 an application of this degree to hamiltonian systems.

\bigskip

At last, let me acknowledge the support I have been given during the
last four years at S.I.S.S.A. and the hospitality of the universities I
have been visiting. In particular, I which to thank Dr. Silvia
Cingolani (Politecnico di Bari), Prof. Charles Stuart (\'Ecole
Polytechnique F\'ed\'erale, Lausanne), and Prof. Michel Willem
(Universit\'e Catholique de Louvain-la-Neuve).

\vspace{1cm}
\begin{flushright}\noindent
Trieste\hfill {\it Simone Secchi}\\
August 2002\\
\end{flushright}

\tableofcontents

\mainmatter
\chapter*{Notation}

We collect here a list of notation commonly used in this thesis. We have tried to stick to the rule of \textit{...with an obvious meaning of symbols}.

\begin{itemize}
\item $X$ will be a real Banach space, and $E$ a real Hilbert space.
\item The Fr\'{e}chet derivative of a differentiable function $J\colon
E\to\R$ at a point $u\in E$ will be denoted by $DJ(u)$. Sometimes, to
have more readable formulas, we use $J'(u)$, and similarly $D^2J$ or
$J''$ for second-order derivatives. Since $DJ(u)$ is a linear map, we frequently omit brackets, {\it viz.} $DJ(u)v$ stands for the action of the linear map $DJ(u)$ on the vector $v$.
 \item Partial derivatives
of functions $u\colon\Rn\to\R$ will be denoted by either $\de_i u$, or
$\de u / \de x_i$ when the name of the independent variable is
important in the context.\footnote{Of course, the first choice is more
intrinsic, though the second is often more explicit.}
\item The symbol $\Delta$ stands for the ordinary Laplace operator: if $u\colon \R^n\to\R$ is a regular function, then
\[
\Delta u = \sum_{i=1}^n \frac{\de^2}{\de x_i^2}u
\]
or, equivalently, the trace of the Hessian matrix of $u$.
\item It is
always a delicate matter to fix notation for multiplication of real
numbers and scalar products in Hilbert spaces. We are not going to say
the ultimate word here. Inner products should be written as $\langle
u,v\rangle$ whenever possible, but we will use $u\bullet v$ when
dealing with geometric problems, and also $u \cdot v$ from time to
time.  \item The complex conjugate of any number $z\in\C$ will be
denoted by $\bar z$. \item The real part of a number $z\in\C$ will be
denoted by $\re z$. \item The ordinary inner product between two
vectors $a,b\in\Rn$ will be denoted by $a \cdot b$. 
\item Integration with respect to the ordinary Lebesgue measure is denoted
by $\int \dots dx$. This is a little abuse of notation, since nobody denotes
this measure by $x$! A more correct $\int \dots d\mathcal{L}$ shoud be used
throughout, but we refrain from this ``calvinism''. Even worse, when no 
confusion can arise, we omit the symbol $dx$ in integrals
over $\Rn$. \item $C$ denotes a generic positive constant, which may
vary inside a chain of inequalities. \item We use the Landau symbols.
For example $O(\e)$ is a generic function such that
$\limsup\limits_{\e\to 0} [O(\e)/\e] < \infty$, and $o(\e)$ is a
function such that $\lim\limits_{\e\to 0} [o(\e)/\e]=0$. Of
course, a function can depend on several parameters, and its behavior
is frequently different according to which parameter is considered. In
these cases, we use a subscript to single out the parameter that 
``tends to..."
\item The end of a proof is denoted by the symbol \qed.\footnote{This
choice is not completely standard. I like it because the symbol \qed
\quad can be immediately seen inside a page.}

\end{itemize}

%
%
%

\part{The importance of being compact}

%
%
%

\chapter{Why compactness?}
\label{intro} 

The r\^{o}le of compactness in Modern Analysis can hardly be
 overrated. From the basic theorem of Weierstrass stating that any
continuous function on a compact space attains a maximum and a minimum,
to the most advanced results in Calculus of Variations, compactness
properties are exploited to prove existence of optimization problems,
ordinary and partial differential equations, regularity theory, and much
more.

\section{Variational methods: how do they work?}

Before we can actually launch ourselves into the study of non--compact
 problems, we need to develop some tools of modern\footnote{The word ``modern"
 is quite abused. Its meaning depends on the context. Roughly speaking, we
 might replace it by ``functional--analytic'', at least in this work.} analysis. Since this
 is not a treatise on nonlinear functional analysis, we shall content
 ourselves with a brief survey of the main ideas. Probably the best
 example is that of elliptic partial differential equations (PDE's).
 While a nice feature of ordinary differential equations is that they
 can often be solved explicitely or numerically, the situation changes
 dramatically for PDE's. For example, take a bounded domain
 $\Omega\subset\Rn$, and a right-hand side $f\in C(\overline{\Omega})$.
 It can be
 shown that the innocent--looking problem
\begin{equation}\label{eq:laplace-example}
\begin{cases}
-\Delta u = f &\text{ in }\Omega\\
\phantom{aaa}u=0 &\text{ on } \de\Omega
\end{cases}
\end{equation}
has, in general,\footnote{With respect to the dimension $n$ and to the
domain $\Omega$.} no solution of class $C^2$. And although there are
 explicit representations of solutions when $\Omega$ has a particular
 shape (a ball, a half-space), nevertheless these formulas are seldom
 useful for both theoretic and practical purposes.

The natural conclusion is that one should probably decide to weaken the
 concept of solution: this is the concept of \emph{weak solutions}.
 Again, instead of presenting a thorough treatment of the general case,
 we prefer to focus on the previous example. Assume $u\in C^2 (\Omega)$
 is a classical solution of \eqref{eq:laplace-example}. Multiply both sides of the
 equation by any $\varphi\in C_0^\infty (\Omega)$, and use the Green
 formula (see Appendix). Since $u$ is zero on the boundary,
\begin{equation}\label{def:weak}
\int_\Omega \nabla u \cdot \nabla \varphi \, dx = \int_\Omega f \varphi
 \, dx.
\end{equation}
Remark that \eqref{def:weak} makes sense for any $u\in L_{loc}^1
 (\Omega)$. We could thus define \emph{distributional solutions} $u\in
 L_{loc}^1 (\Omega)$, but we will find the next concept more useful.

\begin{definition}
We say that $u\in H^1(\Omega)$ is a weak solution of
 \eqref{eq:laplace-example} if \eqref{def:weak} is true for all
 $\varphi\in H^1(\Omega)$.
\end{definition}

Now, we can present the core of the variational approach to problems
 like \eqref{eq:laplace-example}.

\begin{description}
\item[(A)] Look for a weak solution, and possibly for its uniqueness.
\item[(B)] Prove that this solution is actually of class $C^2$.
\item[(C)] Prove that any weak solution of class $C^2$ is a classical
 solution, i.e. \eqref{eq:laplace-example} is true pointwise.
\end{description}

Step (C) is usually easy, whereas step (B) is rather difficult for
 non--linear equations. We won't go into details.

At this stage, how can we treat (A)? At a first glance, \eqref{def:weak}
 is worse than the original equation, but we can remark that it can be
 seen from many viewpoints. For example as an orthogonality relation, or
 as the fact that the differential of the functional $J\colon H_0^1
 (\Omega)\to\R$ is zero at $u$, where
\[
J(u)=\frac{1}{2}\int_\Omega |\nabla u|^2\, dx -\int_\Omega fu \, dx.
\]

The former viewpoint immediately reminds us of the Lax--Milgram theorem,
 which is indeed a powerful tool for linear elliptic equations (see
 \cite{renardy});
the latter viewpoint opens the way to \emph{critical point theory},
 which is the art of finding zeroes of derivatives in Banach or Hilbert
 spaces. We will say some words in the next section, but we refer to the
 literature for extensive tretments
 (\cite{Amb-paris,Chang,MawWil,Struwe}). The main advantage of this
 second approach is that it works for a wide class of semilinear
 elliptic equations.

\section{A standard example of the use of compactness}
\label{sec:1}  

At this stage it is not clear why compactness plays a r\^{o}le in all
 this
 stuff. To clarify the situation, we present a popular approach to some
 equations, i.e. constrained minimization.

\begin{theorem}\label{th:1}
Let $X$ be a reflexive Banach space and a functional $J\colon X\to\R$ of
class $C^1$ be given. Let $u\in X$ be a point of absolute minimum of
$J$. Then
\[
DJ(u)=0 \mbox{\quad in } X^\star,
\]
where $X^\star$ is the dual space of $X$.
\end{theorem}

\begin{definition}
A functional $J\colon X\to\R$ is called boudned from below on $X$ if
there exists a constant $M>0$ such that
\[
J(u)\geq M
\]
for all $u\in X$.
\end{definition}

\begin{definition}
A functional $J\colon X\to\R$ is called coercive on $X$ if\footnote{More
 generally, 
a functional $J\colon X \to \R$ defined on the topological space $X$ is
 coercive if
the set $\{x\in X \mid J(x) \leq t\}$ is relatively compact in $X$ for
 all $t\in\R$.} 
\[
\lim_{\|u\|\to\infty} J(u)=+\infty.
\]
\end{definition}

\begin{theorem}
Let $X$ be a reflexive Banach space. Let a functional $J\colon X\to\R$
be
\begin{enumerate}
\item bounded from below on $X$;
\item weakly lower semicontinuous on $X$;
\item coercive on $X$.
\end{enumerate}
Then $J$ has a minimum.
\end{theorem}

\begin{proof}
Let $\{u_n\}$ be a minimizing sequence for $J$, i.e.
\[
\lim_{n\to\infty} J(u_n)=\inf_{u\in X} J(u) =:I.
\]
In particular, since $J$ is bounded from below, $I>-\infty$, and there
exists a constant $C$ such that $|J(u_n)| \leq C$ for all
$n\in\mathbb{N}$. Since $J$ is coercive, $\{u_n\}$ is bounded in the
norm of $X$. Since $X$ is reflexive, bounded sequences in $X$ are weakly
relatively compact, so
 that, up to a subsequence, $u_n \to u$ weakly as
$n\to\infty$. By the semicontinuity,
\[
I=\lim_{n\to\infty} J(u_n) \geq J(u),
\]
thus proving that $u\in X$ is a minimum point. Notice that the compactness assumptions has given us the element that has turned out to be the minimum of $J$.
\qed \end{proof}

\begin{example}
We come back to \eqref{eq:laplace-example}. Let $\Omega\subset \Rn$ be a
 bounded domain, with smooth boundary
$\de\Omega$. Let $H_0^1(\Omega)$ be the ordinary Sobolev space (see
appendix). For a given $f\in L^2(\Omega)$, define
\begin{equation*}
\begin{aligned}
J \colon &H_0^1(\Omega) \longrightarrow \R \\
&u \mapsto \frac{1}{2} \int_\Omega |\nabla u|^2\, dx
-\int_\Omega fu\, dx.
\end{aligned}
\end{equation*}

It is easy to prove that $J$ is weakly lower semicontinuous. Moreover,
from the Poincar\'e inequality it follows that $J$ is coercive and
bounded from below. Hence it attains a minimum at some $\bar{u}\in
H_0^1(\Omega)$. $J$ is trivially of class $C^1$, and so $DJ(\bar{u})=0$
by Theorem \ref{th:1}. By definition, $\bar{u}$ is a \emph{weak
solution} of the linear problem
\begin{equation*}
\begin{cases}
-\Delta u = f \mbox{\quad in $\Omega$}\\
\phantom{aaa}u=0 \mbox{\quad on $\de\Omega$}.
\end{cases}
\end{equation*}
\end{example}

Clearly, the previous example is almost trivial, but in some sense it is
also the starting point of the variational theory of elliptic equations.

\section{When compactness breaks down}\label{sec:2}

We have seen that weak compactness plays an important r\^ole in
establishing the existence of the minimum point. To get a more precise
feeling, consider the following minimization problem:
\begin{equation}\label{eq:RR}
\inf_{u\in H_0^1(\Omega)\setminus \{0\}} \frac{\int_\Omega |\nabla u|^2
\, dx}{\int_\Omega |u|^2 \, dx}
\end{equation}

By the Lagrange multiplier rule, any minimizer satisfies an equation of
 the form
\begin{equation*}
\begin{cases}
-\Delta u = \lambda u &\text{ in } \Omega\\
\phantom{aaa}u=0 &\text{ on }\de\Omega
\end{cases}
\end{equation*}
for some $\lambda\in\R$.

We can take a minimizing sequence $\{u_n\}$, and it is not restrictive
to assume
\[
\int_\Omega |u_n|^2 \, dx =1 \quad \forall n\geq 1.
\]
Again, $\{u_n\}$ is bounded in $H_0^1 (\Omega)$, and so, up to a
subsequence, it converges \emph{strongly} in $L^2(\Omega)$ (see
appendix) to some $\bar{u}$. Hence
\[
\int_\Omega |u_n|^2 \, dx \to \int_\Omega |\bar u|^2\, dx
\]
and by weak lower semicontinuity,
\[
\int_\Omega |\nabla \bar u|^2 \, dx \leq \lim_{n\to\infty} \int_\Omega |\nabla
 u_n|^2\, dx.
\]
This proves that $\bar u$ is actually a minimum point of the
\emph{Rayleigh--Ritz} quotient \eqref{eq:RR}.

But the compactness of the embedding $H_0^1(\Omega) \subset L^2(\Omega)$
fails if $\Omega$ is unbounded, for instance $\Omega =\Rn$. Essentially,
 there are only three reasons why a weakly convergent sequence
 $\{u_n\}$ in an infinite--dimensional space, for example in
 $L^2(\Omega)$, does not converge strongly\footnote{It is trivially true that weakly and strongly convergent
 sequences coincide in finite dimensional
vector spaces.}
 (see also \cite{Lieb,Willem-hermann}):
\begin{description}
\item[a)] {\bf Oscillation.} If $\{u_n\}$ converges in $L^2(\Omega)$,
 then it has a subsequence converging almost everywhere in $\Omega$. But
 the sequence $u_n(x)=e^{2\pi i n x}$ does not converge in
 $L^2(-1/2,1/2)$, although it converges weakly to zero. Indeed,
 $\|u_n\|_{L^2}=1$.
\item[b)] {\bf Concentration.} Let $u$ be the characteristic function of
 $(-1/2,1/2)$, and let $u_n(x)={n}^{1/2}u(nx)$. For all $g\in C_0^\infty
 (-1/2,1/2)$, we have
\[
\int_{(-1/2,1/2)} u_n g \, dx \to 0.
\]
Since again $\|u_n\|_{L^2}=1$, this implies that $u_n \rightharpoonup 0$
 in $L^2 (-1/2,1/2)$. In the sense of distributions,
\[
u_n^2 \to \delta,
\]
where $\delta$ is the Dirac measure at zero. See \cite{Willem:minimax}
 for a study of distribution theory. 
We remark that $u_n(x)\to
 0$ for all $x$ but $x=0$.
\item[c)] {\bf Vanishing.} Let $u$ be the characteristic function of
 $(-1/2,1/2)$, and let $u_n(x)=u(x-n)$. As before, $u_n \weakto 0$. This
 time $u_n(x)\to 0$ for all $x$.
\end{description}
We shall see later that, roughly speaking, minimizing
 sequence can be classified according to their behavior like a), b) or
 c).

In the next section, we present the basis of the modern variational
 approach to possibly non--compact problems.

\section{Critical points}

As we have seen before, one way to find solutions of differential
 equations is to find global minima of suitable maps $J$ defined on
 Banach or Hilbert spaces. But global extrema are very often too much to
 be sought, and indeed their existence is far from being trivial.
 Moreover, the key point is to find points where the derivative of $J$
 vanishes, i.e. \emph{critical points of }$J$.

Let $E$ be a real Banach space with norm $\|\cdot \|$ and $J\colon
 E\to\R$
 a $C^1$ functional.

\begin{definition}
A critical point of $J$ is a point $u\in E$ such that $DJ(u)=0$ in
 $E^\star$. We say that $c\in\R$ is a critical level for $J$ if there
 exists a critical point $u$ of $J$ such that $J(u)=c$.
\end{definition}

\begin{definition}
An operator $A\colon E\to E$ is called variational if there exists a
 functional $J\in C^1 (E,\R)$ such that $A=DJ$.
\end{definition}
Hence a problem that can be translated into a functional equation like
 $A(u)=0$ is called a variational problem if $A$ is variational.

\begin{example}
Let $\Omega\subset\Rn$ be a bounded domain with smooth boundary
 $\de\Omega$, and let $E=H_0^1(\Omega)$. Let $p\colon \Omega\times
 \R\to\R$ be a continuous function such that
\begin{equation}
|p(x,s)|\leq a_1 + a_2 |s|^\ell
\end{equation}
where $\ell \leq \frac{n+2}{n-2}$ if $n >2$, $\ell$ is unrestricted if
 $n=1,2$. We assume, for simplicity, $n>2$. Let $P(x,s)=\int_0^s
 p(x,\tau)\, d\tau$. Then
\[
|P(x,s)|\leq a_3 |s| + a_4|s|^{\ell +1}.
\]
Since $\ell+1 \leq 2^\star = \frac{2n}{n-2}$, we have from the Sobolev
 embedding $E \subset L^{\ell+1}(\Omega)$ and it makes sense to define
 $\phi\colon E\to\R$ by
\[
\phi (u)=\int_\Omega P(x,u(x))\, dx.
\]
Moreover, it is easy to check that $\phi\in C^1(E,\R)$, and
\[
D\phi (u)\colon v\in E \mapsto \int_\Omega p(x,v(x))\, dx.
\]
We recall that since $\Omega$ is bounded, the embedding $E\subset
 L^{\ell+1}(\Omega)$ is compact, and hence the gradient of $\phi$ is a
 compact map.

Define $J\in C^1 (E,\R)$ by
\[
J(u)=\frac{1}{2}\int_\Omega |\nabla u|^2 \, dx - \phi (u).
\]
If $u\in E$ is a critical point of $J$, then
\[
\int_\Omega [\nabla u \cdot \nabla v - p(x,u(x))v(x)]\, dx =0 \quad
 \forall v\in H_0^1 (\Omega).
\]
Hence $u$ is a weak solution of the semilinear Dirichlet boundary value
 problem
\begin{equation*}
\begin{cases}
-\Delta u = p(x,u) &\text{ in } \Omega\\
\phantom{aaa}u=0 &\text{ on } \de\Omega.
\end{cases}
\end{equation*}
If $p$ is locally H\"{o}lder continuous, it can be proved that $u$ is in
 fact a classical $C^2$--solution.
\end{example}

A common strategy to find critical points is to look for
 \emph{constrained extremum points} over some suitable submanifold.

\begin{definition}
Let $M$ be a $C^1$ Riemannian manifold modeled on a Hilbert space $E$,
 and let $J\in C^1 (M,\R)$. A critical point of $J$ constrained on $M$
 is a point $u\in M$ such that $D_M J(u)=0$. Here $D_M$ stands for the
 derivative of $J$ on $M$.
\end{definition}

Suppose $M$ has codimension $1$ in $E$. This means (see \cite{larsen}) that there exists a
 functional $g\in C^1 (E,\R)$ such that
\[
M=\{ u\in E \colon g(u)=0\}
\]
and $Dg(u)\neq 0$ for all $u\in M$. The tangent space to $M$ at $u$ is
 given by
\[
T_u M=\{ v\in E \mid Dg(u)v=0\}
\]
and a critical point of $J$ on $M$ is a point $u\in M$ such that
 $DJ(u)v=0$ for all $v\in T_u M$. Hence $u$ satisfies
\[
DJ(u)=\lambda Dg(u)
\]
for some $\lambda\in\R$. This is the \emph{Lagrange multiplier rule}.

\begin{example}
By considering $J\colon H_0^1 (\Omega) \to\R$ defined by
\[
J(u)=\frac{1}{2}\int_\Omega |u|^2\, dx,
\]
we can look for critical points constrained on $\overline{B}=\{u\in
 H_0^1 (\Omega) \mid \int_\Omega |\nabla u|^2 \, dx \leq 1\}$. It is
 easy to show that $J_{|\overline{B}}$ has a maximum $u$, which solves
 by the above
\[
\begin{cases}
-\Delta u = \lambda u &\text{ in } \Omega\\
\phantom{aaa}u=0 &\text{ on } \de\Omega.
\end{cases}
\]
The number $\lambda$ is the \emph{first eigenvalue} of $-\Delta$ on
 $H_0^1 (\Omega)$.
\end{example}

\section{Partial compactness: Palais--Smale and mountain passes}

As a rule, to identify a minimum, one constructs a minimizing sequence,
 and tries to prove that it converges (up to a subsequence, if
 necessary). In the previuos examples, this was true because of the
 compact Sobolev embedding. When this tool is no longer available, the
 approach must be refined. Actually, relative compactness of \emph{all}
 bounded sequences is too much.

\begin{definition}[Palais--Smale sequences]
Let $E$ be a Banach space, and let $J\in C^1 (E,\R)$ be a given
 functional. We say that $J$ satisfies the Palais--Smale condition at
 level $c\in \R$, (PS)$_c$ for short, if every sequence $\{u_k\}$ in $E$
 such that $J(u_k)\to c$ and $DJ(u_k)\to 0$ in $E^\star$ is relatively
 compact in $E$.
\end{definition}

Clearly, this is a weaker considition than compactness itself, but it is
 often enough in Critical Point Theory. For the sake of completeness, we remind of the fact that several conditions like (PS) have been introduced over the past years, conditions that cover even weker situations. It is not our purpose to write a treatise on this subject, and we refer to the huge literature for details.

In 1973, Ambrosetti and Rabinowitz (\cite{AmbRab73}) proved the
 following theorem.

\begin{theorem}[Mountain pass]
Let $J$ be a $C^1$ functional on a Banach space $X$. Suppose
\begin{description}
\item[(i)] there exist a neighborhood $U$ of $0$ in $X$ and a constant
 $\varrho$ such that $J(u) \geq \varrho$ for every $u$ on the boundary
 of $U$,
\item[(ii)] $J(0) < \varrho$ and $J(v) < \varrho$ for some $v\notin U$.
\end{description}
Set
\[
c=\inf_{P\in \mathcal{P}} \max_{w\in P} J(w) \geq \varrho,
\] 
where $\mathcal{P}$ denotes the class of continuous paths joining $0$ to
 $v$.
The conclusion is: there is a sequence $\{u_j\}$ in $X$ such that $J(u_j)\to c$
 and $DJ(u_j)\to 0$ in $X^\star$.
If $J$ satisfies (PS)$_c$, then $c$ is a critical level for $J$.
\end{theorem}

This result is very useful to find critical points of $J$ even when $J$
 is strongly undefined, i.e. when $\inf_E J = -\infty$ and $\sup_E J = +
 \infty$. The geometry of $J$ is that $0$ is local minimum, but not a
 global one. See the figure below, showing a profile with two maxima, one local minimum, and a path that a clever walker should follow to minimize effort.

\begin{figure}
\centering
\includegraphics[height=34mm]{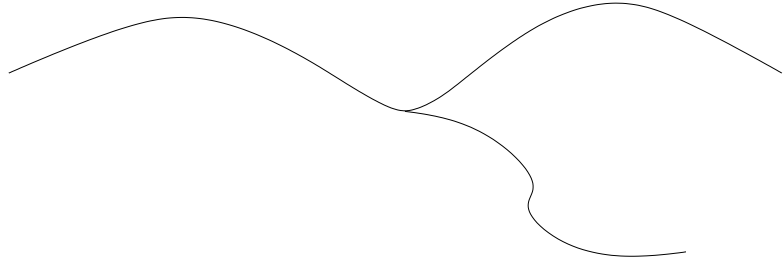}
\caption{The mountain--pass geometry}
\end{figure}
 
Needless to say, the difficult task is to prove that $J$
 satisfies the Palais--Smale condition at the ``right'' level. Good
 functionals satisfy it all levels, but many others satisfy it only in
 subregions of $\R$.

\section{Critical exponents}

Apart from unboundedness of the domain, another source of
 non--compactness is the presence of the \emph{critical Sobolev
 exponent} (see the Appendix). The embedding $H_0^1 (\Omega) \subset
 L^{2^\star}(\Omega)$ fails to be compact even if $\Omega$ is
bounded.\footnote{See the appendix.}  Hence differential equations with 
a critical behavior appear to be very  interesting. Indeed, they can
have only trivial solutions, as the  following example shows.

\begin{example}
Consider
\begin{equation}\label{eq:nonex}
\begin{cases}
-\Delta u = |u|^{q-1} u   &\text{ in } \Omega\\
\phantom{aaa}u=0 &\text{ on } \de\Omega,
\end{cases}
\end{equation}
where $\Omega$ is a bounded domain in $\Rn$. For $x\in\de\Omega$, we
 denote by $\nu_x$ the outward normal vector to $\Omega$.

\begin{proposition}
If $x \cdot \nu_x \geq 0$ on $\de\Omega$, then \eqref{eq:nonex} has only
 the trivial solution $u\equiv 0$, whenever $q \geq \frac{n+2}{n-2}$.
\end{proposition}

The proof is an immediate consequence of the next result, due to
 Pohozaev.

\begin{lemma}
If $u$ is a smooth solution of
\begin{equation}
\begin{cases}
-\Delta u = p(x,u)   &\text{ in } \Omega\\
\phantom{aaa}u=0 &\text{ on } \de\Omega,
\end{cases}
\end{equation}
then
\[
n\int_\Omega P(u(x))\, dx + \frac{2-n}{2} \int_\Omega u(x)p(u(x))\, dx =
 \frac{1}{2} \int_{\de\Omega} (x\cdot \nu_x)\left| \frac{\de
 u}{\de\nu}\right|^2 \, d\sigma.
\]
\end{lemma}
\end{example}

Probably the first, but surely the most impressive paper to study
exstensively the loss of compactness in differential equations is that
 by
Brezis and Nirenberg (\cite{BN83}). In a previous section, we showed
 that equations on unbounded domains usually present additional 
difficulties, due to the loss of compact embedding in $L^p$ spaces. 
A different phenomenon leading to ``non--compact'' problems is that of 
\emph{limiting Sobolev exponent}, for which the embedding is not
 compact,
even on bounded domains.

In their paper, Brezis and Nirenberg consider a bounded domain 
$\Omega\subset \Rn$ and the critical elliptic equation
\begin{equation}\label{eq:BN}
\begin{cases}
-\varDelta u = u^p + \lambda u &\text{ in }\Omega\\
u>0 &\text{ in }\Omega\\
u=0 &\text{ on }\partial\Omega,
\end{cases}
\end{equation}
where $p=(n+2)/(n-2)$, $n\geq 3$, and $\lambda$ is a real constant. We
 remark explicitely that the same problem with subcritical exponent is
 completely solvable:

\begin{theorem}
Assume that $\Omega$ is a bounded domain, and $2 < p < 2^\star$. Then
 the problem
\begin{equation}
\begin{cases}
-\varDelta u = u^p + \lambda u &\text{ in }\Omega\\
u>0 &\text{ in }\Omega\\
u=0 &\text{ on }\partial\Omega,
\end{cases}
\end{equation}
has a nontrivial solution if and only if $\lambda > - \lambda_1
 (\Omega)$, the
 first eigenvalue of $-\Delta$ in $H_0^1 (\Omega)$.
\end{theorem}

\begin{proof}
See \cite{Willem:minimax}
\qed \end{proof}

Their main result is

\begin{theorem}
Assume $n\geq 4$. Then for every $\lambda \in (0,\lambda_1)$, there
 exists a solution of \eqref{eq:BN}. Here $\lambda_1$ is the first
 eigenvalue of $-\varDelta$ with Dirichlet boundary conditions on
 $\Omega$.
\end{theorem}

Since their argument has opened the way to hundreds of papers on
 non-compact problems, we sketch the main ideas. First of all, they
 choose a variationa approach, and set
\[
S_\lambda = \inf_{\substack{u\in H_0^1 (\Omega)\\\|u\|_{L^p}=1}} 
 \int_\Omega \left( |\nabla u|^2 - \lambda |u|^2 \right)\, dx,
\]
so that $S_0$ corresponds to the best constant in the Sobolev embedding
 theorem $H_0^1(\Omega) \subset L^{p+1}(\Omega)$ When $\Omega=\Rn$, it
 was proved by Aubin and Talenti that $S_0$ is attained by the functions
\[
U_\ge (x)=C_\ge (\ge + |x|^2)^{-(n-2)/2},
\]
where $C_\ge >0$ is a normalization constant.

The second step is to show a \emph{a priori} existence result:

\begin{lemma}\label{lem:BN1}
If $S_\lambda < S$, then $S_\lambda$ is attained.
\end{lemma}
\begin{proof}
Let $\{u_j\}$ be a minimizing sequence in $H_0^1$, i.e.
\begin{equation}
\int_\Omega |u_j|^{p+1}\, dx =1,
\end{equation}
\begin{equation}
\int_\Omega |\nabla u_j|^2 \, dx - \lambda \int_\Omega |u_j|^2 \, dx =
 S_\lambda + o(1) \quad \mbox{as $j\to\infty$.}
\end{equation}
Up to a subsequence,
\[
u_j \rightharpoonup u \mbox{\quad weakly in } H_0^1(\Omega),
\]
\[
u_j\to u \mbox{\quad strongly in } L^2(\Omega),
\]
\[
u_j \to u \mbox{\quad almost everywhere in } \Omega,
\] 
with $\|u\|_{L^{p+1}}\leq 1$. Set $v_j = u_j - u$. Hence we have
 $\|\nabla u_j\|_{L^2} \geq S_0$. It follows that $\lambda \|u\|_{L^2}^2
 \geq S_0 - S_\lambda > 0$, and so $u\not\equiv 0$. By construction,
\begin{equation}
\|\nabla u\|_{L^2}^2 +\|\nabla v_j\|_{L^2}^2 -\lambda \|u\|_{L^2}^2 = S_
\lambda + o(1)
\end{equation}
since $v_j \rightharpoonup 0$ weakly. By the Brezis--Lieb Lemma (see the
 next section),
\[
\|u+v_j\|_{L^{p+1}}^{p+1}=\|u\|_{L^{p+1}}^{p+1} +
 \|v_j\|_{L^{p+1}}^{p+1}+o(1).
\]
Thus we have
\[
1=\|u\|_{L^{p+1}}^{p+1}+\|v_j\|_{L^{p+1}}^{p+1}+o(1)
\]
and therefore
\[
1\leq \|u\|_{L^{p+1}}^{2}+\|v_j\|_{L^{p+1}}^{2}+o(1)
\]
which leads to
\begin{equation}
1\leq \|u\|_{L^{p+1}}^{2}+\frac{1}{S_0}\|\nabla
 v_j\|_{L^{p+1}}^{2}+o(1).
\end{equation}
The following claim can be proved by distinguishing the cases (i)
 $S_\lambda > 0$ and (ii) $S_\lambda \leq 0$:
\[
\|\nabla u\|_{L^2}^2 -\lambda \|u\|_{L^2}^2 \leq S_\lambda
 \|u\|_{L^{p+1}}^2.
\]
Now the proof is complete, since $u$ is not identically zero.
\qed \end{proof}

Hence we should compare $S_\lambda$ and $S_0$, and this is done by some
 \emph{cut-off} procedure:

\begin{lemma}\label{lem:BN2}
We have
\[
S_\lambda < S_0 \qquad \mbox{for all $\lambda >0$}.
\]
\end{lemma}

\begin{proof}
The task is to estimate the ratio
\[
Q_\lambda (u)=\frac{\int_\Omega |\nabla u|^2 \, dx - \lambda \int_\Omega
 |u|^2 \, dx}{\left( \int_\Omega |u|^{p+1}\, dx\right)^{2/(p+1)}}.
\]
This can be done with the functions
\[
u(x)=u_\ge (x) = \frac{\varphi (x)}{(\ge + |x|)^{(n-2)/2}},
\]
where $\varphi$ is $C_0^\infty$ function concentrated near the origin of
 $\Rn$. After some calculations, it is possible to prove that
 $Q_\lambda (u_\ge) < 0$ if $\ge$ is small enough.
\qed \end{proof}

Now, a solution of \eqref{eq:BN} is found by rescaling any minimizer in
 Lemma \ref{lem:BN1}. By the same token, the following result can be
 established.

\begin{lemma}
Let
\[
S_0 = \inf \{ \|u\|_2^2 \mid (u\in H_0^1 (\Omega))\land
 (\|u\|_{2^\star}^{2^\star}=1)\}
\]
be the best Sobolev constant. Then for any
\[
c < \frac{1}{n}S_0^{n/2}
\]
the functional
\[
J(u)=\frac{1}{2}\int_\Omega (|\nabla u|^2 -\lambda |u|^2)\, dx -
 \frac{1}{2^\star}\int_\Omega |u|^{2^\star}\, dx
\]
satisfies the (PS)$_c$ condition.
\end{lemma}

More general situations can be studied by means of the
 \emph{mountain--pass} theorem, for example
\begin{equation}
\begin{cases}
-\varDelta u = u^p + f(x,u) &\text{ in } \Omega\\
u>0 &\text{ in } \Omega\\
u=0 &\text {on }\Omega,
\end{cases}
\end{equation}
where the main assumption of $f\colon \Omega\times [0,\infty)\to\R$ is
 that
\[
\lim_{u\to+\infty} \frac{f(x,u)}{u^p}=0.
\]

The case $n=3$ is more delicate, since the estimates in Lemma
 \ref{lem:BN2} are no longer true. Anyway, the following result is
 given.

\begin{theorem}
Assume $\Omega$ is a ball. Then \eqref{eq:BN} has a solution if and only
 if $\frac{1}{4}\lambda_1 < \lambda < \lambda_1$.
\end{theorem}

\subsection{The Brezis--Lieb lemma}

In this appendix we present a very useful result about convergence in
 $L^p$ spaces. The situation is this: a bounded sequence $\{u_n\}$ in
 $L^p$ tends pointwise to some $u\in L^p$. If $\|u_n\|_p\to \|u\|_p$,
 and $p > 1$, it is a consequence of the uniform convexity that $u_n \to
 u$ strongly. Is it possible to precise the behavior of $\|u_n-u\|_p$?

\begin{lemma}[\cite{BrezisLieb83}]
Let $\Omega$ be an open subset of $\Rn$, and let $\{u_n\}$ be a bounded
 sequence of $L^p(\Omega)$ functions tending a.e. to $u\in L^p$. Then
\[
\|u\|_p=\lim_{n\to\infty} \left( \|u_n\|_p -\|u_n-u\|_p \right).
\]
\end{lemma}

\begin{proof}
Let $M=\sup_n \|u_n\|_p<\infty$. By the limit
\[
\lim_{|s|\to\infty} \frac{|s+1|^p-|s|^p-1}{|s|^p}=0,
\]
we see that for every $\ge >0$ there is a constant $C_\ge$ such that
\[
\left| |s+1|^p -|s|^p -1 \right| \leq \ge |s|^p + C_\ge \quad \forall
 s\in\R.
\]
Hence, for all $a$,$b\in\R$,
\[
\left| |a+b|^p -|a|^p -|b|^p \right| \leq \ge |a|^p + C_\ge |b|^p.
\]
Let now $f_n = ||f_n|^p - |f_n-f|^p - |f|^p|$ and $Z_n = \max \{0,f_n -
 \ge |u_n-u|^p\}$. We know that $Z_n$ tends to zero a.e., and $0\leq Z_n
 \leq C_\ge |f|^p$. By the Dominated Convergence theorem, $Z_n \to 0$ in $L^1$. But
 $0\leq f_n \leq \ge |u_n-u|^p + Z_n$, and so
\[
\|f_n\|_1 \leq \ge \|u_n-u\|_1 + \|Z_n\|_1\leq 2^p M^p + \|Z_n\|_1.
\]
Finally, we deduce that $\|f_n\|_1\to 0$.
\qed \end{proof}

\section{The Concentration--Compactness principle}

As seen before, there are some evident reasons why a bounded sequence
 in $L^2(\Rn)$ cannot converge strongly. Roughly speaking, we can say
 $\Rn$ is very large, and sequences can ``move off'' to infinity, or
 concentrate at finite points even though they are constrained to minimize
 some functional. Hence it would be desirable to have some tool to work
 out these features, to regain compactness. In this section we briefly
 survey the celebrated Concentration--Compactness principle by
 P.~L.~Lions (\cite{Lions85-1,Lions85-2}). A completely different
 strategy will be one of the most important topics of our thesis in the
 next part.

\begin{theorem}[Locally compact case]
Let $\{\varrho_n\}$ be a sequence in $L^1 (\Rn)$ such that
\begin{description}
\item[(1)] $\varrho_n \geq 0$ almost everywhere
\item[(2)] $\int_\Rn \varrho_n (x)\, dx = \lambda > 0$.
\end{description}
Then $\{\varrho_n\}$ has a subsequence $\{\varrho_{n_k}\}$ satisfying one of
 the following properties:
\begin{itemize}
\item {\bf Compactness.} There is a sequence $\{y_k\}$ in $\Rn$ such
 that 
\[
(\forall \ge >0)(\exists R>0)\colon \int_{y_k + B_R} \varrho_{n_k}(x)\, dx
 \geq \lambda - \ge.
\]
\item {\bf Vanishing.}
\[
\lim_{k\to\infty} \sup_{y\in\Rn}\int_{y+B_R}\varrho_{n_k}(x)\, dx =0 \quad
 \forall R>0.
\]
\item {\bf Dichotomy.} There exists $\alpha \in (0,\lambda)$ such that
 for all $\ge >0$ there are an integer $k_0 \geq 1$ and
 $\rho_k^1$,$\varrho_k^2\in L^1 (\Rn)$ such that
\begin{enumerate}
\item $\varrho_k^1$, $\varrho_k^2 \geq 0$
\item $\|\varrho_{n_k}-(\varrho_k^1+\rho_k^2)\|_{L^1} \leq \ge$
\item $\left| \int_\Rn \varrho_k^1 (x)\, dx - \alpha \right| \leq \ge$
\item $\left| \int_\Rn \varrho_k^2 (x)\, dx - \alpha \right| \leq \ge$
\item
 $\operatorname{dist}(\operatorname{supp}\varrho_k^1,\operatorname{supp}\varrho_k^2)\to + \infty$ as $k\to + \infty$.
\end{enumerate}
\end{itemize}
\end{theorem}

The following result is often useful to exclude the vanishing case.

\begin{lemma}[\cite{Lions85-1}]
Let $1<p\leq\infty$, $1\leq q < \infty$, and $N\geq 1$ an integer. When
 $N>p$, we suppose that $q\neq p^\star$. Let $\{u_k\}$ be a bounded
 sequence in $L^q (\Rn)$ such that $\{ |\nabla u_k|\}$ is bounded in
 $L^p (\Rn)$; if there is $R>0$ such that
\[
\lim_{k\to\infty} \sup_{y\in\Rn} \int_{B(y,R)} |u_k(x)|\, dx = 0,
\]
then $u_k\to 0$ in $L^r(\Rn)$ provided that $\min\{q,p^\star\} < r <
 \max\{q,p^\star\}$.
\end{lemma}

To present the result about the non--locally compact case, i.e. when
 some critical exponent is involved, we need some preliminaries.

\begin{definition}
$\mathcal{M}(\Omega)$ denotes the space of bounded Radon measures on the
 open set $\Omega\subset\Rn$, endowed with the norm
\[
\|\mu\|=\sup_{\substack{u\in C_0 (\Omega)\\ \|u\|_\infty = 1}} \left|
 \mu (u)\right| .
\]
we say that a sequence $\{\mu_k\}$ in $\mathcal{M}(\Omega)$ converges
 weakly\footnote{The reader may remark that this convergence should be
rather called \textit{weak}${}^*$. Anyway, such a simplification is
customary in Probability theory.} to $\mu$ if \[ \lim_{k\to\infty} 
\mu_k (u) = \mu (u) \quad \forall u\in C_0 (\Omega). \]
and $\{\mu_k\}$ converges tightly to $\mu$ if
\[
\lim_{k\to\infty} \mu_k (u) = \mu (u) \quad \forall u\in L^\infty
(\Omega). \]
\end{definition}

\begin{theorem}[Critical case]
Let $m \geq 1$ be an integer, $p\in [1,n/m)$, and denote
\[
p^{*m}=\frac{np}{n-mp}
\]
\[
\theta = 1-\frac{mp}{n}.
\]
Let
\[
S(m,p,n)=\inf_{\substack{u\in D^{m,p}(\Rn)\\ \|u\|_{p^{*m}}=1}} \|D^m u \|_p^p.
\]
Assume the sequence $\{u_k\}$ in $D^{m,p}(\Rn)$ converges weakly to $u\in D^{m,p}(\Rn)$, and that there exist two bounded measures $\lambda$ and $\mu$ such that
\[
|D^m u_k| \weakto \mu \mbox{\quad in }\mathcal{M}(\Rn)
\]
\[
|u_k|^{p^{*m}}\weakto \lambda \mbox{\quad tightly.}
\]
Then
\begin{enumerate}
\item There exist a set $J$, at most countable, some points 
$\{x_j\}_{j\in J}$, and some real numbers $\{a_j\}_{j\in J}$ such that
$a_j >0$ for all $j\in J$ and \[ \lambda = |u|^{p^{*m}}+\sum_{j\in J}
a_j \delta_{x_j}. \]
\item There exist real numbers $b_j > 0$ such that $S(m,p,n)a_j^\theta \leq b_j$ and
\[
|D^m u|^p + \sum_{j\in J} b_j \delta_{x_j} \leq \mu,
\]
\[
\sum_{j\in J} a_j^\theta < \infty.
\]
\item If $v\in D^{m,p}(\Rn)$ and $|D^m (u_k +v)|^p \weakto \tilde{\mu}$ in $\mathcal{M}(\Rn)$, then $\tilde{\mu}-\mu\in L^1 (\Rn)$ and
\[
|D^m (u+v)|^p + \sum_{j\in J} b_j \delta_{x_j} \leq \tilde{\mu}.
\]
\item If $u\equiv 0$ and $0 < \mu(\Rn) < S(m,p,n)\lambda (\Rn)$, then $J$ is a singleton, and there exist $C>0$, $x_0\in\Rn$ such that
\[
\lambda = C \delta_{x_0},\qquad \mu=S(m,p,n) C^{1+\theta} \delta_{x_0}.
\]
\end{enumerate}
\end{theorem}

Finally, we present a recent improvement appearing in
\cite{BN-T-Willem}, which unifies the concentration at infinity and the
concentration at finite points. Despite its simple appearence, it can be
adapted to a lot of different situations.

\begin{theorem}
Let $\Omega\subset\Rn$ an open set, and let $1\leq p<\infty$. Assume
that
\begin{description}
\item[a)] $\{u_k\}$ is a bounded sequence in $L^p(\Omega)$;
\item[b)] $u_k\to u$ a.e. in $\Omega$,
\item[c)] $|u_k-u|^p \weakto \nu$ in $\mathcal{M}(\Omega)$,
\end{description}
and define
\[
\nu_\infty = \lim_{R\to\infty}\limsup_{k\to\infty} \int_{|x| > R}
|u_k|^p\, dx.
\]
Then it follows that
\[
\limsup_{k\to\infty} |u_k|^p_p = |u_k|_p^p + \|\nu\| + \nu_\infty.
\]
\end{theorem}

\begin{proof}
1) We write $v_k=u_k-u$. According to the Brezis--Lieb lemma, we have,
for every non-negative $h\in C_0^\infty (\Omega)$,
\[
\int h |u|^p \, dx = \lim_{k\to\infty} \left( \int h |u_k|^p \, dx -
\int h|v_k|^p \, dx \right).
\]
It follows that
\begin{equation}
|u_k|^p \weakto |u|^p + \nu \mbox{\quad in } \mathcal{M}(\Omega).
\end{equation}

2) For $R>1$, let $\Psi_R\in C_0^\infty(\Rn)$ be such that $\Psi_R
(x)=1$ for $|x| > R+1$, $\Psi_R(x) = 0$ for $|x|<R$, and $0\leq \Psi_R
\leq 1$ on $\Rn$. It is easy to check that
\[
\nu_\infty = \lim_{R\to\infty} \limsup_{k\to\infty} \int \Psi_R
|u_k|^p\, dx.
\]

3) For every $R>1$, we have by (1.14)
\begin{multline*}
\limsup_{k\to\infty} |u_k|^p\, dx = \limsup_{k\to\infty} \left( \int
\Psi_R^p |u_k|^p \, dx + \int (1-\Psi_R^p)|u_k|^p\, dx\right) \\
=\limsup_{k\to\infty} \int \Psi_R^p |u_k|^p\, dx + \int
(1-\Psi_R^p)|u|^p\, dx + \int (1-\Psi_R^p) \, d\nu (x).
\end{multline*}
When $R\to\infty$, we obtain
\begin{align*}
\limsup_{k\to\infty} \int |u_k|^p\, dx &= \nu_\infty + \int |u|^p\, dx
+ \int \, d\nu (x)\\
&= \nu_\infty +|u|_p^p + \|\nu\|.
\end{align*}
\qed\end{proof}

\begin{remark}
Clearly, the terms $\|\nu\|$ and $\nu_\infty$ are responsible for lack of compactness. The r\^{o}le of the latter is evident, while the first measures the concentration of the sequence $\{u_{k}\}$ at finite points.
\end{remark}

\section{A Sobolev-Hardy type inequality}

In this section we use Theorem 1.19 to find the best constant of a Hardy--type inequality. We follow \cite{SeSmWi}.
Set $\R^N=\R^k\times\R^{N-k}$, with $2\leq k\leq N$, and write $x=(x',z)\in \R^k\times\R^{N-k}$. For given numbers $q$, $s$ such that $1<q<N$, $0\leq s\leq q$, and $s<k$, set $q_\star (s,N,q)=\frac{q(N-s)}{N-q}$. We also write $q_\star (s)$ instead of $q_\star (s,N,q)$.

We are able to proof the following statement, which extends a previous result contained in [13].

\begin{theorem}
Let $2\leq k \leq N$, $1\leq p<k$, $\alpha + k > 0$, $u\in C_0^\infty (\Rn)$. Then
\begin{equation}\label{eq:sob-har}
\int_\Rn |u|^p |x'|^\alpha \, dx \leq \frac{p^p}{(\alpha + k)^p} \int_\Rn |\nabla u|^p |x'|^{\alpha+p}\, dx.
\end{equation}
Moreover, the constant $\frac{p^p}{(\alpha + k)^p}$ is optimal.
\end{theorem}

It is convenient to state the following Lemma.

\begin{lemma}
If $1<p<\infty$ and $\alpha + N >0$, then
\begin{align*}
\int_\Rn |u|^p |x|^\alpha \, dx &\leq \frac{p^p}{(\alpha +N)^p} \int_\Rn
|x\cdot \nabla u|^p |x|^\alpha \, dx \\
&\leq \frac{p^p}{(\alpha+N)^p}
\int_\Rn |\nabla u|^p |x|^{\alpha +N}\, dx 
\end{align*}
for all $u\in C_0^\infty
(\Rn)$. Moreover, the constant $\frac{p^p}{(\alpha+N)^p}$ is the best
possible. \end{lemma}

\begin{proof}
We approximate $x\mapsto |x|^\alpha$ by the ``sequence'' $x\mapsto (|x|^2+\ge)^{\alpha /2}$, so that
\[
\operatorname{div}[|x|^\alpha x] = (\alpha +N) |x|^\alpha,
\]
at least in a weak sense. Moreover, for all $u\in C_0^\infty (\Rn)$,
\[
\int_\Rn \operatorname{div}[|u|^p |x|^\alpha x]\, dx =0.
\]
We get from the H\"{o}lder inequality that
\begin{multline*}
\int_\Rn |u|^p |x|^\alpha \, dx \leq \frac{p^p}{(\alpha+N)^p}\int_\Rn
|x\cdot \nabla u| |u|^{p-1} |x|^\alpha \, dx \\ 
\leq
\frac{p^p}{(\alpha+N)^p} \left( \int_\Rn |u|^{(p-1)p'} |x|^\alpha \, dx
\right)^{1/p'} \left( \int_\Rn |x\cdot \nabla u|^p |x|^\alpha \, dx
\right)^{1/p}\\ 
= \frac{p^p}{(\alpha+N)^p} \left( \int_\Rn |u|^p
|x|^\alpha \, dx \right)^{1-1/p} \left( \int_\Rn |x\cdot \nabla u|^p
|x|^\alpha \, dx \right)^{1/p}.
\end{multline*} 
For the optimality of the
constant, we refer to \cite{SeSmWi}. \qed\end{proof}

\begin{proof}[of Theorem 10]
Consider first the case $p=2$. Then
\begin{align*}
\int_{\R^{N-k}\times \R^k} |\nabla u|^2 |x'|^{\alpha + 2}\, dx &\geq \int_{\R^{N-k}}dz \int_{\R^k} |\nabla_{x'}u|^2 |x'|^{\alpha +2}\, dx' \\
&\geq \frac{(\alpha +2)^2}{4} \int_{\R^{N-k}}dz \int_{\R^k} |u|^2 |x'|^\alpha \, dx' \\
&= \frac{(\alpha +2)^2}{4} \int_{\Rn}|u|^2 |x'|^\alpha\, dx
\end{align*}
where we have used the previous Lemma. Choose now
\[
u\colon (x',z)\mapsto v(x')w(z),
\]
where $v\in C_0^\infty (\R^k)$ and $w\in C_0^\infty (\R^{N-k})$. In the rest of the proof, we identify the gradients of $v$ and $w$ with two vectors in $\Rn$:
\[
\begin{cases}
\nabla v = \left( \frac{\de v}{\de x'},0 \right)\\
\nabla w=\left( 0,\frac{\de w}{\de z}\right).
\end{cases}
\]
Hence, $\nabla v \perp \nabla w$. By definition,
\begin{multline*}
\frac{\int_\Rn |\nabla u|^2 |x'|^{\alpha +2}\, dx}{\int_\Rn |u|^2 |x'|^\alpha \, dx} = \frac{\int_\Rn |x'|^{\alpha +2} \left( |\nabla v|^2 |w|^2 + |v|^2 |\nabla w|^2 \right)\, dx}{\int_\Rn |x'|^\alpha |v|^2 |w|^2 \, dx}\\
=
\frac{\int_{\R^k} |x'|^{\alpha+2} |\nabla v|^2 \, dx'}{\int_{\R^k} |v|^2 |x'|^\alpha \, dx'}+\frac{\int_{\R^k} |x'|^{\alpha +2} |v|^2 \, dx'}{\int_{\R^k} |x'|^\alpha |v|^2 \, dx'}\frac{\int_{\R^{N-k}}|\nabla w|^2}{\int_{\R^{N-k}}|w|^2}.
\end{multline*}

We deduce
\begin{multline*}
\inf_{u\in C_0^\infty (\Rn)\setminus \{0\}} \frac{\int_\Rn |\nabla u|^2 |x'|^{\alpha +2}\, dx}{\int_\Rn |u|^2 |x'|^\alpha \, dx} \\ \leq \inf_{\substack{v\in C_0^\infty (\R^k)\setminus \{0\} \\ w\in C_0^\infty (\R^{N-k}\setminus \{0\}}} \frac{\int_\Rn |x'|^{\alpha+2} \left( |w|^2 |\nabla v|^2 + |v|^2 |\nabla w|^2 \right) \, dx}{\int_\Rn |x'|^2 |v|^2 |w|^2 \, dx}\\
\leq \inf_{v\in C_0^\infty (\R^k)\setminus \{0\}} \frac{\int_{\R^k} |x'|^{\alpha+2}|\nabla v|^2\, dx}{\int_{\R^k} |x'|^\alpha |v|^2 \, dx}=\frac{(\alpha +2)^2}{4}
\end{multline*}
because
\[
\inf_{w\in C_0^\infty (\R^{N-k})\setminus \{0\}} \frac{\int_{\R^k} |\nabla w|^2\, dx}{\int_{\R^{N-k}} |w|^2 \, dx}=0
\]
and the optimality result contained in the previous Lemma.

The general case $p\neq 2$ is slightly more involved, since we work in a non--hilbertian setting. The validity of the inequality is checked as before. Anyway, the optimality of $p^p / (\alpha +k)^p$ requires the following modifications. Split again $u(x)=v(x')w(z)$, so that
\[
\frac{\int_\Rn |\nabla u|^p |x'|^{\alpha +p}\, dx}{\int_\Rn |u|^p |x'|^\alpha \, dx}
=\frac{\int_\Rn |x'|^{\alpha +p} \left( |\nabla v|^2 |w|^2 + |\nabla w|^2 |v|^2  \right)^{p/2}\, dx}{\int_\Rn |x'|^\alpha |v|^p |w|^p \, dx}.
\]
We introduce the map
\begin{align*}
f \colon \R\times\R &\to \R \\
(s,t) &\mapsto \left( s^2 +t^2 \right)^{p/2}.
\end{align*}
Clearly $f$ is a convex function. For each $\ge \in (0,1)$, we see that
\begin{multline*}
f(s,t)=f \left( \ge (\frac{s}{\ge} ,0)+(1-\ge)(0,\frac{t}{1-\ge}) \right) \leq
\\
\ge f \left( \frac{s}{\ge},0 \right) + (1-\ge) f \left( 0,\frac{t}{1-\ge}\right) = \ge \frac{s^p}{\ge^p} + (1-\ge) \frac{t^p}{(1-\ge)^p}.
\end{multline*}

Hence
\begin{multline*}
\frac{\int_\Rn |\nabla u|^p |x'|^{\alpha+p}\, dx}{\int_\Rn |u|^p |x'|^\alpha \, dx}\leq \ge^{1-p} \frac{\int_{\R^k} |x'|^{\alpha +p} |\nabla v|^p \, dx'}{\int_{\R^k} |x'|^\alpha |v|^p\, dx'}+ (1-\ge)^{1-p} \frac{\int_{\R^{N-k}} |\nabla w|^p \, dz}{\int_{\R^{N-k}} |w|^p\, dz}.
\end{multline*}

We deduce that
\begin{equation*}
\inf_{u\in C_0^\infty (\Rn)\setminus \{0\}} \frac{\int_\Rn |\nabla u|^p |x'|^{\alpha +p} \, dx}{\int_\Rn |x'|^\alpha |u|^p \, dx} \leq \ge^{1-p} \inf_{v\in C_0^\infty (\R^k)\setminus \{0\}} \frac{\int_{\R^k} |x'|^{\alpha +p} |\nabla v|^p \, dx'}{\int_{\R^k} |x'|^\alpha |v|^p \, dx'}.
\end{equation*}
We conclude by letting $\ge\to 0$.
\qed
\end{proof}

It is a natural question if the best constant $C$ in \eqref{eq:sob-har} is attained. More precisely, let
\begin{equation}\label{eq:S}
S=\inf \left\lbrace \int_{\R^N} |\nabla u|^q\, dx \colon u\in D^{1,q}(\R^N) \quad\mbox{and}\quad \int_{\R^N} \frac{|u(x',z)|^{q_\star (s)}}{|x'|^s} \, dx' dz =1 \right\rbrace .
\end{equation}
Here $S$ depends on all the constants of the problem, i.e. $k$, $q$, $s$ and $N$. For $s=0$, it reduces to the best Sobolev constant, and to the best Hardy constant for $s=q$ and $k=N$. It is well-known that the best Hardy constant is never attained. However, the following theorem states that in the whole family of inequalities between the Sobolev and the Hardy's one, the latter is the only negative case.

\begin{theorem}
Assume $1\leq k \leq N$, $1<q<N$, $x=(x',z)$, $0\leq s < q$ and $s<k$. Then the extremal problem \eqref{eq:S} attains its infimum at a function $u\in D^{1,q}(\R^N)$ which satisfies
\[
\int_{\R^N} |\nabla u|^q\, dx =S, \quad \int_{\R^N} \frac{|u(x',z)|^{q_\star (s)}}{|x'|^s} \, dx' dz =1.
\]
\end{theorem}

\begin{proof}
For a standard proof using Theorems 7 and 8, we refer to \cite{BadTar}. Here we give a rather short proof by means of a suitable modification of Theorem 9.

Let $\{u_n\}$ be a minimizing sequence for $S$. Hence $\{u_n\}$ is bounded in $D^{1,q}(\R^N)$. By Remark 7 in \cite{BadTar}, it is known that we can always find sequences $\lambda_n >0$ and $\zeta_n\in\R^{N-k}$ such that the new sequence $\{v_n\}$ defined by
\[
v_n(x',z)=\lambda_n^{\frac{N-q}{q}} u_n(\lambda_n x',\lambda_n (z-\zeta_n))
\]
is still a minimizing sequence for $S$, and moreover
\[
\sup_{\zeta\in\R^{N-k}} \int_{|z-\zeta| < 1} \int_{|x'|<1} \frac{|v_n|^{q_\star (s)}}{|x'|^s}\, dx = \int_{|z| < 1} \int_{|x'|<1} \frac{|v_n|^{q_\star (s)}}{|x'|^s}\, dx =\frac{1}{2}.
\]
Hence, without loss of generality, we can replace $u_n$ by $v_n$. By the Sobolev embedding theorem, we may also assume, up to a subsequence, that
$$v_n \weakto v \mbox{\quad in } D^{1,q}(\R^N)$$
$$|\nabla (v_n-v)|^q \to \mu \mbox{\quad in } \mathcal{M}(\R^N)$$
$$\left| \frac{v_n-v}{|x'|^{s/q_\star (s)}}\right|^{q_\star (s)}\to \nu \mbox{\quad in }\mathcal{M}(\R^N)$$
$$v_n \to v \mbox{\quad a.e. in } \R^N$$
We introduce now the two quantities
\[
\mu_\infty = \lim_{R\to+\infty} \limsup_{n\to\infty} \int_{|x|\geq R} |\nabla v_n|^q \, dx
\]
\[
\nu_\infty = \lim_{R\to+\infty} \limsup_{n\to\infty} \int_{|x|\geq R} \frac{|v_n|^{q_\star (s)}}{|x'|^s}\, dx.
\]
 From the definition of $S$ we deduce
\begin{align*}
\nu_\infty^{q/q_\star (s)} &\leq S^{-1} \mu_\infty\\
\|\nu\|^{q/q_\star (s)} &\leq S^{-1} \|\mu\|\\
|\nabla v|_q^q + \|\mu\| + \mu_\infty & \leq S \\
\left| \frac{v}{|x'|^{s/q_\star (s)}}\right|_{q_\star (s)}^{q_\star (s)} + \|\nu\| + \nu_\infty &= 1.
\end{align*}

Hence,
\begin{align*}
S \geq  |\nabla v|_q^q + \|\mu\| + \mu_\infty &\geq |\nabla v|_q^q + S \|\nu\|^{q/q_\star (s)} + S \nu_\infty^{q/q_\star (s)} \\
&\geq S \left| \frac{v}{|x'|^{s/q_\star (s)}}\right|_{q_\star (s)}^{q_\star (s)} + S \|\nu\|^{q/q_\star (s)} + S \nu_\infty^{q/q_\star (s)},
\end{align*}
and finally
\[
\left| \frac{v}{|x'|^{s/q_\star (s)}}\right|_{q_\star (s)}^{q_\star (s)} + \|\nu\|^{q/q_\star (s)} + \nu_\infty^{q/q_\star (s)} \leq 1.
\]
This implies that $\|\nu\|$, $\nu_\infty$ and $\left| \frac{v}{|x'|^{s/q_\star (s)}}\right|_{q_\star (s)}^{q_\star (s)}$ are all equal to either 0 or 1. But $\nu_\infty =0$, since it follows from the choice of $\lambda_n$ and $\zeta_n$ that $\nu_\infty \leq 1/2$. Suppose that $v=0$. Then
\[
1=\|\nu\|^{q/q_\star (s)} = S^{-1} \|\mu\|.
\]
Hence, for all open subsets $\omega\subset\Rn$, $\nu (\omega) \leq \nu(\Rn)$, and this immediately implies that $\nu$ and $\mu$ are concentrated on a single point $x_0=(x_0',z_0)\in\R^{k}\times\R^{N-k}$. Once again, by the choice of $\lambda_n$ and $\zeta_n$, there is a number $r > 0$ such that $|x_0| = \max \{|x_0'|,|z_0|\} \geq r$. We then reach the contradiction
\[
0= \lim_{n\to\infty} \int_{|x-x_0| < r} \frac{|v_n|^{q_\star (s)}}{|x'|^{s}} \, dx = \|\nu\| = 1.
\]
The only possibility is that $\|\nu\|=0$, so that
\[
\frac{|v|^{q_\star (s)}}{|x'|^{s}} \, dx = 1.
\]
In other words, $v$ is a minimizer for $S$.
\qed \end{proof} 

\begin{remark}
The previous theorems are contained in the forthcoming paper \cite{SeSmWi}, where we also prove that all minimizers corresponding to the best constant have a rather precise symmetric shape.
\end{remark}

%
%
%

\part{...provided $\varepsilon$ is small enough}

\chapter{The perturbation technique}

\section{Where is \emph{Critical point theory} going?}

Over the last two decades, a lot
of attention has been paid to \emph{Critical point theory}, because of its flexibility and also intuitiveness. 

Besides the Concentration--Compactness Principle already described, a new technique to deal with problems without compactness  has been introduced in the last years by Ambrosetti and Badiale (\cite{AmbBad1,AmbBad2,ABC}), which improves some previous results contained in \cite{ACE}. It allows us to treat problems in which a ``small parameter'' appears, and permits to find existence and multiplicity results even in some situations where the concentration-compactness methods fails or requires heavy calculations. For the sake of completeness, we present below this technique in some detail. Roughly speaking, we can say that one exploits some very precise properties of solutions to \emph{unperturbed equations} to find one or more solutions of a \emph{perturbed equation}.

Suppose we are given a family $\{f_\ge\}$ of $C^2$ functionals over a real Hilbert space, parametrized by $\ge \geq 0$. When $\ge=0$, we have a ``distinguished'' functional $f_0$, and we suppose that $f_0$ has a whole manifold of critical points.\footnote{Very recently, this assuption has been highly relaxed, to cover situations where the points of this manifold are far from being critical. We refer to \cite{AMN}.} Under suitable assumptions, it is possible to prove the existence of critical points of each $f_\ge$, provided $\ge$ is small enough. Let us formalize what we have said.

We want to find critical points of functionals of the form
\begin{equation}\label{eq-abstr}
f_\eps (u)=\frac{1}{2}\|u\|^2 -F(u)+G(\eps,u),
\end{equation}
where $\| \cdot \|$ denotes the norm in an Hilbert space $E$, $F\colon E\to\R$ and $G:\R \times E\to\R$. We assume  that $f_0=\frac{1}{2}\|\cdot\|^2 -F$ has a nondegenerate critical manifold\footnote{We always assume manifold to be smooth enough. As a rule, $C^2$ \emph{is} enough.} $Z$, namely we require that
\begin{description}
\item[(F0)]$F \in C^2$ and $G(0,u)=0$ for all $u\in E$;
\item[(F1)]there exists a $d$-dimensional $C^2$ manifold $Z$ at level $b$ of critical points of $f_0$;
\item[(F2)]$F''(z)$ is compact, for all $z\in Z$;
\item[(F3)]$T_z Z=\ker f_0''(z)$, for all $z\in Z$.
\end{description}
Here $T_zZ$ stands for the tangent space to $Z$ at the point $z\in Z$.
Indicating with $G'$ and $G''$ respectively the first and the second derivative of $G$ with respect to $u$, on $G$ we assume 
\begin{description}
\item[(G0)]$G$ is continuous in $\R \times E$ and $G(0,u)=0$, for all $u \in E$;
\item[(G1)]$G$ is of class $C^2$ with respect to $u \in E$;
\item[(G2)]the maps $(\eps, u) \longmapsto G'(\eps, u)$ and $(\eps, u) \longmapsto G''(\eps, u)$ are continuous; 
\item[(G3)]there exist $\alpha >0$ and a continuous function $\Gamma:Z\to \R$ such that
\[
\Gamma (z)=\lim_{\eps \to 0}\frac{G(\eps,z)}{\eps^\alpha}
\]
and
\[
G'(\eps,z)=o(\eps^{\alpha/2}).
\]
\end{description}

To avoid technicalities, we will assume that $Z=\zeta(\R^d)$ where $\zeta\in C^2(\R^d, E)$. We will set $Z^r=\zeta(B_r)$ where $B_r=\{ x\in \R^d : |x|\le r\}$. 

The main idea is to use the implicit function theorem to produce a good local deformation of $Z$ in the normal direction. More precisely we have the following.

\begin{lemma}\label{lemma:w}
Given $R>0$, there exist $\ge_0 >0$ and a smooth function
\[
w\colon M=Z^R\times (-\ge_0 ,\ge_0) \to E
\]
such that
\begin{description}
\item[(i)] $w(z,0)=0$ for all $z\in Z^R$;
\item[(ii)] $f'_\ge (z+w(z,\ge))\in T_zZ$ for all $(z,\ge)\in M$;
\item[(iii)] $w(z,\ge)$ is orthogonal to $T_zZ$ for all $(z,\ge)\in M$.
\end{description}
\end{lemma}
\begin{proof}
See \cite{AmbBad1}.
\end{proof} 

\begin{lemma}\label{lemma:natural}
The manifold
\[
Z_\ge = \{z+w(z,\ge)\mid (z,\ge)\in M\}
\]
is a natural constraint for $f_\ge$, namely: if $u\in Z_\ge$ and ${f_\ge}_{|Z_\ge}' (u)=0$, then $f_\ge'(u)=0$.
\end{lemma}

\begin{proof}
Suppose that ${f_\ge}_{|Z_\ge}' (u)=0$. Then $f_\ge'(u)$ is orthogonal to $T_uZ_\ge$. On the other hand, $f_\ge(u)\in T_z Z$ and $T_uZ_\ge$ is close to $T_zZ$ for $\ge$ small. This implies $f_\ge'(u)=0$.
\end{proof}

We have the following theorem (see Theorem 1.3 of \cite{AmbBad2}):

\begin{theorem}\label{th1.3}
Suppose (F0-3) and (G0-3) hold and assume that there exist $\delta>0$ and $z^* \in  Z$ such that 
\[
{\text either} \quad \min_{\|z-z^*\|=\delta}\Gamma (z)>\Gamma(z^*) \quad {\text or} \quad \max_{\|z-z^*\|=\delta}\Gamma (z)<\Gamma(z^*). 
\]
Then, for $\eps$ small, $f_\eps $ has a critical point in $Z$.
\end{theorem}

\begin{proof}
We give only the sketch of the proof divided into three steps.

\vspace{0.2cm}
\textbf{Step 1} \quad It is not difficult to prove that the function $w$ built in the previous Lemma is such that
\begin{equation}\label{AB-2.1}
\|w\|=o(\eps^{\alpha/2})\qquad \text{and} \qquad f'_\eps(z+w)\in T_z Z. 
\end{equation}

\vspace{0.2cm}
\textbf{Step 2} \quad Using the Taylor expansion, for $u \in Z_\eps$, one can easily see that 
\[
f_\eps (u)=b+\eps ^\alpha \Gamma (z)+o(\eps ^ \alpha) \qquad (\eps \to 0).
\] 

\vspace{0.2cm}
\textbf{Step 3} \quad It readily follows that $f_\eps$ has a local constrained minimum (or maximum) on $Z_\eps$ at some $u_\eps$. According to Lemma \ref{lemma:natural}, such a $u_\eps$ is a critical point of $f_\eps$.
\end{proof}

Under some additional assumptions, we can have more information about the bifurcation set: we can say that this set is a curve.

Hereafter, we will denote $z_\theta,\,\theta \in \R^d$, the elements of $Z$. We will indicate with $\partial_i z$ and $\partial_{ij} z$ respectively $\partial z / \partial \theta_i$ and $\partial^2 z / \partial \theta_i\partial \theta_j$.

Now we consider two case:

\vspace{0.2cm}
\textbf{First case} \quad Suppose that
\begin{description}
\item[(F4)]$F$ is of class $C^4$;
\item[(G4)]$G$ is of class $C^4$ with respect to $u \in E$, moreover the maps $(\eps,u)\longmapsto G'''(\eps,u)$ is continuous;
\item[(G5)]$G'(\eps,z_\theta)=o(\eps^{\alpha/2})$ and there exist continuous functions $\gamma_{ij},\tilde{\gamma}_{ij}:Z\to\R $ such that if $\theta_\eps \to \theta$ as $\eps \to 0$, there results
\begin{equation}\label{G5-1}
\eps^{-\alpha}\left( G'(\eps,z_{\theta_\eps})|\partial_{ij}z_{\theta_\eps} \right)\to \gamma_{ij}(\theta),
\end{equation}
\begin{equation}\label{G5-2}
\eps^{-\alpha}\left( G''(\eps,z_{\theta_\eps}) \partial_{i}z_\theta |\partial_{j}z_{\theta} \right)\to \tilde{\gamma}_{ij}(\theta),
\end{equation}
\begin{equation}\label{G5-3}
\eps^{-\alpha/2} G''(\eps,z_{\theta_\eps})\partial_{ij}z_{\theta_\eps} \to 0,
\end{equation}
\begin{equation}\label{G5-4}
\eps^{-\alpha/2}G'''(\eps,z_{\theta_\eps})[ \partial_{i}z_\theta,\partial_{j}z_{\theta}] \to 0,
\end{equation}
as $\eps \to 0$.
\end{description}
\vspace{0.2cm}
\textbf{Second case} \quad Suppose that
\begin{description}
\item[(F4)$'$]$F$ is of class $C^3$;
\item[(G4)$'$]$G$ is of class $C^3$ with respect to $u \in E$ and $G'''(\eps, u) \to 0$ as $\eps \to 0$, uniformly in $u \in Z^r$; moreover the map $(\eps,u)\longmapsto G'''(\eps,u)$ is continuous;
\item[(G5)$'$]$G'(\eps,z_\theta)=O(\eps^\alpha)$ and there exist continuous functions $\gamma_{ij},\tilde{\gamma}_{ij}:Z\to\R $ such that if $\theta_\eps \to \theta$ as $\eps \to 0$, there results
\[
\eps^{-\alpha}\left( G'(\eps,z_{\theta_\eps})|\partial_{ij}z_{\theta_\eps} \right)\to \gamma_{ij}(\theta),
\]
\[
\eps^{-\alpha}\left( G''(\eps,z_{\theta_\eps}) \partial_{i}z_\theta |\partial_{j}z_{\theta} \right)\to \tilde{\gamma}_{ij}(\theta),
\]
as $\eps \to 0$.
\end{description}

The next theorem shows that the Morse index of the critical points $u_\ge$ is related to the nature of $u_\ge$ as a critical point of the Melnikov function $\Gamma$. Sometimes, this allows us to conclude that a particular critical point could not be found directly by an application of constrained minimization or of the mountain pass theorem. The Morse index of critical points of \emph{mountain pass type} has been calculated. We refer to \cite{Chang}.

Despite its length, we present the proof, since in the original paper \cite{AmbBad2} there is a mistake. We wish to express our gratitude to A.~Pomponio for pointing out this fact  to us and supplying us with the correct statement and proof. See also \cite{BadPom}.

\begin{theorem}\label{th3.2}
Let $F \in C^3$, (F0-3) and (G0-3) hold. Suppose either (F4) and (G4-5) or (F4)$'$ and (G4-5)$'$ are satisfied. Let $\theta \in \R^d$ be given and let $u_\eps=z(\theta_\eps) +w(\eps,\theta_\eps) \in Z_\eps$ be a critical point of $f_\eps$ on $Z_\eps$ with $\theta_\eps \to \theta$ as $\eps \to 0$. Assume that $z_\theta = \lim_\eps z(\theta_\eps)$ is nondegenerate for the restriction of $f_0$ to $\left(T_{z_\theta}Z  \right)^\perp$, with Morse index equal to $m_0$, and that the Hessian $D^2 \Gamma(\theta)$ is positive (negative) definite.

Then $u_\eps$ in a nondegenerate critical point for $f_\eps$ with Morse index equal to $m_0$, if $D^2 \Gamma(\theta)>0$, and to $m_0+d$, if $D^2 \Gamma(\theta)<0$. As a consequence, the critical points of $f_\eps$ form a continuous curve.
\end{theorem}

\begin{proof}
The first part of the proof works for all the two cases.

We write $E=E^{+} \oplus E^{0} \oplus E^{-}$ where $E^{0}=T_{z}Z$, $\dim (E^{-})=m_{0}$ and
\begin{equation}\label{><0}
\begin {cases}
D^{2}f_{0}(z)[v,v]>0 \quad \forall v\in E^{+}\\
D^{2}f_{0}(z)[v,v]<0 \quad \forall v\in E^{-}.
\end {cases}
\end{equation}
Moreover we have 
\begin{equation*}
D^{2}f_{0}(z)[\partial_{i}z ,\partial_{j}z]=0.
\end{equation*}
We denote with $\varphi_{i}^{0}=\partial_{i}z$, for $i=1, \ldots , d$, and with $t_{i}^{0}$, for $i=1, \ldots , m_{0}$, respectively the orthonormal base of $E^{0}$ and $E^{-}$.

We call now $E_{\eps}^{0}=T_{u_{\eps}}Z$ and we will denote with $\varphi^{\eps}_{i}$, for $i=1, \ldots , d$, its orthonormal base. 

For $i=1, \ldots , m_{0}$, we want to find $\alpha_{j}^{\eps}$ such that if $t_{i}^{\eps}=t_{i}^{0}+\alpha_{i}^{\eps}$ then $(t_{i}^{\eps} \vert \varphi_{j}^{\eps})=0$, for all $j=1,\ldots, d$. So we must have 
\begin{equation*}
0=(t_{i}^{\eps} \vert \varphi_{j}^{\eps})=(t_{i}^{0}+ \alpha_{i}^{\eps} \vert \varphi_{j}^{0}+(\varphi_{j}^{\eps}-\varphi_{j}^{0}))=(t_{i}^{0}\vert \varphi_{j}^{\eps}-\varphi_{j}^{0}) +(\alpha_{i}^{\eps} \vert \varphi_{j}^{\eps}),
\end{equation*}
therefore
\begin{equation*}
(\alpha_{i}^{\eps} \vert \varphi_{j}^{\eps})=-(t_{i}^{0}\vert \varphi_{j}^{\eps}-\varphi_{j}^{0})
\end{equation*}
and we take $\alpha_{i}^{\eps}=\sum_{j=1}^{d} -(t_{i}^{0}\vert \varphi_{j}^{\eps}-\varphi_{j}^{0})\varphi ^{\eps}_{j}$. 

Let $E^{-}_{\eps}$ be the space spanned by $\left\{\frac{t^{\eps}_{1}}{\vert t^{\eps}_{1}\vert},\ldots, \frac{t^{\eps}_{m_{0}}}{\vert t^{\eps}_{m_{0}}\vert}\right\}$. We have
\begin{gather*}
D^{2}f_{\eps}(u_{\eps})[t_{i}^{\eps} ,t_{j}^{\eps}]=D^{2}f_{\eps}(u_{\eps})[t_{i}^{0}+\alpha^{\eps}_{i} ,t_{j}^{0}+\alpha^{\eps}_{j}]=
\\
=D^{2}f_{\eps}(u_{\eps})[t_{i}^{0} ,t_{j}^{0}]+D^{2}f_{\eps}(u_{\eps})[t_{i}^{0} ,\alpha_{j}^{\eps}]
+D^{2}f_{\eps}(u_{\eps})[\alpha_{i}^{\eps},t_{j}^{0}]+D^{2}f_{\eps}(u_{\eps})[\alpha_{i}^{\eps} ,\alpha_{j}^{\eps}].
\end{gather*}
Since $\alpha^{\eps}_{i}$ goes to $0$, as $\eps \to 0$, by (\ref{><0}) we have
\begin{equation}\label{<0}
D^{2}f_{\eps}(u_{\eps})[t_{i}^{\eps} ,t_{j}^{\eps}]<0
\end{equation}
and therefore $D^{2}f_{\eps}(u_{\eps})$ is negative definite on $E^{-}_{\eps}$, for $\eps $ small.

Let us put now $E^{+}_{\eps}=(E^{0}_{\eps}\oplus E^{-}_{\eps})^{\perp}$. Of course we have
\[
E=E^{+}_{\eps}\oplus E^{0}_{\eps}\oplus E^{-}_{\eps}.
\]
We want to show that $D^{2}f_{\eps}(u_{\eps})$ is positive definite on $E^{+}_{\eps}$. Let $v_{\eps}\in E^{+}_{\eps}$ 
with $\vert v_{\eps}\vert =1$. Let $P^{+}$ indicate the orthogonal projection on $E^{+}$, then
\begin{gather*}
D^{2}f_{\eps}(u_{\eps})[v_{\eps},v_{\eps}]=D^{2}f_{\eps}(u_{\eps})[P^{+}v_{\eps}+(v_{\eps} -P^{+}v_{\eps}),P^{+}v_{\eps}+(v_{\eps} -P^{+}v_{\eps})]=
\\
=D^{2}f_{\eps}(u_{\eps})[P^{+}v_{\eps},P^{+}v_{\eps}]+2D^{2}f_{\eps}(u_{\eps})[P^{+}v_{\eps},v_{\eps} -P^{+}v_{\eps}]+
\\
+D^{2}f_{\eps}(u_{\eps})[v_{\eps} -P^{+}v_{\eps},v_{\eps} -P^{+}v_{\eps}].
\end{gather*}
We claim that there exist $\delta>0$ and $\eps_{0}>0$ such that for all $\eps<\eps_{0}$ and for all $v_{\eps}\in E^{+}_{\eps}$ we have
\begin{equation}\label{>0}
D^{2}f_{\eps}(u_{\eps})[v_{\eps},v_{\eps}] \geq \delta.
\end{equation}
By contradiction, suppose the contrary, i.e. for all $k\in \N$ there exist $\eps_{k}$ and $v_{\eps_{k}}\in E^{+}_{\eps_{k}}$ such that
\[
D^{2}f_{\eps_{k}}(u_{\eps_{k}})[v_{\eps_{k}},v_{\eps_{k}}] < \frac{1}{k}.
\] 
So we should have
\begin{gather*}
\frac{1}{k}>D^{2}f_{\eps_{k}}(u_{\eps_{k}})[P^{+}v_{\eps_{k}},P^{+}v_{\eps_{k}}]+2D^{2}f_{\eps_{k}}(u_{\eps_{k}})[P^{+}v_{\eps_{k}},v_ {\eps_{k}}-P^{+}v_{\eps_{k}}]+
\\
+D^{2}f_{\eps_{k}}(u_{\eps_{k}})[v_ {\eps_{k}}-P^{+}v_{\eps_k},v_{\eps_k} -P^{+}v_{\eps_k}]
\end{gather*}
Recall that 
\[
v_{\eps_k}-P^{+}v_{\eps_k}=\sum_{i=1}^{m_{0}}(v_{\eps_k}|t_{i}^{0})t_{i}^{0} + \sum_{i=1}^{d}(v_{\eps_k}|\varphi_{i}^{0})\varphi_{i}^{0}.
\]
and notice that, since $v_{\eps_k}\in E^{+}_{\eps_{k}}$, $(v_{\eps_k}|t_{i}^{\eps})=0$ and $(v_{\eps_k}|\varphi_{i}^{\eps})=0$. We have
\[
(v_{\eps_k}|t_{i}^{0})=(v_{\eps_k}|t_{i}^{\eps})+(v_{\eps_k}|t_{i}^{0}-t_{i}^{\eps}),
\]
the first term is zero, the second one goes to zero since $\{v_{\eps_k}\}$ is a bounded sequence. Similarly we have that also $(v_{\eps_k}|\varphi_{i}^{0})$ goes to zero and so
\[
v_{\eps_k}-P^{+}v_{\eps_k} \to 0.
\] 
This is a contradiction since, for $k$ sufficiently large, by (\ref{><0}) we have
\[
D^{2}f_{\eps_{k}}(u_{\eps_{k}})[v_{\eps_k},v_{\eps_k}]\geq \bar{\delta}>0.
\]

At this moment, by (\ref{<0}) and (\ref{>0}), we have shown that $D^{2}f_{\eps}(u_{\eps})$ is negative definite on $E^{-}_{\eps}$ and positive definite on $E^{+}_{\eps}$. We now want to investigate the behavior of $D^{2}f_{\eps}(u_{\eps})$ on $E^{0}_{\eps}$. 

\vspace{0.2cm}
\textbf{First case} \quad We have
\begin{multline}\label{D2}
D^{2}f_{\eps}(u_{\eps})[\partial_{i} z_{\theta_{\eps}},\partial_{j} z_{\theta_{\eps}}]=
\\
=(\partial_{i} z_{\theta_{\eps}}|\partial_{j} z_{\theta_{\eps}})
-F''(u_\eps)[\partial_{i} z_{\theta_{\eps}},\partial_{j} z_{\theta_{\eps}}]
+G''(\eps, u_\eps)[\partial_{i} z_{\theta_{\eps}},\partial_{j} z_{\theta_{\eps}}]
\end{multline}
We know by \textbf{(F3)} that $\partial_{i} z_{\theta_{\eps}}\in \ker [I_E-F''(z_{\theta_{\eps}})]$ and so recalling that $w(0,z)=0$ and expanding $F''(u_\eps)$ and $G''(\eps,u_\eps)$, (\ref{D2}) becomes 
\begin{multline*}
D^{2}f_{\eps}(u_{\eps})[\partial_{i} z_{\theta_{\eps}},\partial_{j} z_{\theta_{\eps}}]=
\\
=(\partial_{i} z_{\theta_{\eps}}|\partial_{j} z_{\theta_{\eps}})
-F''(z_{\theta_{\eps}})[\partial_{i} z_{\theta_{\eps}},\partial_{j} z_{\theta_{\eps}}]
-\left(F'''(z_{\theta_{\eps}})[\partial_{i} z_{\theta_{\eps}},\partial_{j} z_{\theta_{\eps}}]|w_\eps\right)+
\\
+G''(\eps, z_{\theta_{\eps}})[\partial_{i} z_{\theta_{\eps}},\partial_{j} z_{\theta_{\eps}}]
+\left(G'''(\eps,z_{\theta_{\eps}})[\partial_{i} z_{\theta_{\eps}},\partial_{j} z_{\theta_{\eps}}]|w_\eps\right)+O(\|w_\eps\|^2)=
\\
=-\left(F'''(z_{\theta_{\eps}})[\partial_{i} z_{\theta_{\eps}},\partial_{j} z_{\theta_{\eps}}]|w_\eps\right)
+G''(\eps, z_{\theta_{\eps}})[\partial_{i} z_{\theta_{\eps}},\partial_{j} z_{\theta_{\eps}}]+
\\
+\left(G'''(\eps,z_{\theta_{\eps}})[\partial_{i} z_{\theta_{\eps}},\partial_{j} z_{\theta_{\eps}}]|w_\eps\right)+O(\|w_\eps\|^2)
\end{multline*}
We observe explicitly that, since $w_\eps =o(\eps^{\alpha/2})$ and by (\ref{G5-4}) of \textbf{(G5)} we can write
\begin{multline}\label{def-pos0}
D^{2}f_{\eps}(u_{\eps})[\partial_{i} z_{\theta_{\eps}},\partial_{j} z_{\theta_{\eps}}]= 
\\
=-\left(F'''(z_{\theta_{\eps}})[\partial_{i} z_{\theta_{\eps}},\partial_{j} z_{\theta_{\eps}}]|w_\eps\right)
+G''(\eps, z_{\theta_{\eps}})[\partial_{i} z_{\theta_{\eps}},\partial_{j} z_{\theta_{\eps}}]+o(\eps^\alpha) 
\end{multline}
 
By (\ref{AB-2.1}) we have 
\[
z_{\theta_{\eps}}+w_\eps -F'(z_{\theta_{\eps}}+w_\eps)+G'(\eps, z_{\theta_{\eps}}+w_\eps)=\sum_l a_l \partial_l z_{\theta_{\eps}},
\]
and so a development of $F'(z_{\theta_{\eps}}+w_\eps)$ and $G'(\eps, z_{\theta_{\eps}}+w_\eps)$ yields
\begin{equation}\label{def-pos1}
z_{\theta_{\eps}}+w_\eps -F'(z_{\theta_{\eps}})-F''(z_{\theta_{\eps}})w_\eps+G'(\eps, z_{\theta_{\eps}})+G''(\eps, z_{\theta_{\eps}})w_\eps=\sum_l a_l \partial_l z_{\theta_{\eps}}.
\end{equation}
Without loss of generality, we can assume that the parametrization of $Z$ is normal to $\theta$, namely that $\partial_{ij} z_{\theta_{\eps}}\perp T_{z_{\theta_{\eps}}} Z$. Taking the scalar product of (\ref{def-pos1}) with $\partial_{ij} z_{\theta_{\eps}}$ and, since $z_{\theta_{\eps}}=F'(z_{\theta_{\eps}})$, we find
\begin{multline}\label{def-pos2}
(w_\eps|\partial_{ij} z_{\theta_{\eps}}) -(F''(z_{\theta_{\eps}})w_\eps|\partial_{ij} z_{\theta_{\eps}})+\\ (G'(\eps, z_{\theta_{\eps}})|\partial_{ij} z_{\theta_{\eps}})+(G''(\eps, z_{\theta_{\eps}})w_\eps|\partial_{ij} z_{\theta_{\eps}})=0
\end{multline}
By (\ref{G5-3}) of \textbf{(G5)}, (\ref{def-pos2}) becomes
\begin{equation}\label{def-pos3}
(w_\eps|\partial_{ij} z_{\theta_{\eps}}) -(F''(z_{\theta_{\eps}})w_\eps|\partial_{ij} z_{\theta_{\eps}})+(G'(\eps, z_{\theta_{\eps}})|\partial_{ij} z_{\theta_{\eps}})+o(\eps^\alpha)=0
\end{equation}

From $z_{\theta_{\eps}}=F'(z_{\theta_{\eps}})$, it follows
\begin{equation}\label{def-pos4}
\partial_{ij} z_{\theta_{\eps}}=F''(z_{\theta_{\eps}})\partial_{ij} z_{\theta_{\eps}}+F'''(z_{\theta_{\eps}})[\partial_{i} z_{\theta_{\eps}},\partial_{j} z_{\theta_{\eps}}].
\end{equation}
The scalar product of (\ref{def-pos4}) with $w_\eps$ is
\begin{equation}\label{def-pos5}
(\partial_{ij} z_{\theta_{\eps}}|w_\eps)=(F''(z_{\theta_{\eps}})\partial_{ij} z_{\theta_{\eps}}|w_\eps)+(F'''(z_{\theta_{\eps}})[\partial_{i} z_{\theta_{\eps}},\partial_{j} z_{\theta_{\eps}}]|w_\eps).
\end{equation}

Substituting (\ref{def-pos5}) in (\ref{def-pos3}), we have
\begin{equation}\label{def-pos6}
-(F'''(z_{\theta_{\eps}})[\partial_{i} z_{\theta_{\eps}},\partial_{j} z_{\theta_{\eps}}]|w_\eps)=(G'(\eps, z_{\theta_{\eps}})|\partial_{ij} z_{\theta_{\eps}})+o(\eps^\alpha).
\end{equation}

Dividing by $\eps^{\alpha}$ and passing to the limit for $\eps \to 0$, by (\ref{G5-1}) and (\ref{G5-2}) of \textbf{(G5)} we have:
\begin{gather*}
-\eps^{-\alpha}\left(F'''(z_{\theta_{\eps}})[\partial_{i} z_{\theta_{\eps}},\partial_{j} z_{\theta_{\eps}}]|w_\eps\right)=\eps^{-\alpha} (G'(\eps, z_{\theta_{\eps}})|\partial_{ij} z_{\theta_{\eps}})+o(1)\to \gamma_{ij}(\theta),
\\
\eps^{-\alpha} G''(\eps, z_{\theta_{\eps}})[\partial_{i} z_{\theta_{\eps}},\partial_{j} z_{\theta_{\eps}}] \to \tilde{\gamma}_{ij}(\theta).
\end{gather*}
Therefore, by (\ref{def-pos0})
\[
\lim_{\eps \to 0}\frac{1}{\eps^{\alpha}}D^{2}f_{\eps}(u_{\eps})[\partial_{i} z_{\theta_{\eps}},\partial_{j} z_{\theta_{\eps}}]=\gamma_{ij}(\theta)+\tilde{\gamma}_{ij}(\theta).
\]
On the other side we have
\begin{gather*}
\partial_{ij}\Gamma(\theta)=\lim_{\eps \to 0}\frac{1}{\eps^{\alpha}}\left(G'(\eps,z)|\partial_{ij} z\right)
+\lim_{\eps \to 0}\frac{1}{\eps^{\alpha}}G''(\eps,z)[\partial_{i} z,\partial_{j} z]=\gamma_{ij}(\theta)+\tilde{\gamma}_{ij}(\theta).
\end{gather*}

\vspace{0.2cm}
\textbf{Second case} \quad To study the behavior of $D^{2}f_{\eps}(u_{\eps})$ on $E^{0}_{\eps}$, we note that, by \textbf{(G5)$'$} and recalling that $w_\eps =O(\eps^\alpha)$ (see Lemma 2.2 of \cite{AmbBad2}), one can show easily that (\ref{def-pos0}) and (\ref{def-pos6}) still hold and so the proof goes as in our case.
\end{proof}

In \cite{AmbBad1}, the case of a perturbation like
\[
f_\ge = f_0 + \ge G
\]
was considered, leading to to easier statements and proofs. For some applications, we refer to the bibliography.

\vspace{2cm}

We present in the next chapters a few applications of the perturbation technique. The reader will see that the strategy is rather clear, despite some modifications due to the specific shape of the perturbation terms. In the case of the \emph{Schr\"{o}dinger equation}, the perturbation term is not globally smooth in the variables $(\ge ,u)\in \R\times H^1(\Rn)$, so that some trick is needed in order to apply the finite--dimensional reduction.

\chapter{Closed geodesics on cylinders}

In this chapter, we concentrate on the existence of closed geodesics on a non--compact manifold $M$. First, we gather some results from Riemannian geometry.

\section{Tools from differential geometry}

We assume some knowledge of what a smooth differentiable manifold is, together with the main structures: tangent and cotangent spaces, inner products, and so on.

\begin{definition}
Let $M$ be a smooth differentiable manifold. A \emph{Riemannian metric} on $M$ is a tensor field $g\colon C_2^\infty (TM)\to C_2^\infty (TM)$ such that for each $p\in M$ the restriction $g_p = g_{|T_p M \otimes T_p M} \colon T_p M \otimes T_p M \to \R$ with
\[
g_p \colon (X_p ,Y_p) \mapsto g(X,Y)(p)
\]
is an inner product on the vector space $T_pM$ The pair $(M,g)$ is acalled a \emph{Riemannian manifold}.
\end{definition}

\begin{definition}
Let $(E,M,\pi)$ be a smooth vector bundle over $M$. A \emph{connection} on $(E,M,\pi)$ is a map $\nabla \colon C^\infty (TM)\times C^\infty (E)\to C^\infty (E)$ such that
\begin{enumerate}
\item $\nabla_X (\lambda \cdot v + \mu \cdot w)=\lambda \cdot \nabla_X v + \mu\cdot\nabla_X w$;
\item $\nabla_X (f\cdot v)=X(f)\cdot v + f\cdot \nabla_X v$;
\item $\nabla_{(f\cdot X + g\cdot Y)} v = f\cdot \nabla_X v + g\cdot \nabla_Y v$
\end{enumerate}
for all $\lambda$, $\mu\in\R$, $X$, $Y\in C^\infty (TM)$, $v$, $w\in C^\infty (E)$, and $f$, $g\in C^\infty (M)$. A section $v\in C^\infty (E)$ is called \emph{parallel} with respect to the connection if $\nabla_X v=0$ for all $X\in C^\infty (TM)$.
\end{definition}

\begin{definition}
Let $M$ be a smooth manifold, and $\nabla$ be a connection on the tangent bundle $(TM,M,\pi)$. The \emph{torsion} of $\nabla$ is the map $T\colon C_2^\infty (TM)\to C_1^\infty (TM)$ defined by
\[
T(X,Y)=\nabla_X Y - \nabla_Y X - [X,Y],
\]
where $[\, , \, ]$ is the Lie bracket on $C^\infty (M)$. A connection $\nabla$ is \emph{torsion--free} if $T\equiv 0$.

If $g$ is a Riemannian metric on $M$, then $\nabla$ is said to be \emph{metric} if
\[
X(g(Y,Z))=g(\nabla_X Y,Z)+g(Y,\nabla_X Z)
\]
for all $X$, $Y$, $Z\in C^\infty (TM)$.
\end{definition}

\begin{theorem}
Let $(M,g)$ be a Riemannian manifold. There exists one and only one connection $\nabla$ on $M$ such that
\begin{multline*}
g(\nabla_X Y,Z)=\frac{1}{2} \left( X(g(Y,Z))+Y(g(Z,X))-Z(g(X,Y))+\right.\\
\left. + g(Z,[X,Y])+g(Y,[Z,X])-g(X,[Y,Z])\right)
\end{multline*}
for all $X$, $Y$, $Z\in C^\infty (TM)$. This connection is called the \emph{Levi-Civita connection} on $M$.
\end{theorem}

\begin{definition}
Let $M$ be a smooth manifold and $(TM,M,\pi)$ be its tangent bundle. A \emph{vector field} $X$ \emph{along a curve} $\gamma \colon I \to M$ is a curve $X\colon I \to TM$ such that $\pi \circ X=\gamma$. By $C_\gamma^\infty (TM)$ we denote the set of all smooth vector fields along $\gamma$.
\end{definition}

\begin{definition}
Let $(M,g)$ be a Riemannian manifold and $\gamma \colon I \to M$ be a $C^1$--curve. A vector field $X$ along $\gamma$ is said to be \emph{parallel} along $\gamma$ if
\[
\nabla_{\dot \gamma} X=0.
\]
A $C^2$--curve $\gamma \colon I \to M$ is said to be a \emph{geodesic} if the vector field $\dot\gamma$ is parallel along $\gamma$, i.e.\[
\nabla_{\dot\gamma}\dot\gamma =0.
\]
\end{definition}

Geodesics possess a variational characterization that we now clarify.

\begin{definition}
Let $(M,g)$ be a Riemannian manifold and $\gamma \colon I \to M$ be a $C^2$--curve on $M$. For every compact interval $[a,b]\subset I$ we define the energy functional
\[
E_{[a,b]}(\gamma)=\frac{1}{2}\int_a^b g(\dot\gamma (t),\dot\gamma (t))\, dt.
\]
\end{definition}

\begin{theorem}
A $C^2$--curve $\gamma$ is a critical point for the energy functional if and only if it is a geodesic.
\end{theorem}

\begin{definition}
A geodesic is \emph{closed} if its image is diffeomorphic to $S^1$.
\end{definition}

When $M$ is compact, closed geodesics can be found by standard variational arguments. For a survey of standard and advanced results on the existence of closed geodesics, we refer to \cite{kli3}.

On the other hand, when $M$ is not compact, the existence of critical points of the energy $E$ is no longer evident. Indeed, very few results are known in the literature, see \cite{benci,ku,tho}.
In particular, Tanaka deals with the manifod $M=\mathbb{R} \times S^{N}$, endowed with a metric $g(s,\xi)=g_0(\xi)+h(s,\xi)$, where $g_0$ is the standard product metric on $\mathbb{R} \times S^{N}$.
Under the assumption that $h(s,\xi)\to 0$ as $|s| \to \infty$, he proves the existence of a closed geodesic, found as a critical point of the energy functional
\begin{equation}
E(u)=\frac{1}{2} \int_0^1 g(u)[\dot u , \dot u] \, dt,
\end{equation}
defined on the loop space $\Lambda = \Lambda (M) = H^1 (S^1,M)$. \footnote{We will identify $S^1$ with $[0,1] / \{ 0,1 \}$.} The lack of compactness due to the unboundedness of $M$ is overcome by a suitable use of the concentration--compactness principle. To carry out the proof, the fact that $M$ has the specific form $M = \mathbb{R} \times S^{N}$ is fundamental, because this permits to compare $E$ with a {\em functional at infinity} whose behavior is explicitely known.

Here we consider a perturbed metric $g_\ge = g_0 + \ge h$, and extend Tanaka's result in two directions. First, we show the existence of at least $N$, in some cases $2N$, closed geodesics on $M=\mathbb{R} \times S^N$, see Theorem \ref{th:2.1}. Such a theorem can also be seen as an extension to cilyndrical domains of the result by Carminati \cite{carminati}. Next, we deal with the case in which $M=\mathbb{R} \times M_0$ for a general compact $N$--dimensional manifold $M_0$. \footnote{By manifold we mean a smooth, connected manifold.} The existence result we are able to prove requires that either $M_0$ possesses a non--degenerate closed geodesic, see Theorem \ref{th:3.1}, or that $\pi_1 (M_0) \ne \{0\}$ and the geodesics on $M_0$ are {\em isolated}, see Theorem \ref{th:4.1}.

\begin{figure}
\centering
\includegraphics[height=34mm]{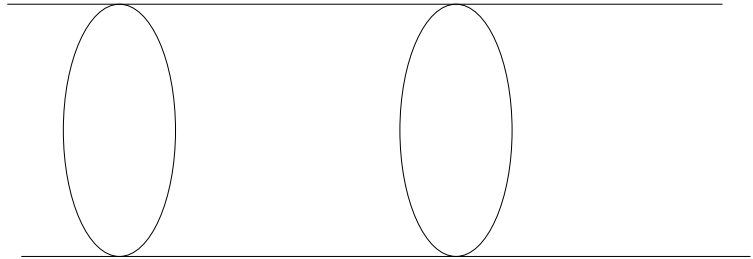}
\caption{Unperturbed cylinder}
\end{figure}

\begin{figure}
\centering
\includegraphics[height=34mm]{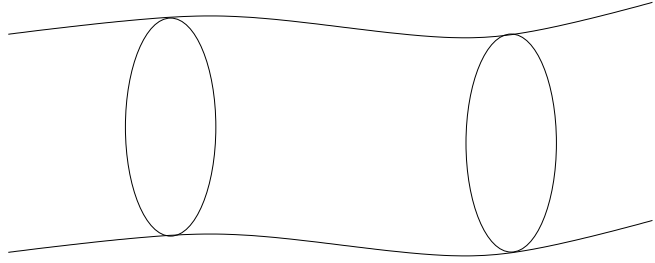}
\caption{Perturbed cylinder}
\end{figure}

The approach we use is different that Tanaka's one, and relies on the perturbation result discussed before. Roughly, the main advantages of using this abstract perturbation method are that
\begin{description}
\item[(i)] we can obtain sharper results, like the multeplicity ones;
\item[(ii)] we can deal with a general manifold like $M=\mathbb{R} \times M_0$, not only $M= \mathbb{R} \times S^N$, when the results -- for the reasons indicated before -- cannot be easily obtained by using Tanaka's approach.
\end{description}

\section{Spheres}

In this section we assume $M=\mathbb{R} \times S^N$, where $S^N=\{ \xi \in \mathbb{R}^{N+1} \colon | \xi | = 1 \}$ \footnote{Hereafter, we use the notation $\xi \bullet \eta = \sum_i \xi_i \eta_i$ for the scalar product in $\mathbb{R}^{N+1}$, and $|\xi|^2 = \xi \bullet \xi$.}.
For $s \in \mathbb{R}$, $r \in T_s \mathbb{R} \approx \mathbb{R}$, $\xi \in S^N$, $\eta \in  T_\xi S^N$, let
\begin{equation}
g_0(s,\xi)((r,\eta),(r,\eta))=|r|^2 + |\eta|^2
\end{equation}
be the standard product metric on $M=\mathbb{R} \times S^N$. We consider a perturbed metric
\begin{equation}
g_\ge (s,\xi)((r,\eta),(r,\eta))=|r|^2+|\eta|^2 + \ge h(s,\xi)((r,\eta),(r,\eta)),
\end{equation}
where $h(s,\xi)$ is a bilinear form, not necessarily positive definite.

Define the space of closed loops
\begin{equation} \label{Lambda}
\Lambda = \{ u = (r,x) \in H^1 (S^1,\mathbb{R}) \times H^1 (S^1 , S^N)\}
\end{equation}

Closed geodesics on $(M,g_\ge)$ are the critical points of $E_\ge \colon \Lambda \to \mathbb{R}$ given by
\begin{equation}\label{energia:perturbata}
E_\ge (u)=\frac{1}{2}\int_0^1 g_\ge (u)[\dot u , \dot u] \, dt.
\end{equation}
One has that
\begin{equation}\label{eq:vabene}
E_\ge (u)=E_\ge(r,x) = E_0(r,x)+\ge G(r,x),
\end{equation}
where
\begin{equation*}
E_0(r,x)=\frac{1}{2}\int_0^1 \left( |\dot r |^2 + | \dot x |^2 \right) \, dt
\end{equation*}
and
\begin{equation}\label{eq:G}
G(r,x)=\frac{1}{2} \int_0^1 h(r,x)[(\dot r,\dot x),(\dot r, \dot x)] \, dt.
\end{equation}
In particular, we split $E_0$ into two parts, namely
\begin{equation}
E_0(r,x)=L_0(r)+E_{M_0}(x),
\end{equation}
where
\[
L_0(r)=\frac{1}{2} \int_0^1 |\dot r|^2 \, dt, \quad E_{M_0}(x)=\frac{1}{2} \int_0^1 | \dot x|^2 \, dt.
\]

The form of $E_\ge$ suggests to apply the perturbative results of \cite{AmbBad1} in the following form.

\begin{theorem}\label{th:berti:bolle}
Let $H$ be a real Hilbert space, $E_\ge \in C^2 (H)$ be of the form
\begin{equation}\label{eq:vabeneastratto}
E_\ge (u)=E_0(u)+\ge G(u),
\end{equation}
 where $G \in C^2 (E)$.. Suppose that there exists a finite dimensional manifold $Z$ such that
\begin{description}
\item[(AS1)] $E_0^\prime (z) = 0$ for all $z \in Z$;
\item[(AS2)] $E_0^{\prime \prime} (z)$ is a compact perturbation of the identity, for all $z \in Z$;
\item[(AS3)] $T_z Z= \ker E_0^{\prime \prime} (z)$ for all $z \in Z$.
\end{description}
There exist a positive number $\ge_0$ and a smooth function $w \colon Z \times (-\ge_0 , \ge_0) \to~H$ such that the critical points of
\begin{equation}\label{phi}
\Phi_\ge (z)=E_\ge (z+w(z,\ge)), \quad z \in Z,
\end{equation}
are critical points of $E_\ge$.
\end{theorem}

Moreover, it is possible to show that
\begin{equation}\label{eq:gamma}
\Phi_\ge (z)=b+\ge \Gamma (z) + o(\ge),
\end{equation}
where $b=E_0 (z)$ and $\Gamma = G_{|Z}$.
From this ``first order'' expansion, one infers

\begin{theorem}[\cite{AmbBad1}]
Let $H$ be a real Hilbert space, $E_\ge \in C^2 (H)$ be of the form \eqref{eq:vabeneastratto}. Suppose that
(AS1)--(AS3) hold.
Then any strict local extremum of $G_{|Z}$ gives rise to a critical point of $E_\ge$, for $|\ge|$ sufficiently small.
\end{theorem}

\vspace{10pt}

In the present situation, the critical points of $E_0$ are nothing but the great circles of $S^N$, namely
\begin{equation}\label{zetapq}
z_{p,q}=p \cos 2 \pi t + q \sin 2 \pi t,
\end{equation}
where $p,q \in \mathbb{R}^{N+1}$, $p \bullet q = 0$, $|p|=|q|=1$.
Hence $E_0$ has a ``critical manifold'' given by
\begin{equation*}
Z=\{z(r,p,q)=(r,z_{p,q}(\cdot)) \mid r \in \mathbb{R}, \; z_{p,q} \mbox{\; as in \eqref{zetapq}} \}.
\end{equation*}

\begin{lemma}
$Z$ satisfies (AS2)--(AS3).
\end{lemma}

\begin{proof}
The first assertion is known, see for instance \cite{kli1}.

For the second statement, we closely follow \cite{carminati}.

For $z \in Z$, of the form $z(t)=(r,z_{p,q}(t))$, it turns out that
\[
 E_0^{\prime \prime} (z)[h,k]=\int_0^1 \left[ \dot h \bullet \dot k - |\dot z|^2 h \bullet k \right] \, dt
\]
for any $h,k \in T_z Z$.

Let $e_i \in \mathbb{R}^{N+1}$, $i=2, \dots,N+1$, be orthonormal vectors such that $\{ \frac{1}{2 \pi} {\dot z}_{p,q}, e_2, \dots e_{N+1}\}$ is a basis of $T_z Z$, and set
\begin{equation*}
e_i (t)=
\begin{cases}
\dot z_{u^1,u^2}(t) / 2 \pi &\text{ if } i=1 \\
e_i &\text{ if } i > 1,
\end{cases}
\end{equation*}

Then, for $h,k$ as before, we can write a ``Fourier--type'' expansion
\begin{equation} \label{acca}
h(t)=h_0 (t) \frac{d}{dt} + \sum_{i=1}^{N-1} h_i (t) e_i(t), \quad k(t)=k_0(t) \frac{d}{dt} + \sum_{i=1}^{N-1} k_i (t) e_i (t).
\end{equation}

Assume now that $h \in \ker E_0^{\prime \prime}(z_{p,q})$, i.e.
\[
\int_0^1 \dot h \bullet \dot k \, dt = \int_0^1 |\dot z|^2 h \bullet k \, dt \quad \forall k \in T_{z_{p,q}} Z .
\]

We plug \eqref{acca} into this relations, and we get the system

\begin{equation}
\begin{cases}
\overset{..}{h_1} = 0 \\
\overset{..}{h_j} + 4 \pi^2 h_j = 0 \quad j=2, \dots , N-1 \\
\overset{..}{h_0} = 0.
\end{cases}
\end{equation}

Recalling that $h_0$ and $h_1$ are periodic, we find
\begin{equation}
\begin{cases}
h_0 = \lambda_0, \quad h_1 = \lambda_1 \\
h_j = \lambda_j \cos 2 \pi t + \mu_j \sin 2 \pi t \quad j=2, \dots N-1 .
\end{cases}
\end{equation}

Therefore, $h \in T_z Z$. This shows that $\ker E_0^{\prime \prime}(z_{p,q}) \subset T_{z_{p,q}}Z$. Since the converse inclusion is is always true, the lemma follows.
\end{proof}

\begin{lemma} \label{lem:due}
Suppose
\begin{description}
\item[(h1)] $h(r,\cdot) \to 0$ pointwise on $S^N$, as $|r| \to \infty$,
\end{description}
then
\[
\Phi_\ge \to b \equiv E_0(z).
\]
Recall that $\Phi_\ge$ was defined in \eqref{phi}.
\end{lemma}

\begin{proof}
This is proved as in \cite{AM,bertibolle}. We just sketch the argument. The idea is to use the contraction mapping principle to characterize the function $w(\ge ,z)$ (see Theorem 1). Indeed, define
\begin{equation*}
H(\alpha , w,z_r,\ge )=
\begin{pmatrix}
E_\ge^\prime (z_r +w)-\alpha \dot z \\
w \bullet \dot z.
\end{pmatrix}
\end{equation*}
So $H=0$ if and only if $w \in (T_{z_r} Z)^\bot$ and $E_\ge^\prime (z_r+w) \in T_{z_r}Z$. Now,
\begin{multline*}
H(\alpha , w , z_r , \ge)=0 \Leftrightarrow H(0,0,z_r,0)+ \frac{\de H}{\de (\alpha ,w)} (0,0,z_r,0)[\alpha ,w] + \\ R(\alpha ,w,z_r, \ge ) = 0,
\end{multline*}
where $R(\alpha ,w,z_r, \ge ) = H(\alpha ,w,z_r,\ge)-\frac{\de H}{\de (\alpha ,w)} (0,0,z_r,0)[\alpha ,w]$.

Setting
\[
R_{z_r,\ge}(\alpha ,w)=- \left[ \frac{\de H}{\de (\alpha ,w)} (0,0,z_r,0) \right]^{-1}R(\alpha ,w,z_r, \ge),
\]
one finds that
\[
H(\alpha , w , z_r , \ge)=0 \Leftrightarrow (\alpha ,w) = R_{z_r,\ge}(\alpha ,w).
\]
By the Cauchy--Schwarz inequality, it turns out that $R_{zr,\ge}$ is a contraction mapping from some ball $B_{\rho (\ge)}$ into itself. If $|\ge|$ is sufficiently small, we have proved the existence of $(\alpha , w)$ uniformly for $z_r \in Z$. We want to study the asymptotic behavior of $w=w(\ge,z_r)$ as $|r| \to + \infty$. We denote by $R_\ge^0$ the functions $R_{z_r,\ge}$ corresponding to the unperturbed energy functional $E_0=E_{M_0}$. It is easy to see (\cite{bertibolle}, Lemma 3) that the function $w^0$ found with the same argument as before satisfies $\|w^0 (z_r)\| \to 0$ as $|r| \to + \infty$. Thus, by the continuous dependence of $w(\ge ,z_r)$ on $\ge$ and the characterization of $w(\ge ,z_r)$ and $w^0$ as fixed points of contractive mappings, we deduce as in \cite{AM}, proof of Lemma 3.2, that $\lim_{r \to \infty} w(\ge ,z_r)=0$. In conclusion, we have that $\lim_{|r| \to + \infty} \Phi_\ge (z_r + w(\ge ,z_r))=E_{M_0}(z_0)$.

\end{proof}

\begin{remark} \label{rem:3}
There is a natural action of the group $O(2)$ on the space $\Lambda$, given by
\begin{align*}
\{ \pm 1 \} \times S^1 \times \Lambda &\longrightarrow \Lambda \\
(\pm 1 , \theta , u) & \mapsto u(\pm t + \theta ),
\end{align*}
under which the energy $E_\ge$ is invariant. Since this is an isometric action under which $Z$ is left unchanged, it easily follows that the function $w$ constructed in Theorem \ref{th:berti:bolle} is invariant, too.
\end{remark}

\begin{theorem}\label{th:2.1}
Assume that the functions $h_{ij}=h_{ji}$'s are smooth, bounded, and (h1) holds
Then $M=\mathbb{R} \times M_0$ has at least $N$ non-trivial closed geodesics, distinct modulo the action of the group $O(2)$.
Furthermore, if
\begin{description}
\item[(h2)] the matrix $[h_{ij}(p,\cdot)]$ representing the bilinear form $h$ is  positive definite for $p \to + \infty$, and negative definite for $p \to - \infty$,
\end{description}
then $M$ possesses at least $2N$ non--trivial closed geodesics, geometrically distinct.
\end{theorem}

\begin{proof}
Observe that $Z=\mathbb{R} \times Z_0$, where  $Z_0=\{z_{p,q} \mid |p| = |q| =1, \; p \bullet q = 0\}$.
According to Theorem 2.1, it suffices to look for critical points of $\Phi_\ge$. From Lemma \ref{lem:due}, it follows that either $\Phi_\ge = b$ everywhere, or it has a critical point $(\bar{r},\bar{p},\bar{q})$. In any case such a critical point gives rise to a (non--trivial) closed geodesic of $(M,g_\ge)$.

From Remark \ref{rem:3}, we know that $\Phi_\ge$ is $O(2)$--invariant. This allows us to introduce the $O(2)$--category $\operatorname{cat}_{O(2)}$. One has
\begin{equation*}
\operatorname{cat}_{O(2)} (Z) \geq \operatorname{cat}(Z/O(2)) \geq \operatorname{cuplength} (Z/O(2)) +1.
\end{equation*}
Since $\operatorname{cuplength}(Z/O(2)) \geq N-1$, (see \cite{sch}), then $\operatorname{cat}_{O(2)} (Z) \geq N$. Finally, by the Lusternik--Schnirel'man theory, $M$ carries at least $N$ closed geodesics, distinct modulo the action $O(2)$. This proves the first statement.

Next, let
\begin{equation}
\Gamma (r,p,q)=G((r,z_{p,q}))=\frac{1}{2}\int_0^1 h(r,z_{p,q}(t))[{\dot z}_{p,q},{\dot z}_{p,q}] \, dt
\end{equation}
Then (h) immediately implies that
\begin{equation}
\Gamma (r,p,q) \to 0 \mbox{\quad as } |r| \to \infty,
\end{equation}
Moreover, if (h2) holds, then $\Gamma (r,p,q) > 0$ for $r > r_0$, and $\Gamma (r,p,q) < 0$ for $r < -r_0$. Since (recall equation \eqref{eq:gamma})
\begin{equation} \label{phi:gamma}
\Phi_\ge (r,p,q)=b+\ge \Gamma (r,p,q) + o(\ge),
\end{equation}
it follows that
\begin{equation*}
\begin{cases}
\Phi_\ge (r,p,q) > b &\text{ for } r > r_0 \\
\Phi_\ge (r,p,q) < b &\text{ for } r < - r_0.
\end{cases}
\end{equation*}
We can now exploit again the $O(2)$ invariance.

By assumption, and a simple continuity argument, $\{ \Phi_\ge > b \} \supset [R_0 , \infty ) \times Z_0$, and similarly $\{ \Phi_\ge < b \} \supset [-\infty , -R_0 ) \times Z_0$, for a suitably large $R_0 > 0$. Hence $\operatorname{cat}_{O(2)} ( \{ \Phi_\ge > b \}) \geq \operatorname{cat}_{O(2)}(Z_0) = N$. The same argument applies to $\{\Phi_\ge < b\}$. This proves the existence of at least $2N$ closed geodesics.
\end{proof}

\begin{remark}
\begin{description}
\item[(i)] In \cite{carminati}, the existence of $N$ closed geodesics on $S^N$ endowed with a metric close to the standard one is proved. Such a result does not need any study of $\Phi_\ge$ and its behavior. The existence of $2N$ geodesics is, as far as we know, new. We emphasize that it strongly depends on the {\em form} of $M=\mathbb{R} \times M_0$.

\item[(ii)] In \cite{tanaka}, the metric $g$ on $M$ is possibly not perturbative. No multiplicity result is given.
\end{description}
\end{remark}

\section{The general case}

In this section we consider a compact riemannian manifold $(M_0,g_0)$,
and in analogy to the previous section, we put
\begin{equation}
g_\ge (s,\xi)((r,\eta),(r,\eta))=|r|^2+g_0(\xi)(\eta ,\eta)+\ge h(s,\xi)((r,\eta),(r,\eta)).
\end{equation}
Again, we define $\Lambda = \{u=(r,x) \mid r \in H^1(S^1,\mathbb{R}), \; x \in H^1(S^1,M_0) \}$,
\[
E_{M_0}(x)=\frac{1}{2} \int_0^1 g_0(x)(\dot x , \dot x) \, dt,\quad E_0(r,x)=\frac{1}{2} \int_0^1 | \dot r|^2 \, dt + E_{M_0}(x),
\]
and finally
\[
E_\ge (r,x)=E_0(r,x)+\ge G(r,x),
\]
with $G$ as in \eqref{eq:G}.
It is well known (\cite{kli2}) that $M_0$ has a closed geodesic $z_0$. The functional $E_{M_0}$ has again a critical manifold $Z$ given by
\begin{equation*}
Z = \{ u(\cdot)=(\rho,z_0(\cdot + \tau)) \mid \rho \mbox{\; constant, } \tau \in S^1 \}.
\end{equation*}
Let $Z_0 = \{ z_0 ( \cdot + \tau ) \mid \tau \in S^1 \}$. It follows that $Z \approx \mathbb{R} \times Z_0$. The counterpart of $\Gamma$ in \eqref{eq:gamma} is
\begin{equation}\label{eq:gammagenerale}
\Gamma (r,\tau)=\frac{1}{2} \int_0^1 h(r,z_\tau)[\dot z_\tau, \dot z_\tau] \, dt.
\end{equation}

Let us recall some facts from \cite{kli1}.

\begin{remark}
There is a linear operator $A_z \colon T_z \Lambda (M_0) \to T_z \Lambda$, which is a compact perturbation of the identity, such that
\[
E_{M_0}^{\prime \prime} (z) [h,k]=\langle A_z h \mid k \rangle_1 = \int_0^1 \overbrace{A_z h}^{\cdot} \bullet \dot k \, dt.
\]
In particular, $E_0$ satisfies (AS2).
\end{remark}

\begin{definition} \label{defi:nondegen}
Let
\[
\ker E_{M_0}^{\prime \prime}(z_0)=\{ h \in T_{z_0} \Lambda (M_0) \mid \langle A_{z_0} h \mid k \rangle_1 =0 \quad \forall k \in T_{z_0} \Lambda (M_0) \}.
\]
We say that a closed geodesic $z_0$ of $M_0$ is {\em non--degenerate}, if
\[
\dim \ker E_{M_0}^{\prime \prime}(z_0)=1.
\]
\end{definition}

\begin{remark}
For example, it is known that when $M_0$ has negative sectional curvature, then all the geodesics of $M_0$ are non--degenerate. See \cite{DNF}.
Moreover, it is easy to see that the existence of non--degenerate closed geodesics is a {\em generic} property.
\end{remark}

\begin{lemma} \label{lem:3.1}
If $z_0$ is a non--degenerate closed geodesic of $M_0$, then $Z$ satisfies (AS2).
\end{lemma}

\begin{proof}
It is always true that $T_{z_r} Z \subset \ker E_0^{\prime \prime}(z_r)$. By \eqref{eq:nucleoE0},  we have that  $\dim T_{z_r} Z = \dim \ker E_0^{\prime \prime}(z_r)$. This implies that $T_{z_r} Z = \ker E_0^{\prime \prime}(z_r)$. 
A generic element of $Z$ has the form $(\rho,z^\tau)$ for $\rho \in \mathbb{R}$ and $z^\tau = z(\cdot + \tau)$; then
\[
T_{(\rho ,z^\tau)} M=\mathbb{R} \times T_{z^\tau} M_0,
\]
and any two vector fields $Y$ and $W$ along a curve on $M=\mathbb{R} \times M_0$ can be decomposed into
\begin{equation} \label{eq:Y}
Y=h(t)\frac{d}{dt} + y(t) \in \mathbb{R} \oplus T_{z^\tau} Z_0,
\end{equation}
\begin{equation} \label{eq:W}
W=k(t)\frac{d}{dt} + w(t) \in \mathbb{R} \oplus T_{z^\tau} Z_0.
\end{equation}

In addition, there results (see \cite{hebey})
\begin{equation}\label{eq:14}
E_{M_0}^{\prime \prime} (z_0)[y,w]=\int_0^1 \left[ g_0 (D_t y,D_t w)-g_0(R_{M_0}y(t),\dot z_0 (t))\dot z_0(t) \mid w(t)) \right] \, dt,
\end{equation}
and
\begin{equation}
R_M(r,z)=R_{\mathbb{R}} (r) + R_{M_0} (z) = R_{M_0}(z),
\end{equation}
where $R_M$, $R_{M_0}$, etc.  stand for the curvature tensors of $M$, $M_0$, etc. By \eqref{eq:14}, \eqref{eq:Y} and \eqref{eq:W}, as in the previous section, $E_0^{\prime \prime}(\rho ,z_\tau)[Y,W]=0$ is equivalent to the system
\begin{equation}\label{eq:nondegen}
\begin{cases}
\ddot h =0 \\
\int_0^1 g_0 (z)[D_t y,D_t w]- \langle R_{M_0}(y(t) , \dot z_r(t))\dot z_r(t) \mid w(t) \rangle  \, dt =0.
\end{cases}
\end{equation}
As in the case of the sphere, the first equation implies that $h$ is constant. The second equation in \eqref{eq:nondegen} implies that $y \in \ker E_{M_0}^{\prime \prime} (z^\tau)=\ker E_{M_0}^{\prime \prime}(z_0)$.  Hence,
\begin{equation} \label{eq:nucleoE0}
\ker E_0^{\prime \prime}(z_r)=\{ (h,y) \mid h \mbox{\; is constant, and } y \in \ker E_{M_0}^{\prime \prime} (z_0) \}.
\end{equation}
This completes the proof.
\end{proof}

\vspace{7pt}

\begin{theorem} \label{th:3.1}
Let $M \sb 0$ be a compact, connected manifold of dimension $N < \infty$. Assume that $M_0$ admits a non--degenerate closed geodesic $z$, and that (in local coordinates) $h_{ij}(p,\cdot) \to a_{-}$ as $p \to - \infty$, and $h_{ij}(p,\cdot) \to a_{+}$ as $p \to + \infty$.
\begin{enumerate}
\item If $a_{-}=a_{+}$ and $h_{ij}(p, \cdot)$ satisfies (h2), then $M$ has at least one closed geodesic.
\item  If $a_{-} \leq a_{+}$ and $h_{ij}(p,\cdot)[u,v]-a_{-}(u \mid v )$ is negative definite for $p \to - \infty$ and $h_{ij}(p,\cdot)[u,v]-a_{+}(u \mid v )$ is positive definite for $p \to + \infty$, then $M$ has at least two non-trivial closed geodesic.
\end{enumerate}
\end{theorem}

\begin{proof}
Lemma \ref{lem:3.1} allows us to repeat all the argument in Theorem \ref{th:2.1}, and the result follows immediately.
\end{proof}

\section{Isolated geodesics}

In this final section, we discuss one situaion where the critical manifold $Z$ may be degenerate. 
Here, the non--degeneracy condition (AS3) fails, and $T_z Z \subset \ker E_0^{\prime \prime}(z)$ strictly. Fix a closed geodesic $Z_0$ for $M_0$, and put $\tilde{W}=(T_{z_0} Z)^\bot$. Since $T_z Z \subset \ker E_0^{\prime \prime}(z)$ strictly, there exists $k > 0$ such that $\tilde{W}=(\ker E_0^{\prime \prime}(z_0))^\bot \oplus \mathbb{R}^k$. Repeating the preceding finite dimensional reduction, one can find again a unique map $\tilde{w}=\tilde{w}(z,\zeta)$, where $z \in Z$ and $\zeta \in \mathbb{R}^k$, in such a way that $E_\ge^\prime =0$ reduces to an equation like
\[
\nabla A(z+\zeta + \tilde{w}(z,\zeta))=0.
\]
If $z_0$ is an isolated minimum of the energy $E_{M_0}$ over some connected component of $\Lambda (M_0)$, then it is possible to show that there exists again a function $\Gamma \colon Z \to \mathbb{R}$ such that
\[
\nabla A(z+\zeta + \tilde{w}(z,\zeta))=0 \iff \frac{\de\Gamma}{\de r}(-R,\tau)\frac{\de \Gamma}{\de r}(R,\tau) \ne 0
\]
for some $R \in \mathbb{R}$ and all $\tau \in S^1$. For more details, see \cite{berti}. In particular, we will use the following result.

\begin{theorem} \label{th:berti}
Let $H$ be a real Hilbert space, $f_\ge \colon H \to \mathbb{R}$ is a family of $C^2$--functionals of the form $f_\ge = f_0 + \ge G$, and that:
\begin{description}
\item[(f0)] $f_0$ has a finite dimensional manifold $Z$ of critical points, each of them being a minimum of $f_0$;
\item[(f1)] for all $z \in Z$, $f_0^{\prime \prime}(z)$ is a compact perturbation of the identity.
\end{description}
Fix $z_0 \in Z$, put $W=(T_{z_0}Z)^\bot$, and suppose that $(f_0)_{|W}$ has an isolated minimum at $z_0$. Then, for $\ge$ sufficiently small, $f_\ge$ has a critical point, provided $\deg (\Gamma^\prime , B_R , 0 ) \ne 0$.
\end{theorem}

\begin{remark}
Theorem \ref{th:berti} has been presented in a linear setting. For Riemannian manifold, we can either 
reduce to a {\em local} situation and then apply the exponential map, or directly resort to the slightly 
more general degree theory on Banach manifold developed in \cite{ET}.
\end{remark}

\begin{theorem}\label{th:4.1}
Assume that $\pi_1 (M_0) \ne \{0\}$, and that all the critical points of $E_0$, the energy functional of $M_0$, are isolated. Suppose the bilinear form $h$ satisfies (h1), and
\begin{description}
\item[(h3)] $\dfrac{\de h}{\de r}(R,\xi)\dfrac{\de h}{\de r}(-R,\xi) \ne 0$ for some $R>0$ and all $\xi \in S^1$.
\end{description}
Then, for $\ge >0$ sufficiently small, the manifold $M=\mathbb{R} \times M_0$ carries at least one closed geodesic.
\end{theorem}

\begin{proof}
We wish to use Theorem \ref{th:berti}. Since $\pi_1 (M_0) \ne \{0\}$, then $E_0$ has a geodesic $z_0$ such that 
$E_0 (z_0) = \min E_0$ over some component $C$ of $\Lambda(M_0)$. See \cite{kli2}.

We consider the manifold
\[
Z = \{ u \in \Lambda \mid u(t)=(\rho,z_0(t+\tau)), \; \rho \mbox{\; constant, } \tau \in S^1 \}.
\]
Here we do not know, a priori, if $Z$ is non--degenerate in the sense of condition (AS2).
But of course $(E_0)_{W}$ has a minimum at the point $(\rho , z_0)$, where
$W=(T_{\rho , z^\tau)}Z)^\bot$. We now check that it 
is isolated for $(E_0)_{W}$. We still know that $Z = \mathbb{R} \times Z_0$. Take any point $(\rho , z_\tau) \in Z$, and observe that 
$T_{(\rho,z^\tau)} Z = \{ (r,y) \mid r \in \mathbb{R}, \; y \in T_{z^\tau}
Z_0 \}$. 
For all $(r,y)\in W$ sufficiently close to $(\rho,z_0)$, it holds in particular that $y \bot z_\tau$. Hence
\[
E_0(r,y) = L_0(r)+E_{M_0}(y) \geq E_{M_0}(y) >
E_{M_0}(z^\tau)=E_{M_0}(z_0) = E_0(\rho , z_0)
\]
since $L_0 \geq 0$ and $z_0$ (and hence $z^\tau$) is a minimum of $E_{M_0}$ by assumption.

Finally, thanks to assumption (h3), $\frac{\de \Gamma}{\de r}(-R,\tau)\frac{\de\Gamma}{\de r}(R,\tau) \ne 0$.

This concludes the proof.
\end{proof}


\chapter{Semilinear Schr\"{o}dinger equations}


In this chapter, following closely \cite{ams}, we present another application of the perturbation method to the existence of multiple positive solutions for a class
of
nonlinear Schr\"{o}dinger Equation (NLS in short)
\begin{equation}\tag{NLS}\label{eq:00}
\left\{
\begin{array}{lll} -\e^2 \Delta u + u +V(x)u=K(x)u^p, && \\
u \in W^{1,2}(\Rn),\,u>0 &&
\end{array}
\right.
\end{equation}
where $\Delta$ denotes the Laplace operator and
$$
1<p<\left\{ \begin{array}{lll} &\frac{n+2}{n-2},&\text{ if } \,
n\geq 3, \\
 &+\infty, & \text{ if } \, n=2.
\end{array}
\right.
$$

This equation has been studied under several viewpoints, and for a good review of
methods and results we refer to the monograph \cite{Chab} and its references. 

As a first step, we introduce some motivation for the study of the \emph{stationary Schr\"{o}dinger equation}, which will be the object of our interest in the next chapters. We shall assume some knowledge of linear functional analysis and operator theory. We refer to \cite{Weidmann} for details.

\begin{definition}
Given $V\in L^\infty (\Rn)$, we define the Schr\"{o}dinger operator $S\colon D(S)\subset L^2(\Rn)\to L^2(\Rn)$ generated by the potential $V$ by
\begin{equation*}
\begin{cases}
D(S)=H^2(\Rn)\\
Su=-\Delta u + Vu.
\end{cases}
\end{equation*}
\end{definition}

\begin{lemma}
\begin{enumerate}
\item $S$ is self--adjoint.
\item $\lambda \in \rho (S)$ if and only if $(a)$ $S_\lambda I$ is injective and $(b)$ $S-\lambda I)^{-1}\colon R(S-\lambda I)\subset L^2(\Rn)\to H^2(\Rn)$ is bounded.
\item $S$ has no eigenvalues if $V\equiv 0$.
\end{enumerate}
\end{lemma}

We now give some properties of the first eigenvalue of $S$. Let
\begin{equation}\label{eq:Lambda}
\Lambda = \inf \left\{ \int_\Rn (|\nabla u|^2 + Vu)\, dx \mid u\in H^1(\Rn) \mbox{\quad and }\int_\Rn |u|^2 =1 \right\} .
\end{equation}

\begin{definition}
Any minimizer of the right-hand side of \eqref{eq:Lambda} is called a \emph{ground state} of the Schr\"{o}dinger equation
\end{definition}

\begin{lemma}
Let $V\in L^\infty(\Rn)$. Then
\begin{enumerate}
\item $\Lambda \geq -\|V\|_\infty > -\infty$.
\item $\Lambda = \inf \left\{ \int_\Rn (|\nabla u|^2 + Vu)\, dx \mid u\in D(S) \mbox{\quad and }\int_\Rn |u|^2 =1 \right\}$ and so we also have
\[
\Lambda = \inf \left\{ \int_\Rn (Su)u\, dx \mid u\in H^1(\Rn) \mbox{\quad and }\int_\Rn |u|^2 =1 \right\} .
\]
\item If $u\in H^1(\Rn)$ with $\int_\Rn |u|^2\, dx =1$ and $\int_\Rn (|\nabla u^2 + Vu)\, dx=\Lambda$, then
\[
u\in H^2(\Rn),\, u\in \ker (S-\Lambda I) \text{\quad and } \Lambda\in \sigma_p(S).
\]
\end{enumerate}
\end{lemma}

\begin{lemma}[Properties of eigenfunctions]
\begin{enumerate}
\item Let $V\in L^\infty(\Rn)$ and consider $u\in\ker (S-\lambda I)$ for some $\lambda \in\Rn$. Then $u\in C(\Rn)\cap W^{1,s}(\Rn)\cap H^2(\Rn)$ for all $2\leq s\leq \infty$ and
\[
\lim_{|x|\to\infty} u(x)=0.
\]
\item Let $V\in L^\infty(\Rn)$ and choose $\xi <l= \lim_{R\to\infty} \essinf_{|x|\geq R} |V(x)|$. For any $\mu\in (0,\sqrt{l-\xi})$, there is a constant $C$, depending only on $\xi$ and $\mu$, such that
\[
|u(x)|\leq C \|u\|_\infty \, \E^{-\mu |x|}\mbox{\quad for all $x\in\Rn$}
\]
provided that $u\in\ker (S-\lambda I)$ for some $\lambda \leq \xi$.
\end{enumerate}
\end{lemma}

In the sequel we will always assume that $V,K\colon\Rn\to \R$ satisfy
\begin{description}
\item[$(V1)$] $V\in C^{2}(\Rn)$, $V$ and $D^{2}V$ are bounded;
\item[$(V2)$] $\inf_{x\in\Rn}(1+V(x))>0$;
\item[$(K1)$] $K\in C^{2}(\Rn)$, $K$ is bounded
and $K(x)>0$ for all $x\in \Rn$.
\end{description}
We seek solutions $u_\e$ of (NLS) that concentrate, as $\e \to 0$, near some
point $x_0\in \Rn$ (semiclassical standing waves).  By this we mean that
for all $x\in \Rn\setminus\{x_{0}\}$ one has that $u_\e(x)\to 0$ as $\e\to 0$.

 When $K$ equals a positive constant, say $K(x)\equiv 1$,
 (NLS) has been widely investigated, see
 \cite{ABC,ABer,dPF,FW,Li,Oh,W} and references therein.
 Moreover, the existence of multibump solutions has also been
studied in \cite{CN,G}, see also
\cite{ABer} where solutions with infinitely many bumps has been proved.
 It has been also pointed out,
 see e.g. \cite[Section 6]{ABC}, that the results contained in the
forementioned papers
 can be  extended to equations where $u^p$ is substituted by a
function $g(u)$, which behaves like $u^p$. Nonlinearities depending upon
$x$ have been
handled in \cite{Gr,WZ} where the existence of one-bump solutions is proved.

In a group of papers Equation (NLS) is studied
by perturbation arguments. For example, in \cite{ABC} a Liapunov-Schmidt
type procedure
is used to
reduce, for $\e$ small, (NLS) to a finite dimensional equation, that
inherits the
variational structure of the original problem. So, one looks for the
critical points of a finite
dimensional functional, which leading term is strictly related to the
behaviour of $V$ near its stationary points. For instance, by a direct application of the general theory described in Chapter 3, Ambrosetti, Badiale and Cingolani (\cite{ABC}) obtained the following result, in which one solution is found whenever $K\equiv 1$ and $V$ has a proper extremum at some point:

\begin{theorem}\label{th:ABC}
Let $K\equiv 1$, and assume that $V\colon \Rn \to\R$ is of class $C^2$ and possesses a non-degenerate maximum (or minimum) at a pont $x_0\in\Rn$. Then the equation
\[
-\Delta u + V(\e x)u=|u|^{p-1}u
\]
has at least a positive solution ocncentrating near $x_0$.
\end{theorem}

It should be clear that more than one solution can be found under the assumption that $V$ has well--separated proper extrema. One can indeed apply a rather standard gluing technique, similar to that appearing in \cite{ABer}.

As we shall see in the next sections, our approach is slightly different, and in some sense more general, since we can also find work in the ``non generic'' situation in which the critical points of $V$ arise in a compact manifold like a torus or a sphere. The price we pay is that the technical framework gets a bit more involved. As a counterpart, Theorem \ref{th:ABC} can be considered as a corollary or a special case of our main results.

A different approach has been carried out in another group of papers, like
\cite{dPF,W,WZ}. First, one finds a solution $u$ of (NLS) as a Mountain
Pass critical point of
the Euler functional;  after, one proves the concentration as $\e \to 0$.
Although this second approach allows to deal with potentials $V$ such that
$\sup V =+\infty$, it
works, roughly speaking, near minima of $V$, only. Actually, as pointed out
in \cite{ABC}, see also \cite[Theorem 3.2]{AmbBad1}, the Morse index of $u$ is
equal to
$1\,+$ the Morse index of the stationary point $V$  where concentration
takes place.
Hence $u$ can be found as a Mountain Pass only in the case of minima of $V$.
Such a severe restriction is even more apparent when one deals with
nonlinearities depending
upon $x$, see \cite{WZ}.  A case dealing with saddle-like potentials
is studied in \cite{dPFM}, but under some technical conditions
involving the saddle level and the supremum of $V$.

The second approach sketched above has been used in \cite{CL1} to obtain
the existence of multiple \emph{one--bump} solutions.
Letting $V_0=\inf [1+V(x)]$, it is assumed that
$K(x)\equiv 1$ and $\liminf_{|x|\to\infty}(1+V(x))>V_0 >0$. It has been
shown that, for $\e>0$
small,
(NLS) has at least $\cat(M,M_\d)$ solutions, where
 $M=\{x\in \Rn\mid 1+V(x)=V_0\}$, $M_\d$ is a $\d$-neighbourhood of $M$ and
 $\cat$ denotes the Lusternik-Schnirelman category. This result has been extended
 in \cite{CL2} to  nonlinearities like $K(x)u^p + Q(x)u^q$, with
 $1<q<p$.

 In the the present paper we will improve the above multiplicity result.
 For example, when $K(x)\equiv 1$, we can obtain the results summarized below.

 \begin{description}
 \item[$(i)$]  suppose that $V$ has a manifold $M$ of stationary points,
and  that $M$
 is nondegenerate in the sense of Bott, see \cite{Bott}. We prove, see
Theorem \ref{th:main},
 that (NLS) has at least $l(M)$ critical points concentrating near points
of $M$. Here
 $l(M)$ denotes the cup long of $M$.  This kind of result is new because
 \cite{CL1,CL2} deal only with (absolute) minima.

 \item[$(ii)$] If the points of $M$ are local minima or maxima of $V$
 the above result can be sharpened because $M$ does not need to be a
 manifold and $l(M)$ can be subsituted by $\cat(M,M_\d)$, see Theorem
 \ref{th:CL}.  Unlike \cite{CL1,CL2}, we can handle sets of {\it local}
 minima.  Furthermore, we do not require any condition at infinity.
 The case of maxima is new.

 \item[$(iii)$] When $K(x)$ is not identically constant, the preceding result
 holds true provided that $V$ is substituted by
 $A=(1+V)^{\theta}K^{-2/(p-1)}$, where $\theta=-(n/2)+(p+1)/(p-1)$, see
 Theorem \ref{th:K}.
 \end{description}
A similar multiplicity result holds for problems involving more
general nonlinearities, as in \cite{Gr}, see Remark \ref{rem:CLW}.

Our approach is based on the perturbation technique introduced above. However, a straight application of the
 general arguments, would only provide the existence of one
solution of (NLS) because that procedure leads to finding critical points
of a finite dimensional functional which does not inherit the
topological features of $M$.  So, in order to find multiplicity
results like those described above, it is necessary to use a different
finite dimensional reduction.  With this new approach one looks for
critical points of a finite dimensional functional $\Phi_{\e}$ which
is defined in a $\delta$--neighbourhood $M_\delta$ of $M$ and is
close to $V$ (or $A$) in the $C^1$ norm. This readily permits us to prove the multiplicity results in the case of maxima or minima, using the Lusternik--Schnirelman category. In the general case, we can use Conley theoretical arguments, see \cite{Con}. Rouglhy, it turns out that $M_\delta$ is an isolating block for the flow of $\nabla \Phi_\ge$ and we can apply an abstract critical point result of \cite{Chang} to find $l(M)$ solutions.

\

{\bf Notation}
\begin{itemize}
\item $\Wn$ denotes the usual Sobolev space, endowed with the norm
$$
\|u\|^{2}=\int_{\Rn}|\nabla u|^{2}dx + \irn u^2 dx.
$$

        \item  $o_h (1)$ denotes a function that tends to $0$ as $h\to 0$.

        \item  $c, C, c_i$ etc. denote (possibly different) positive constants,
        independent of the parameters.
\end{itemize}

\section{Preliminaries}\label{sec:prel}
In this section we will collect some preliminary material which will be
used in the rest of the paper.
 In order to semplify the notation, we will discuss in detail
equation (NLS) when $K\equiv 1$.  The general case will be handled in
Theorem \ref{th:K} and requires some minor changes, only.

 Without loss of generality we can assume that $V(0)=0$.
Performing the change of variable $x\mapsto \e x$, equation $(NLS)$ becomes
\begin{equation}\tag{$P_{\e}$}
 - \Delta u + u +V(\e x)u=u^p.
\label{eq:P}
\end{equation}
Solutions of (\ref{eq:P}) are the critical points $u\in \Wn$ of
 \[
 f_\e (u)=f_0(u)+ \frac{1}{2}\int_{\Rn}V(\e x)u^2dx,
 \]
 where
\[
f_0(u)=\frac{1}{2}\|u\|^2 - \frac{1}{p+1}\int_{\Rn}u^{p+1}dx.
\]
The solutions of (\ref{eq:P}) will be found near a solution of
\begin{equation}
- \Delta u + u +V(\e \xi)u=u^p,
\label{eq:xi}
\end{equation}
for an appropriate choice of $\xi\in \Rn$.  The solutions of
(\ref{eq:xi}) are critical points
of the functional
\begin{equation}
F^{\e\xi}(u)=f_{0}(u)+\frac{1}{2}\,V(\e
\xi)\,\int_{\Rn}u^2dx-\frac{1}{p+1}\int_{\Rn}u^{p+1}dx
\label{eq:F}
\end{equation}
and can be found explicitly.  Let $U$ denote the unique, positive,
radial solution of
\begin{equation}\label{eq:unp}
 - \Delta u + u =u^p ,\qquad u \in W^{1,2}(\Rn).
\end{equation}
Then a straight calculation shows that $\a U(\b x)$ solves (\ref{eq:P})
whenever
$$
\b=\b(\e\xi)= [1+V(\e\xi)]^{1/2}\quad {\rm and}\quad \a = \a(\e\xi)=[\b
(\e\xi)]^{2/(p-1)}.
$$
We set
\begin{equation}
        z^{\e\xi}(x)=\a(\e\xi)U(\b(\e\xi)x)
\label{eq:zU}
\end{equation}
and
$$
Z^{\e}=\{z^{\e\xi}(x-\xi)\mid \xi\in \Rn\}.
$$
 When there is no possible misunderstanding, we will write $z$, resp. $Z$,
 instead of $z^{\e\xi}$, resp $Z^{\e}$.
 We will also use the notation $z_{\xi}$ to denote the function
 $z_{\xi}(x):=z^{\e\xi}(x-\xi)$.
   Obviously all the
 functions in $z_{\xi}\in Z$ are solutions of (\ref{eq:unp}) or, equivalently,
 critical points of $F^{\e\xi}$.
 For future references let us point out some estimates. First of all, we
evaluate:
\begin{multline*}
\partial_{\xi}z^{\e\xi}(x-\xi)= \partial_{\xi}
\left[\a(\e\xi)U(\b(\e\xi)(x-\xi))\right] = \\
=  \e\a'(\e\xi)U(\b(\e\xi)(x-\xi))+
 \e\a(\e\xi)\b'(\e\xi)U'(\b(\e\xi)(x-\xi))-\\ \a(\e\xi)U'(\b(\e\xi)(x-\xi)).
\end{multline*}
Recalling the
 definition of $\a$, $\b$ one finds:
 \begin{equation}
 \partial_{\xi}z^{\e\xi}(x-\xi)=-\partial_{x}z^{\e\xi}(x-\xi)+O(\e|\nabla
V(\e\xi)|).     
        \label{eq:partialz}
 \end{equation}

\noindent The next Lemma shows that $\nabla f_{\e}(z_\xi)$ is close to
zero when
$\e$ is small.

 \begin{lemma}\label{lem:1}
For all $\xi\in \Rn$ and all $\e>0$ small, one has that
$$
\|\nabla f_{\e}(z_{\xi})\|\leq C\left(\e |\nabla
V(\e\xi)|+\e^{2}\right),\quad C>0.
$$
 \end{lemma}
\begin{proof}
From
$$
        f_{\e}(u)  =  F^{\e\xi}(u)+\frac{1}{2}
        \int_{\Rn}\left[V(\e x)-V(\e\xi)\right]u^{2}dx
$$
and since $z_{\xi}$ is a critical point of $F^{\e\xi}$, one has
\begin{multline*}
(\nabla f_{\e}(z_{\xi})|v)       = (\nabla F^{\e\xi}(z_{\xi})|v)+
\int_{\Rn}\left[V(\e
x)-V(\e\xi)\right]z_{\xi} v\,dx  \\
= \int_{\Rn}\left[V(\e x)-V(\e\xi)\right]z_{\xi} v\,dx .
\end{multline*}
Using the H\"{o}lder inequality, one finds
\begin{equation}
|(\nabla f_{\e}(z_{\xi})|v)|^{2}\leq \|v\|^{2}_{L^{2}} \int_{\Rn}|V(\e
x)-V(\e\xi)|^{2}z_{\xi}^{2}dx.
\label{eq:1.1}
\end{equation}
From the assumption that $|D^{2}V(x)|\leq {\rm const.}$ one infers
\[
|V(\e x)-V(\e\xi)|\leq \e |\nabla V(\e\xi)|\cdot |x-\xi|+c_{1}\e^{2}|x-\xi|^{2}.
\]
This implies
\begin{multline}\label{eq:1.3}
 \int_{\Rn}|V(\e x)-V(\e\xi)|^{2}z_{\xi}^{2}dx\leq  \nonumber \\
 c_1 \e^{2}|\nabla V(\e\xi)|^{2}\int_{\Rn}|x-\xi|^{2}z^{2}(x-\xi)dx +
 c_{2}\e^{4}\int_{\Rn}|x-\xi|^{4}z^{2}(x-\xi)dx.
\end{multline}
Recalling (\ref{eq:zU}), a direct calculation yields
\begin{eqnarray*}
\int_{\Rn}|x-\xi|^{2}z^{2}(x-\xi)dx & = &
\a^{2}(\e\xi)\int_{\Rn}|y|^{2}U^{2}(\b(\e\xi) y)dy \\
& = & \a^{2} \b^{-n-2}\int_{\Rn}|y'|^{2}U^{2}(y')dy'\leq c_{3}.
\end{eqnarray*}
From this (and a similar calculation for for the last integral in the
above formula) one infers
\begin{equation}
\int_{\Rn}|V(\e x)-V(\e\xi)|^{2}z_{\xi}^{2}dx\leq c_{4}\e^{2}|\nabla
V(\e\xi)|^{2} + c_{5}\e^{4}.
\label{eq:1.4}
\end{equation}
Putting together (\ref{eq:1.4}) and (\ref{eq:1.1}), the Lemma follows.
 \end{proof}

\section{Invertibility of $D^{2}f_{\e}$ on $TZ^{\perp}$}\label{sec:inv}

In this section we will show that $D^{2}f_{\e}$ is invertible on
$TZ^{\perp}$.  This will be the main tool to perform the finite
dimensional reduction, carried out in Section \ref{sec:fdr}.

Let $L_{\e,\xi}:(T_{z_{\xi}}Z^{\e})^{\perp}\to
(T_{z_{\xi}}Z^{\e})^{\perp}$ denote the operator defined by setting
$(L_{\e,\xi}v|w)= D^{2}f_{\e}(z_{\xi})[v,w]$.  We want to show
\begin{lemma}\label{lem:inv}
Given $\overline{\xi}>0$ there exists $C>0$ such that for $\e$ small enough
one has that
\begin{equation}
        |(L_{\e,\xi}v|v)|\geq C \|v\|^{2},\qquad \forall\;|\xi|\leq
\overline{\xi},\;\forall\; v\in
        (T_{z_{\xi}}Z^{\e})^{\perp}.
\label{eq:inv}
\end{equation}
\end{lemma}
\begin{proof}
From (\ref{eq:partialz}) it follows that every element
$\zeta\in T_{z_{\xi}}Z$ can be written in the form $\zeta =
z_{\xi}-\partial_{x}z^{\e\xi}(x-\xi)+O(\e)$.  As a consequence,
\begin{description}
        \item[$(*)$]  it suffices to prove (\ref{eq:inv}) for all
$v\in {\rm span}\{z_{\xi},\phi\}$, where $\phi$ is orthogonal to
 ${\rm span}\{z_{\xi},\partial_{x}z^{\e\xi}(x-\xi)\}$.
\end{description}
Precisely we shall prove that there exist $C_{1},C_{2}>0$ such that
for all $\e>0$ small and all $|\xi|\leq |\overline{\xi}|$
one  has:
\begin{eqnarray}
        (L_{\e,\xi}z_{\xi}|z_{\xi})& \leq & - C_{1}< 0.
        \label{eq:neg} \\
(L_{\e,\xi}\phi|\phi)&\geq & C_{2} \|\phi\|^2.
        \label{eq:claim}
\end{eqnarray}
It is clear that the Lemma immediately follows from $(*)$,
(\ref{eq:neg}) and (\ref{eq:claim}).

\

\noindent {\bf Proof of} (\ref{eq:neg}). First, let us recall that, since
$z_{\xi}$ is a Mountain Pass critical point of $F$, then given
$\overline{\xi}$ there exists $c_0>0$ such that for all $\e>0$ small
and all $|\xi|\leq |\overline{\xi}|$ one finds:
\begin{equation}
         D^2 F^{\e\xi}(z_{\xi})[z_{\xi},z_{\xi}] < -c_0< 0.
        \label{eq:D2Fa}
\end{equation}
One has:
$$
        (L_{\e,\xi}z_{\xi}|z_{\xi})=D^2 F^{\e\xi}(z_{\xi})[z_{\xi},z_{\xi}]+
        \int_{\Rn}\left[V(\e x)-V(\e\xi)\right]z_{\xi}^2 dx.
$$
The last integral can be estimated as in (\ref{eq:1.4}) yielding
\begin{equation}
        (L_{\e,\xi}z_{\xi}|z_{\xi})\leq D^{2} F^{\e\xi}(z_{\xi})[z_{\xi},z_{\xi}]
        + c_1 \e |\nabla V(\e\xi)|+c_2 \e^2.
        \label{eq:D2zz}
\end{equation}
From (\ref{eq:D2Fa}) and (\ref{eq:D2zz}) it follows that
(\ref{eq:neg}) holds.

\

\noindent {\bf Proof of} (\ref{eq:claim}).  As before, the fact that
$z_{\xi}$ is a Mountain Pass critical point of $F$ implies that
\begin{equation}
 D^2 F^{\e\xi}(z_{\xi})[\phi,\phi]>c_1 \|\phi\|^2.
        \label{eq:D2Fb}
\end{equation}
Let $R\gg 1$ and consider a radial smooth function
$\chi_{1}:\R^{n}\to \R$ such that
\begin{equation}\tag{$\chi'$}\label{eq:c1}
\chi_{1}(x) = 1, \quad \hbox{ for } |x| \leq R; \qquad
\chi_{1}(x) = 0, \quad \hbox{ for } |x| \geq 2 R;
\end{equation}
\begin{equation}\tag{$\chi''$}\label{eq:c2}
|\n \chi_{1}(x)| \leq \frac{2}{R}, \quad \hbox{ for } R \leq |x| \leq 2 R.
\end{equation}
We also set $ \chi_{2}(x)=1-\chi_{1}(x)$.
Given $\phi$ let us consider the functions
$$
\phi_{i}(x)=\chi_{i}(x-\xi)\phi(x),\quad i=1,2.
$$
A straight computation yields:
$$
\irn \phi^2 = \irn \phi_1^2 + \irn \phi_2^2 + 2\irn \phi_{1} \,
\phi_{2},
$$
$$
\irn |\n \phi|^2 = \irn |\n \phi_1|^2 + \irn |\n \phi_2|^2 + 2\irn
\n\phi_{1} \cdot \n \phi_{2},
$$
and hence
$$
\| \phi \|^2 = \| \phi_1 \|^2 + \| \phi_2 \|^2+
2 \irn\left[ \phi_{1} \, \phi_{2}+\n\phi_{1} \cdot \n
\phi_{2}\right].
$$
Letting $I$ denote the last integral,
one immediately finds:
$$
I=\underbrace{\irn \chi_{1}\chi_{2}(\phi^{2}+|\n \phi|^{2})}_{I_{\phi}} +
\underbrace{\irn\phi^{2}\n\chi_{1}\cdot \n\chi_{2}}_{I'}+
\underbrace{\irn\phi_{1}\n\phi\cdot\n\chi_{2}+\phi_{2}\n
\phi\cdot\n\chi_{1}}_{I''}.
$$
Due to the definition of $\chi$, the two integrals $I'$ and $I''$
reduce to integrals from $R$ and $2R$, and thus they are $o_{R}(1)\|\phi\|^{2}$.
As a consequence we have that
\begin{equation}\label{eq:d}
\| \phi \|^2 = \| \phi_1 \|^2 + \| \phi_2 \|^2 + 2I_\phi + o_R(1)\| \phi \|^2,
\end{equation}
After these preliminaries, let us evaluate the three terms in the
equation below:
$$
(L_{\e,\xi}\phi|\phi)=
\underbrace{(L_{\e,\xi}\phi_{1}|\phi_{1})}_{\a_{1}}+
\underbrace{(L_{\e,\xi}\phi_{2}|\phi_{2})}_{\a_{2}}+
2\underbrace{(L_{\e,\xi}\phi_{1}|\phi_{2})}_{\a_{3}}.
$$
One has:
\begin{equation} \label{eq:alfa1}
\a_{1}=(L_{\e,\xi}\phi_{1}|\phi_{1})=D^{2}F^{\e\xi}[\phi_{1},\phi_{1}]+\irn\left[V(\e x)-V(\e\xi)\right]\phi_{1}^{2}.
\end{equation}
In order to use (\ref{eq:D2Fb}), we introduce the function
$\overline{\phi}_{1}=\phi_{1}-\psi$, where
$$
\psi=(\phi_{1}|z_{\xi})z_{\xi}+(\phi_{1}|\partial_{x}z_{\xi})\partial_{x}z_{
\xi}.
$$
Then we have:
\begin{equation}
D^{2}F^{\e\xi}[\phi_{1},\phi_{1}]=D^{2}F^{\e\xi}[\overline{\phi}_{1},\overline{\phi}_{1}]+
D^{2}F^{\e\xi}[\psi,\psi]+2D^{2}F^{\e\xi}[\overline{\phi}_{1},\psi]
\label{eq:alfa2}
\end{equation}
 Let us explicitely point out that $\overline{\phi}_{1}\perp {\rm
span}\{z_{\xi},\partial_{x}z^{\e\xi}(x-\xi)\}$ and hence
(\ref{eq:D2Fb}) implies
\begin{equation}
D^{2}F^{\e\xi}[\overline{\phi}_{1},\overline{\phi}_{1}]\geq c_{1}
\|\overline{\phi}_{1}\|^{2}.
\label{eq:alfa3}
\end{equation}
On the other side, since $(\phi|z_{\xi})=0$ it follows:
\begin{eqnarray*}
        (\phi_{1}|z_{\xi}) & = & (\phi|z_{\xi})-(\phi_{2}|z_{\xi})=
-(\phi_{2}|z_{\xi}) \\
& = & -\irn\phi_{2}z_{\xi}dx-\irn \n z_{\xi}\cdot \n \phi_{2}dx \\
&=& -\irn\chi_{2}(y)z(y)\phi(y+\xi)dy-\irn \n z(y)\cdot \n
\chi_{2}(y)\phi(y+\xi)dy .
\end{eqnarray*}
Since $\chi_{2}(x)=0$ for all $|x|<R$, and since $z(x)\to 0$ as
$|x|=R\to\infty$, we infer
$(\phi_{1}|z_{\xi})=o_{R}(1)\|\phi\|$.  Similarly one shows that
$(\phi_{1}|\partial_{x}z_{\xi})=o_{R}(1)\|\phi\|$ and it follows that
\begin{equation}
         \|\psi\|=o_{R}(1)\|\phi\|.
\label{eq:alfa4}
\end{equation}
We are now in position to estimate the last two terms in Eq. (\ref{eq:alfa2}).
Actually, using (\ref{eq:alfa4}) we get
\begin{equation}
D^{2}F^{\e\xi}[\psi,\psi]=\|\psi\|^{2}+V(\e\xi)\irn\psi^{2}-p\irn
z_{\xi}^{p-1}\psi^{2}=o_{R}(1)\|\phi\|^{2}.     
\label{eq:alfa5}
\end{equation}
The same arguments readily imply
\begin{equation}
D^{2}F^{\e\xi}[\overline{\phi}_{1},\psi]=(\overline{\phi}_{1}|\psi)
+V(\e\xi)\irn \overline{\phi}_{1}\psi -p\irn
z_{\xi}^{p-1}\overline{\phi}_{1}\psi =o_{R}(1)\|\phi\|^{2}.
\label{eq:alfa6}
\end{equation}
Putting together (\ref{eq:alfa3}), (\ref{eq:alfa5}) and  (\ref{eq:alfa6})
we infer
\begin{equation}
D^{2}F^{\e\xi}[\phi_{1},\phi_{1}]\geq \|\phi_{1}\|^{2}+o_{R}(1)\|\phi\|^{2}.
\label{eq:alfa7}
\end{equation}
Using arguments already carried out before, one has
\begin{eqnarray*}
\irn|V(\e x)-V(\e\xi)|\phi_{1}^{2}dx & \leq
&c_{2}\irn|x-\xi|\chi^{2}(x-\xi)\phi^{2}(x)dx \\
& \leq & \e c_{3}\irn y\chi(y)\phi^{2}(y+\xi)dy \\
& \leq & \e c_{4}\|\phi\|^{2}.
\end{eqnarray*}
This and (\ref{eq:alfa7}) yield
\begin{equation}
        \a_{1}=(L_{\e,\xi}\phi_{1}|\phi_{1})\geq c_{5}\|\phi_{1}\|^{2}-
        \e c_{4}\|\phi\|^{2}+o_{R}(1)\|\phi\|^{2}.
\label{eq:alfa8}
\end{equation}
Let us now estimate $\a_{2}$.  One finds
$$
\a_{2}=(L_{\e,\xi}\phi_{2}|\phi_{2})=\irn |\n \phi_{2}|^{2}+
\irn (1+V(\e x))\phi_{2}^{2} -p\irn z_{\xi}^{p-1}\phi_{2}^{2}
$$
and therefore, using $(V2)$,
$$
\a_{2}  \geq  c_{6} \|\phi_{2}\|^{2} -p\irn z_{\xi}^{p-1}\phi_{2}^{2}.
 $$
As before, $\phi_{2}(x)=0$ for all $|x|<R$ and $z(x)\to 0$ as
$|x|=R\to\infty$ imply that
\begin{equation}
\a_{2}  \geq  c_{6} \|\phi_{2}\|^{2}+o_{R}(1)\|\phi\|^{2}.
\label{eq:alfa9}
\end{equation}
In a quite similar way one shows that
\begin{equation}
\a_{3} \geq c_{7}I_{\phi}+o_{R}(1)\|\phi\|^{2}.
\label{eq:alfa10}
\end{equation}
Finally, (\ref{eq:alfa8}), (\ref{eq:alfa9}), (\ref{eq:alfa10}) and the fact
that $I_{\phi}\geq 0$, yield
\begin{eqnarray*}
(L_{\e,\xi}\phi|\phi) & = & \a_{1}+\a_{2}+2\a_{3}\\
& \geq & c_{8}\left[\|\phi_{1}\|^{2}+\|\phi_{2}\|^{2}
+2 I_{\phi}\right]-c_{9}\e\|\phi\|^{2}+ o_{R}(1)\|\phi\|^{2}.
\end{eqnarray*}
Recalling (\ref{eq:d}) we infer that
$$
(L_{\e,\xi}\phi|\phi)\geq c_{10}\|\phi\|^{2}-c_{9}\e\|\phi\|^{2}+
o_{R}(1)\|\phi\|^{2}.
$$
Taking $\e$ small and $R$ large, eq.  (\ref{eq:claim}) follows.
This completes the proof of Lemma \ref{lem:inv}.
\end{proof}

\section{The finite dimensional reduction}\label{sec:fdr}

 In this Section we will show that the existence of critical points of $f_{\e}$
 can be reduced to the search of critical points of an auxiliary finite
 dimensional functional.  The proof will be carried out in 2
 subsections dealing, respectively, with  a Liapunov-Schmidt reduction, and
 with the behaviour of the auxiliary finite dimensional functional.
In a final subsection we handle the general case in which $K$ is not
identically constant.

\section{A Liapunov-Schmidt type reduction}\label{subsec:LS} The main
result of this section is the following lemma.

\begin{lemma}\label{lem:w}
For $\e>0$ small and $|\xi|\leq \overline{\xi}$ there exists a unique
$w=w(\e,\xi)\in
(T_{z_\xi} Z)^{\perp}$ such that
$\nabla f_\e (z_\xi + w)\in T_{z_\xi} Z$.
Such a $w(\e,\xi)$
 is of class $C^{2}$, resp.  $C^{1,p-1}$, with respect to $\xi$, provided that
 $p\geq 2$, resp.  $1<p<2$.
Moreover, the functional $\Phi_\e (\xi)=f_\e (z_\xi+w(\e,\xi))$ has
the same regularity of $w$ and satisfies:
 $$
\left(\n \Phi_\e(\xi_0)=0\right) \Rightarrow \left( \n
f_\e\left(z_{\xi_0}+w(\e,\xi_0)\right)=0\right)
$$
\end{lemma}
\begin{proof}
Let $P=P_{\e\xi}$ denote the projection onto $(T_{z_\xi} Z)^\perp$. We want
to find a
solution $w\in (T_{z_\xi} Z)^{\perp}$ of the equation
$P\nabla f_\e(z_\xi +w)=0$.  One has that $\n f_\e(z+w)=
\n f_\e (z)+D^2 f_\e(z)[w]+R(z,w)$ with $\|R(z,w)\|=o(\|w\|)$, uniformly
with respect to
$z=z_{\xi}$, for  $|\xi|\leq \overline{\xi}$. Using the notation introduced
in the
preceding Section \ref{sec:inv}, we are led to the equation:
$$
L_{\e,\xi}w + P\n f_\e (z)+PR(z,w)=0.
$$
According to Lemma \ref{lem:inv}, this is equivalent to
$$
w = N_{\e,\xi}(w), \quad \mbox{where}\quad
N_{\e,\xi}(w)=-L_{\e,\xi}^{-1}\left( P\n f_\e (z)+PR(z,w)\right).
$$
From Lemma \ref{lem:1} it follows that

\begin{equation}        \label{eq:N}
        \|N_{\e,\xi}(w)\|\leq c_1 (\e|\n V(\e\xi)|+\e^2)+ o(\|w\|).
        \end{equation}
Then one readily checks that $N_{\e,\xi}$ is a contraction on some ball in
$(T_{z_\xi} Z)^{\perp}$
provided that $\e>0$ is small enough and $|\xi|\leq \overline{\xi}$.
Then there exists a unique $w$ such that $w=N_{\e,\xi}(w)$.  Let us
point out that we cannot use the Implicit Function Theorem to find
$w(\e,\xi)$, because the map $(\e,u)\mapsto P\n f_\e (u)$ fails to be
$C^2$.  However, fixed $\e>0$ small, we can apply the Implicit
Function Theorem to the map $(\xi,w)\mapsto P\n f_\e (z_\xi + w)$.
Then, in particular, the function $w(\e,\xi)$ turns out to be of class
$C^1$ with respect to $\xi$.  Finally, it is a standard argument, see
\cite{AmbBad1,ABC}, to check that the critical points of $\Phi_\e
(\xi)=f_\e (z+w)$ give rise to critical points of $f_\e$.
\end{proof}

\begin{remark}\label{rem:w}
 From (\ref{eq:N}) it immediately follows that:
\begin{equation}        \label{eq:w}
        \|w\|\leq C \left(\e |\n V(\e\xi)|+\e^2\right),
\end{equation}
where $C>0$.
 \end{remark}

\

\noindent  For future references, it is
convenient to estimate the derivative $\partial_\xi w$.

\begin{lemma}\label{lem:Dw}
One has that:
\begin{equation}        \label{eq:Dw}
        \|\partial_\xi w\|\leq c \left(\e |\n
V(\e\xi)|+O(\e^2)\right)^\gamma,\qquad c>0,\;
        \gamma=\min\{1,p-1\}.
\end{equation}
\end{lemma}
\begin{proof}
In the proof we will write $\dw_i$, resp $\dz_i$, to denote the components of
 $\partial_\xi w$, resp.
$\partial_\xi z$; moreover we will set $h(z,w)=(z+w)^p-z^p-pz^{p-1}w$.
With these notations, and recalling that
$L_{\e,\xi}w = -\Delta w +w +V(\e x)w -pz^{p-1}w$, it follows that
$w$ satisfies $\forall \;v\in (T_{z_\xi} Z)^{\perp}$:
\begin{multline*}
(w|v) + \irn V(\e x)wvdx- p\irn z^{p-1}wvdx \\
+\irn [V(\e x)-V(\e\xi)]zv dx-
\irn h(z,w)v dx=0.
\end{multline*}
Hence $\dw_i$ verifies
\begin{multline}
         (\dw_i |v) + \irn V(\e x)\dw_i vdx- p\irn z^{p-1}\dw_i vdx
-p(p-1)\irn z^{p-2}\dz_i wv dx  \\
         +\irn [V(\e x)-V(\e\xi)]\dz_i v dx-\e \partial_{x_i} V(\e\xi)\irn zv dx
          - \irn (h_z\dz_i +
h_w \dw_i )v dx=0.\label{eq:dw}
\end{multline}
Let us set $L'=L_{\e,\xi}-h_w $. Then
 (\ref{eq:dw}) can be written as
\begin{multline}\label{eq:Ltilde}
(L'\dw_i|v) =p(p-1)\irn z^{p-2}\dz_i wv
- \irn [V(\e x)-V(\e\xi)]\dz_i v\\
+\e \partial_{x_i} V(\e\xi) \irn zv + \irn  h_z\dz_i v ,
        \end{multline}
and hence one has
\begin{multline*}
|(L'\dw_i|v)| \leq c_1 \|w\|\cdot\|v\|+\left|\irn [V(\e x)-V(\e\xi)]\dz v
dx\right| \\
+c_2 \e |\n V(\e\xi)|\cdot\|v\|
+ \left|\irn h_z\dz v dx\right|.
\end{multline*}
As in the proof of Lemma \ref{lem:1} one gets
$$
\left|\irn [V(\e x)-V(\e\xi)]\dz v dx\right|\leq\left( c_3 \e |\n
V(\e\xi)|+c_4 \e^2\right)\,\|v\|.
$$
Furthermore, the definition of $h$ immediately yields
$$
\left|\irn h_z\dz v dx\right|\leq c_5 \|w\|^\gamma \|v\|,\quad \mbox{where}\quad
\gamma = \min\{1,p-1\}.
$$
Putting together the above estimates we find
$$
|(L'\dw_i|v)|\leq
\left[c_6 \e |\n V(\e\xi)|+c_4 \e^2 + c_6 \|w\|^\gamma\right]\,\|v\|
$$
Since $h_w\to 0$ as $w \to 0$, the operator
$L'$, likewise $L$, is  invertible for $\e >0$ small and therefore one finds
$$
\|\dw\|\leq \left( c_7 \e |\n V(\e\xi)|+c_8 \e^2\right) + c_9 \|w\|^\gamma,
$$
Finally, using (\ref{eq:w})
the Lemma follows.      
 \end{proof}

\section{The finite dimensional functional}\label{subs:3}
The main purpose of this subsection is to use the estimates
on $w$ and $\partial_\xi w$ estabilished above to find
an expansion of $\n \Phi_\e (\xi)$, where $\Phi_\e (\xi)= f_\e(z_\xi
+w(\e,\xi))$.  In
the sequel, to be short, we will often write $z$ instead of $z_\xi$
and $w$ instead of $w(\e,\xi)$.  It is always understood that
$\e$ is taken in such a way that all the results discussed in the
preceding Section hold.

For the reader's convenience we will divide the arguments in some steps.

\

\noindent {\sl Step 1.} We have:
$$
        \Phi_\e (\xi) =  \frac{1}{2}\|z+w\|^2 + \frac{1}{2}\irn V(\e x)(z+w)^2
        -\frac{1}{p+1}\irn (z+w)^{p+1}.
$$
Since $-\Delta z+z+V(\e\xi)z=z^p$ we infer that
\begin{eqnarray*}
        \|z\|^2 & = &  -V(\e\xi)\irn z^2 + \irn z^{p+1}, \\
        (z|w) & = & -V(\e\xi)\irn zw + \irn z^pw,
\end{eqnarray*}
Then we find :
\begin{eqnarray*}
         \Phi_\e (\xi)& = & \left(\frac{1}{2}-\frac{1}{p+1}\right) \irn z^{p+1} +
         \frac{1}{2}\int_{\Rn}\left[V(\e x)-V(\e\xi)\right]z^{2}\\
         &  & +\int_{\Rn}\left[V(\e x)-V(\e\xi)\right]zw+
         \frac{1}{2}\irn V(\e x)w^2 \\
&&+ \frac{1}{2} \|w\|^2 - \frac{1}{p+1}\irn
\left[(z+w)^{p+1}-z^{p+1}-(p+1)z^pw\right].
\end{eqnarray*}
Since $z(x)=\a(\e\xi)U(\b(\e\xi)x)$, where
$\a=(1+V)^{1/(p-1)}$ and $\b=(1+V)^{1/2}$, see (\ref{eq:zU}),
it follows
$$
\irn z^{p+1}dx= C_0 \left(1+V(\e\xi)\right)^\theta,\quad C_0=
 \irn U^{p+1};\qquad \theta=\frac{p+1}{p-1}-\frac{n}{2}.
$$
Letting $C_1= C_0 [1/2 -1/(p+1)]$ one has
\begin{multline}\label{eq:F1}
        \Phi_\e (\xi)  =  C_1 \left(1+V(\e\xi)\right)^\theta+
        \frac{1}{2}\int_{\Rn}\left[V(\e x)-V(\e\xi)\right]z^{2}+\\ 
        \int_{\Rn}\left[V(\e x)-V(\e\xi)\right]zw
         + \frac{1}{2}\irn V(\e x)w^2 \\ +\frac{1}{2} \|w\|^2 -
        \frac{1}{p+1} \irn \left[(z+w)^{p+1}-z^{p+1}-(p+1)z^pw\right].
\end{multline}

\

\noindent {\sl Step 2.}
Let us now evaluate the derivative of the right hand side
of (\ref{eq:F1}).
For this, let us write:

\begin{equation}        \label{eq:D}
\Phi_\e (\xi)=C_1(1+V(\e\xi))^{\theta}+\Lambda_\e(\xi)+
\Psi_\e(\xi),
        \end{equation}
where
$$
\Lambda_\e(\xi)=\frac{1}{2}\int_{\Rn}\left[V(\e x)-V(\e\xi)\right]z^2+
\int_{\Rn}\left[V(\e x)-V(\e\xi)\right]zw
$$
and
$$
\Psi_\e(\xi)=
        \frac{1}{2}\irn V(\e x)w^2 +\frac{1}{2} \|w\|^2 -
         \irn \left[(z+w)^{p+1}-z^{p+1}-(p+1)z^pw\right].
$$
It is easy to check, by a direct calculation, that
\begin{equation}        \label{eq:D1}
        |\n \Psi_\e(\xi)|\leq c_1\|w\|\,\|\dw\|.
        \end{equation}
Furthermore,
since $V(\e x)-V(\e\xi)=\e\n V(\e\xi)\cdot(x-\xi)+\e^2 D^2V(\e\xi+\tau
\e(x-\xi))[x-\xi,x-\xi]$ for some $\tau\in [0,1]$ we get:
\begin{eqnarray*}
        \int_{\Rn}\left[V(\e x)-V(\e\xi)\right]z^2dx & = & \e \irn \n
V(\e\xi)\cdot(x-\xi)z^2dx\\
&&\quad +\e^2 \irn D^2V(\e\xi+\tau
\e(x-\xi))[x-\xi,x-\xi]z^2 dx  \\
&=&  \e \irn \n V(\e\xi)\cdot y\, z^2(y)dy\\
         &  &\quad +\e^2 \irn D^2V(\e\xi+\tau\e y)[y,y]z^2(y) dy \\
&=& \e^2 \irn D^2V(\e\xi+\tau\e y)[y,y]z^2(y) dy.
\end{eqnarray*}
Similarly, from $V(\e x)-V(\e\xi)=\e\n V(\e\xi+\tau
\e(x-\xi))\cdot (x-\xi)$ one finds
$$
\int_{\Rn}\left[V(\e x)-V(\e\xi)\right]zw=\e \irn \n V(\e\xi+\tau
\e(x-\xi))\cdot (x-\xi) zw
$$
The preceding equations imply:
\begin{equation}
        |\n \Lambda_\e(\xi)|
        \leq c_2 \e^2 + c_3 \e\|w\|.
        \label{eq:D2}
\end{equation}
Taking the gradients in (\ref{eq:D}), using (\ref{eq:D1}-\ref{eq:D2}) and
recalling
the estimates (\ref{eq:w}) and (\ref{eq:Dw}) on $w$ and $\dw$, respectively,
we readily find:
\begin{lemma}\label{lem:DF}
Let $a(\e\xi)=\theta C_1(1+V(\e\xi))^{\theta -1}$. Then one has:
$$
\n \Phi_\e (\xi)= a(\e\xi) \e\n V(\e\xi) + \e^{1+\gamma} R_\e(\xi),
$$
where $|R_\e(\xi)|\leq const$ and $\gamma=\min\{1,p-1\}$.
\end{lemma}

\begin{remark}\label{rem:Phi}
Using similar arguments one can show that
 $$
 \Phi_\e (\xi)= C_1 \left(1+V(\e\xi)\right)^\theta + \rho_\e(\xi),\quad
C_1>0,\quad
 \theta=\frac{p+1}{p-1}-\frac{n}{2},
 $$
 where $|\rho_\e(\xi)|\leq \operatorname{const.} \left(\e |\nabla V(\e\xi)| +\e^2\right)$.
\end{remark}

\section{The general case}\label{subsec:gen}
Let us indicate the counterpart of the above results in the case in which
$K$ is not constant.  Since the arguments are quite similar, we
will only  outline the main modifications that are needed.

\

After rescaling, the solutions of (NLS) are the critical points of
$$
\widetilde{f}_\e(u)=\frac{1}{2}\|u\|^2 +\frac{1}{2}\int_{\Rn}V(\e x)u^2-
\frac{1}{p+1}\int_{\Rn}K(\e x)u^{p+1}.
$$
The solutions of (NLS) will be found near  solutions of
\begin{equation}
- \Delta u + u +V(\e \xi)u=K(\e \xi)u^p,
\label{eq:xiK}
\end{equation}
namely near critical points of
$$
\widetilde{F}^{\e\xi}(u)=\frac{1}{2}\|u\|^2
+\frac{1}{2}\,V(\e \xi)\,\int_{\Rn}u^2 -
\frac{1}{p+1}\,K(\e \xi)\int_{\Rn}u^{p+1}.
$$
If $\wz=\wz^{\e\xi}$ is a solution of (\ref{eq:xiK}), then $\wz (x)=
\tilde{\a}(\e\xi)U(\tilde{\b}(\e\xi)x)$, where
$$
\tilde{\a}(\e\xi)=\left(1+V(\e\xi)\right)^{1/2},\quad \tilde{\b}(\e\xi)=
\left(\frac{1+V(\e\xi)}{K(\e\xi)}\right)^{1/(p-1)}.
$$
This implies that

\begin{equation}        \label{eq:wzA}
\irn \wz ^{p+1}=C_0 A(\e\xi),
        \end{equation}
where
$$
A(x)=\left(1+V(x)\right)^{\theta}[K(x)]^{-2/(p-1)},\quad
\theta=\frac{p+1}{p-1}-\frac{n}{2}.
$$
Define $\widetilde{L}=\widetilde{L}_{\e,\xi}$ on $(T\widetilde{Z})^{\perp}$
by setting $(\widetilde{L} v|w)=D^2\widetilde{f}_\e(\wz_\xi)[v,w]$. As in Lemma
\ref{lem:inv} one shows that $\widetilde{L}$ is invertible for $\e$ small.
Furthermore, one has
$$
\widetilde{f}_\e(u)=\widetilde{F}^{\e\xi}(u)+
\frac{1}{2}\irn [V(\e x)-V(\e \xi)]u^2 -
\frac{1}{p+1}\irn[K(\e x)-K(\e \xi)]u^{p+1}
$$
and, as in Lemma \ref{lem:1} one finds that
$$
\|\widetilde{f}_\e(u)\|\leq c_1 \left(\e(|\n V(\e \xi)|+|\n
K(\e\xi)|)+\e^2\right).
$$
This and the invertibility of $\widetilde{L}$ imply,
as in  Lemma \ref{lem:w}, the existence of
$\widetilde{w}(\e,\xi)$ such that the critical points of the finite
dimensional functional
$$
\widetilde{\Phi}_\e(\xi)=\widetilde{f}_\e(\wz_\xi+\widetilde{w}(\e,\xi))
$$
give rise to critical points of $\widetilde{f}_\e(u)$.
Such a $\widetilde{w}$ is $C^1$ with respect
to $\xi$ and, as in Lemma \ref{lem:Dw}, the following estimate holds
$$
\|\partial_\xi\widetilde{w}\|\leq c_2
\left(\e(|\n V(\e \xi)|+|\n K(\e\xi)|)+O(\e^2)\right)^\gamma.
$$
It remains to study the finite dimensional functional
$\widetilde{\Phi}_\e (\xi)$,
$$
\widetilde{\Phi}_\e(\xi)=\frac{1}{2}\|\wz+\widetilde{w}\|^2 +
\frac{1}{2}\irn V(\e x)(\wz+\widetilde{w})^2
        -\frac{1}{p+1}\irn K(\e x)(\wz+\widetilde{w})^{p+1}.
$$
Since now $\|\wz\|^2=-V(\e\xi)\int \wz^2 + K(\e\xi)\int \wz^{p+1}$
and $(\wz |\widetilde{w})=-V(\e\xi)\int \wz\widetilde{w}+ K(\e\xi)\int \wz^{p}
\widetilde{w}$, one gets
\begin{eqnarray*}
         \widetilde{\Phi}_\e (\xi)& = & \irn\left(\frac{1}{2}K(\e x)-\frac{1}{p+1}
         K(\e\xi)\right) \wz^{p+1} +
         \frac{1}{2}\irn\left[V(\e x)-V(\e\xi)\right]\wz^{2}\\
         &  & +\int_{\Rn}\left[V(\e x)-V(\e\xi)\right]\wz \widetilde{w}+
         \frac{1}{2}\irn V(\e x)\widetilde{w}^2+ \frac{1}{2} \|w\|^2 \\
         && +K(\e\xi)\irn\wz^p \widetilde{w} -
         \frac{1}{p+1}\irn K(\e x) \left[(z+w)^{p+1}-z^{p+1}-(p+1)z^pw\right]
\end{eqnarray*}
From $K(\e x)=K(\e\xi)+\e\n K(\e\xi)\cdot (x-\xi)+O(\e^2)$ and since
$\int \n K(\e\xi)\cdot y \wz^{p+1}=0$ we infer:
$$
\irn\left(\frac{1}{2}K(\e x)-\frac{1}{p+1}
         K(\e\xi)\right) \wz^{p+1}=
         (\frac{1}{2}-\frac{1}{p+1})\irn
         K(\e\xi) \wz^{p+1} + O(\e^2).
 $$
Using (\ref{eq:wzA}), one finds the counterpart of (\ref{eq:F1}), namely
\begin{multline*}
        \widetilde{\Phi}_\e (\xi)  =  C_1 A(\e\xi)+
        \frac{1}{2}\irn\left[V(\e x)-V(\e\xi)\right]\wz^{2}
          +\int_{\Rn}\left[V(\e x)-V(\e\xi)\right]\wz \widetilde{w}\\
        + \frac{1}{2}\irn V(\e x)\widetilde{w}^2
        + \frac{1}{2} \|w\|^2
 +K(\e\xi)\irn\wz^p \widetilde{w} \\- \frac{1}{p+1}\irn K(\e x)
 \left[(z+w)^{p+1}-z^{p+1}-(p+1)z^pw\right]+O(\e^2).
\end{multline*}
Taking the derivative of the above equation and using
the preceding estimates, one finally yields

\begin{equation}\label{eq:wDF}
\n \widetilde{\Phi}_\e (\xi)=C_1 \n A(\e\xi)+\e^{1+\gamma} \widetilde{R}_\e
(\xi),
\quad \gamma=\min\{1,p-1\},
 \end{equation} where $|\widetilde{R}_\e (\xi)|\leq \operatorname{const}$, which is
the counterpart of Lemma \ref{lem:DF}.  Let us also point out that,
like in Remark \ref{rem:Phi}, one has:
\begin{equation}
        \widetilde{\Phi}_\e (\xi)=C_1 A(\e\xi)+O(\e).
        \label{eq:A}
\end{equation}

\section{Main results}\label{sec:main}
In this section we will prove the main results of the present paper.
First, some preliminaries are in order.

Given a set $M\subset \Rn$, $M\ne \emptyset$, we denote by $M_\d$ its $\d$
neighbourhood.
If $M\subset N$, $cat(M,N)$ denotes the Lusternik-Schnirelman category
of $M$ with respect to $N$, namely the least integer $k$ such that
$M$ can be covered by $k$ closed subsets of
$N$, contractible to a point in $N$. We set $cat(M)=cat(M,M)$.

The {\sl cup long} $l(M)$ of $M$ is defined by
$$
l(M)=1+\sup\{k\in \N: \exists\, \a_{1},\ldots,\a_{k}\in
\check{H}^{*}(M)\setminus 1, \,\a_{1}\cup\ldots\cup\a_{k}\ne 0\}.
$$
If no such class exists, we set $l(M)=1$. Here $\check{H}^{*}(M)$ is the
Alexander
cohomology of $M$ with real coefficients and $\cup$ denotes the cup product. For example, if $M=S^{n-1}$, the $(n-1)$--dimensional sphere in $\Rn$, then $l (M)=\cat (M)=2$; however, in general, $l (M) \leq \cat (M)$.

Let us suppose that $V$ has a smooth manifold of critical points $M$. According to Bott \cite{Bott}, we say that $M$ is nondegenerate (for $V$) if every $x\in M$ is a nondegenerate critical point of $V_{|M^\perp}$. The Morse index of any $x\in M$, as a critical point of $V_{|M^\perp}$, is constant and is, by definition, the Morse index of $M$.

We first recall a result from \cite{Chang}, namely Theorem 6.4 in Chapter II, adapted our notatin and purposes. See also \cite{Con}.

\begin{theorem}\label{th:chang}
Let $h\in C^2(\Rn)$ and let $\Sigma\subset \Rn$ be a smooth compact nondegenerate manifod of critical points of $h$. Let $U$ be a neighborhood of $\Sigma$ and let $\ell\in C^1 (\Rn)$. Then, if $\|\ell - h\|_{C^1(U)}$ is sufficiently small, the function $\ell$ possesses at least $l(\Sigma)$ critical points in $U$.
\end{theorem}

Roughly, the previuos theorem fits into the frame of Conley theory and extends an older result by Conley and Zehnder \cite{CZ}.

We are now ready to state our multiplicity results.
Our first result deals with (NLS) with $K(x)\equiv 1$,
namely with the equation (\ref{eq:P}) introduced in Section
\ref{sec:prel}.  
\begin{theorem}\label{th:main}
Let  $(V1-2)$ hold and suppose $V$ has
a nondegenerate smooth manifold of critical points  $M$.
Then for $\e>0$ small, (\ref{eq:P}) has at least $l(M)$ solutions that
concentrate near
points of $M$.
\end{theorem}
\begin{proof}First of all, we fix $\overline{\xi}$ in such a way that
$|x|<\overline{\xi}$ for all $x\in M$.  We will apply the
finite dimensional procedure with such $\overline{\xi}$ fixed. Since $M$ is
a nondegenerate smooth manifold of critical points of $V$, it is obviuosly a nondegenerate manifold of critical points of $C_1 (1+V(\cdot))^\theta$ as well.

In order to use Theorem 20 we set $h(\xi)=C_1 (1+V(\xi))^\theta$, $\Sigma = M$ and $\ell (\xi)=\Phi_\e (\xi/\e)$. Fix a $\delta$-neighborhood $M_\delta$ of $M$ in such a way that $M_\delta\subset \{|x|<\bar\xi\}$ and the only critical points of $V$ in $M_\delta$ are those of $M$, and let $U=M_\delta$. From Lemma 19, the function $\Phi_\e (\cdot / \e)$ converges to $C_1 (1+V(\cdot))^\theta$ in $C^1 (U)$ when $\e\to 0$. Hence Theorem 20 applies and we can infer the existence of at least $l(M)$ critical points of $\ell$, provided $\e>0$ is small enough. Let $\xi_i\in M_\delta$ be any of those critical points. Then $\xi_i /\e$ is a critical point of $\Phi_\e$ and so $u_{\e,\xi_i}(x)=z^{\xi_i}(x-\xi_i / \e)+w(\e,\xi_i)$ is a critical point of $f_\e$ and hence a solution of ($P_\e$). It follows that
\[
u_{\e,\xi_i}(x/\e) \approx z^{\xi_i}\left( \frac{x-\xi_i}{\e}\right)
\]
is a solution of (NLS), with $K\equiv 1$. Any $\xi_i$ converges to some $\xi_i^*\in M_\delta$ as $\e\to 0$ and, using again Lemma 19, it is easy to see that $\xi_i^*\in M$. This shows that $u_{\e,\xi_i}(x/\e)$ concentrates near a point of $M$.
\end{proof}

The next result deals with the more general equation (NLS).
The results will be given using the auxiliary function
$$
A(x)=\left(1+V(x)\right)^{\theta}[K(x)]^{-2/(p-1)},\quad
\theta=\frac{p+1}{p-1}-\frac{n}{2},
$$
introduced in Subsection \ref{subsec:gen}.

\begin{theorem}\label{th:K}
Let  $(V1-2)$ and $(K1)$ hold and suppose $A$ has
a nondegenerate smooth manifold of critical points  $\widetilde{M}$.
Then  for $\e>0$ small, (NLS) has at least $l(\widetilde{M})$ solutions
that concentrate near
points of $\widetilde{M}$.
\end{theorem}
\begin{proof}
The proof is quite similar to the preceding one, by using the results
discussed in Subsection \ref{subsec:gen}.
\end{proof}
\begin{remark}\label{rem:CLW}
The preceding results can be extended to cover a class of nonlinearities
$g(x,u)$ satisfying the same assumptions of \cite{Gr}.  In addition to
some techenical conditions, one roughly assumes that the problem

\begin{equation}
        -\Delta u+u+V(\e\xi)u=g(\e\xi,u),\quad u>0,\quad u\in \Wn,
\label{eq:Gr}
\end{equation}
has a unique radial solution $z=z^{\e\xi}$ such that the linearized problem
at $z$ is invertible on $(T_{z}Z)^{\perp}$.  In such a case one can
obtain the same results as above with the auxiliary function $A$
substituted by
$$
{\cal A}(\e\xi)=\irn\left[\frac{1}{2}g(\e\xi,z^{\e\xi}(x-\xi))-
G(\e\xi,z^{\e\xi}(x-\xi))\right]dx,
$$
where $\partial_{u}G=g$.  When $g(x,u)=K(x)u^p$ one finds that
${\cal A}=A$.  But, unlike such a case, the function $\cal A$ cannot be written
in an explicit way.  This is the reason why
 we have focused our study to the model problem (NLS)
when the auxiliary function $A$ has an explicit and neat form.  Let us
also point out that the class of nonlinearities handled in \cite{CL2}
does not even require that (\ref{eq:Gr}) has a unique solution.
\end{remark}

\

\noindent When we deal with local minima (resp. maxima) of $V$, or $A$, the
preceding results
can be improved because the number of positive solutions of (NLS)
can be estimated by means of the category and $M$ does not need to be
a manifold.

\begin{theorem}\label{th:CL}
Let  $(V1-2)$ and $(K1)$ hold and suppose $A$ has
a compact set $X$ where $A$ achieves a strict local minimum, resp.
maximum.

Then  there exists $\e_{\d}>0$ such that (NLS) has at least $cat(X,X_\d)$
solutions that concentrate near points of $X_{\d}$, provided $\e\in
(0,\e_{\d})$.
\end{theorem}
\begin{proof}
Let $\d>0$ be such that
$$
b:=\inf\{A(x):x\in \partial X_{\d}\}>a:= A_{|X}
$$
and fix again
$\overline{\xi}$ in such a way that $X_{\d}$ is contained in
$\{x\in\Rn : |x|<\overline{\xi}\}$.  We set
$
X^{\e}=\{\xi:\e\xi\in X\}$, $X_{\d}^{\e}=\{\xi:\e\xi\in X_{\d}\}$ and
$Y^{\e}=\{\xi\in X_{\d}^{\e} :\widetilde{\Phi}_{\e}(\xi)\leq C_{1}(a+b)/2\}$.
From
(\ref{eq:A}) it follows that there exists $\e_{\d}>0$ such that
\begin{equation}
        X^{\e}\subset Y^{\e}\subset X^{\e}_{\d},
\label{eq:X}
\end{equation}
 provided $\e\in
(0,\e_{\d})$.  Moreover, if $\xi\in \partial X_{\d}^{\e}$ then
$V(\e\xi)\geq b$ and hence
$$
\widetilde{\Phi}_{\e}(\xi)\geq C_{1}V(\e\xi) + o_{\e}(1) \geq C_{1}b +
o_{\e}(1) .
$$
On the other side, if $\xi\in Y^{\e}$
then $\widetilde{\Phi}_{\e}(\xi)\leq C_{1}(a+b)/2$.  Hence, for
 $\e$ small, $Y^{\e}$ cannot meet $\partial X_{\d}^{\e}$
 and this readily implies that $Y^{\e}$ is compact.
Then $\widetilde{\Phi}_{\e}$ possesses at least $cat(Y^{\e},X^{\e}_{\d})$
critical points in $ X_{\d}$.  Using (\ref{eq:X}) and the properties of
the category one gets
$$
\cat(Y^{\e},Y^{\e})\geq \cat(X^{\e},X^{\e}_{\d})=\cat(X,X_{\d}),
$$
and the result follows.
\end{proof}

\

\begin{remark}
The approach carried out in the present paper
works also in the more standard case when $V$, or $A$, has an isolated
set $S$ of stationary points and $\deg (\n V,\Omega,0)\not= 0$ for some
neighbourhood $\Omega$ of $S$.  In this way we can, for example, recover the
results of \cite{Gr}.
\end{remark}


\chapter{Magnetic fields}


Let us consider the  nonlinear Schr\"{o}dinger equation
\begin{equation}\label{CS:eq:1.1}
\I h \frac{\de \psi}{\de t} = \left( \frac{h}{\I} \nabla -A(x)\right)^2
\psi +U(x)\psi - f(x,\psi), \quad x\in\Rn
\end{equation}
where $t \in \R$, $x \in \Rn$ $(n \geq 2)$.
The function $\psi(x,t)$ takes on complex values,
$h$ is the Planck constant, $\I$ is the imaginary unit.
Here  $A\colon \Rn \to \Rn$ denotes a magnetic potential and
the Schr\"odinger operator is defined by
$$
\left( \frac{h}{\I} \nabla -A(x)\right)^2 \psi
:= -h^2 \varDelta \psi - \frac{2h}{\I} A \cdot \nabla \psi + |A|^2 \psi -
\frac{h}{\I} \psi\, \operatorname{div}A\,  .
$$
Actually, in general dimension $n \geq 2$,
the magnetic field $B$ is a $2$-form where
$B_{i,j}= \partial_j A_k - \partial_k A_j$ ;
in the case $n=3$, $B= \operatorname{curl}A$.
The function $U\colon \Rn \to \R$ represents an electric potential.
In the sequel, for the sake of simplicity, we limit ourselves
to the particular case in which $f(x,t)=K(x)|t|^{p-1}t$, with $p>1$ if
$n=2$ and $1<p<\frac{n+2}{n-2}$ if $n\geq 3$.

It is now well known that the
nonlinear Schr\"{o}dinger equation \eqref{CS:eq:1.1}
arises from a perturbation approximation for strongly nonlinear
dispersive wave systems. Many papers are devoted to the nonlinear
Schr\"{o}dinger equation and its solitary wave solutions.

In this section we seek for standing wave solutions to
\eqref{CS:eq:1.1}, namely waves of the form
$\psi (x,t)=e^{-iEth^{-1}}u(x)$ for some function
$u\colon \Rn \longrightarrow \C$.
Substituting this
{\it ansatz} into \eqref{CS:eq:1.1}, and denoting for
convenience $\ge=h$, one is led to solve the complex equation in $\Rn$
\begin{equation}
\left( \tfrac{\e}{\I} \nabla -A(x)\right)^2 u +
(U(x)-E) u = K(x)|u|^{p-1}u \,.
\tag{NLS}
\end{equation}
Renaming $V(x)+1=U(x)-E$, we assume from now on that $1+V$ is
strictly positive on the whole $\Rn$.
Moreover, by an obvious change of variables, 
the problem becomes that of finding some function
$u\colon\Rn\to\C$ such that
\begin{equation}\label{S}
\left( \tfrac{\nabla}{\I} -A(\e x)\right)^2 u + u +
V(\e x)u = K(\e x)|u|^{p-1}u, \quad x\in\Rn \,.
\tag{$S_\e$}
\end{equation}
Concerning nonlinear Schr\"odinger equation with external magnetic field,
we firstly quote a paper by Esteban and Lions \cite{EL}, where
concentrations and compactness arguments are applied to solve some
minimization problems associated to $(S_\e)$ under
suitable assumptions on the magnetic field.

The purpose of this section is 
to study the time--independent nonlinear
Schr\"{o}dinger equation $(S_\e)$ in the semiclassical limit.
This seems a very interesting problem since
the Correspondence's Principle establishes that
Classical Mechanics is, roughly speaking, contained in
Quantum Mechanics.
The mathematical transition is obtained letting to zero
the Planck constant ($\e \to 0$)
and solutions $u(x)$ of $(S_\e)$
which exist for small value of $\e$ are
usually referred as semi-classical ones (see \cite{RS}).

We remark that in the
linear case, Helffer {\it et al.} in
\cite{hel1,helS} have studied the asymptotic
behavior of the eigenfunctions
of the  Schr\"{o}dinger operators with magnetic fields
in the semiclassical limit.
Note that in these papers
the {\sl wells} of the  Schr\"{o}dinger operators with magnetic fields
are the same as those without magnetic field,
so that one doesn't `see'  the magnetic field in the
definition of the {\sl well}. See also \cite{brum}
for generalization of the results by \cite{helS}
for potentials which degenerate at infinity.

In the case $A=0$, (no magnetic field),
a recent extensive literature is devoted to
study the time--independent nonlinear
Schr\"{o}dinger equation $(S_\e)$ in the semi-classical limit.
We shortly recall the main results in literature.
The first paper is due to Floer and Weinstein which
investigated the one-dimensional nonlinear Schr\"{o}dinger equation
(with $K(x)=1$)  and gave a description of the limit behavior of $u(x)$
as $\e \to 0$.
Really they proved that if the potential $V$ has a non-degenerate
critical point, then $u(x)$ concentrates near this
critical point, as $\e \to 0$.

Later, other authors proved that this problem is really local in nature and
the presence of an isolated critical point of the potential $V$
(in the case $K(x)=1$)
produces a semi-classical solution $u(x)$ of $(S_\e)$ which
concentrates near this point.
Different approaches are used to cover different cases
(see \cite{ABC,dPF,Li,Oh,WZ}).
Moreover when $V$ oscillates, the existence of multibumps solutions has also
been studied in \cite{ABer,CN,dPFM,G}.
Furthermore  multiplicity results are obtained in 
 \cite{CL1,CL2} for potentials $V$ having a set of degenerate
global minima and recently in \cite{ams},
for potentials $V$ having a set of
critical points, not necessarily global minima.

\smallskip
A natural answer arises:
{\sl how does the presence of an external magnetic field 
influence the existence and the
concentration behavior of
standing wave solutions to \eqref{CS:eq:1.1} in the semi-classical limit}?
\smallskip

A first result in this direction is contained in \cite{ku}
where Kurata has proved the existence of
least energy solutions to $(S_\e)$ for any $\e >0$,
under some assumptions linking the magnetic field
$B= (B_{i,j})$ and the electric potential $V(x)$.
The author also investigated the semi-classical limit of the found
least energy solutions and
showed a concentration phenomenon near global
minima of the electric potential in the case $K(x) = 1$
and  $|A|$ is small enough.

Recently in \cite{cingolani}, Cingolani
obtained a multiplicity result of semi-classical standing waves
solutions to $(S_\e)$,
relating the number of solutions to $(S_\e)$ to
the {\sl richness} of a set $M$ of global minima 
of an auxiliary function $\Lambda$ defined by setting
\[
\Lambda (x)=\frac{(1+V(x))^\theta}{K(x)^{-2/(p-1)}},\qquad \theta =
\frac{p+1}{p-1}-\frac{n}{2},
\]
(see $(\ref{def:auxiliary})$ in section 4 for details) depending on $V(x)$ and $K(x)$.
We remark that, if $K(x)=1$ for any $x \in \Rn$, global
minima of $\Lambda$ coincides with global minima of $V$.
The variational approach, used in \cite{cingolani},
allows to deal with unbounded potential $V$
and does not require assumptions on the magnetic field.
However this approach works only near global minima of
$\Lambda$.

In the present section we deal with the more
general case in which the auxiliary function
$\Lambda$ has a manifold $M$ of stationary points,
not necessarily global minima. 
For bounded magnetic potentials $A$, 
we are able to prove a multiplicity result of semi-classical
standing waves of $(S_\e)$,
following the perturbation approach
described in the previous chapter.

Now we briefly describe the proof of the result.
First of all, we highlight that solutions of $(S_\e)$ naturally
appear as {\em orbits}: in fact, equation $(S_\e)$ is invariant under the
multiplicative action of $S^1$. Since there is no danger of confusion, we
simply speak about solutions.
The complex--valued solutions to $(S_\e)$ are
found near least energy solutions of the equation
\begin{equation}\label{risc}
\bigg(\frac{\nabla}{\I} - A(\e \xi)\bigg)^2u + u + V(\e \xi)u=K(\e\xi)|u|^{p-1}u.
\end{equation}
where $\e \xi$ is in a neighborhood of $M$.
The least energy of (\ref{risc}) have the form
\begin{equation}\label{def:zed}
z^{\e\xi ,\sigma} \colon x\in\Rn \mapsto e^{\I\sigma +
\I A(\e\xi)\cdot x}
\left(\frac{1+V(\e\xi)}{K(\e\xi)}\right)^\frac{1}{p-1} 
U((1+V(\e\xi))^{1/2} (x-\xi))
\end{equation}
where $\e \xi$ belongs to $M$ and $\sigma \in [0, 2 \pi]$.
As before, the proof relies on a suitable
finite dimensional reduction,
and critical points of the Euler functional $f_\e$ associated to problem $(S_\e)$
are found near critical point of a finite dimensional
functional $\Phi_\e$ which is defined on a suitable neighborhood of $M$.
This allows to use Lusternik-Schnirelman category in the case $M$ is a
set of local maxima or minima of $\Lambda$.
We remark that the case of maxima cannot be handled by using
direct variational arguments as in \cite{cingolani}.

Again, under suitable assumptions on $M$, more than one solution can be found.

Firstly we present a special case of our results.

\begin{theorem}\label{th:nond}
Assume that
\begin{description}
\item[$(K1)$] \, $K\in L^\infty (\Rn)\cap C^2(\Rn)$ is strictly positive and $K''$ is
bounded; 
\item[$(V1)$] \, $V\in L^\infty (\Rn)\cap C^2(\Rn)$ satisfies $\inf_{x\in\Rn}
(1+V(x)) > 0$, and $V''$ is bounded;
\item[$(A1)$] \, $A\in L^\infty (\Rn,\Rn)\cap C^1(\Rn,\Rn)$, and
the jacobian $J_A$ of $A$ is globally bounded in $\Rn$.
\end{description}

\noindent
If the auxiliary function $\Lambda$ has
a non-degenerate critical point $x_0\in \Rn$, then
for $\e > 0$ small, the problem $(S_\e)$ has at least a
(orbit of) solution concentrating near $x_0$.
\end{theorem}

Actually, we are able to prove the following generalization.

\begin{theorem}\label{CS:th:main}
As in Theorem 1.1, assume again (K1), (V1) and (A1) hold.

\noindent
If the auxiliary function $\Lambda$ has
a smooth, compact, non-degenerate manifold of critical points $M$,
then for $\e > 0$ small, the problem $(S_\e)$ has at least
$\ell (M)$ (orbits of) solutions concentrating near points of $M$.
\end{theorem}

Finally we point out that the presence of
an external magnetic field
produces a phase in the complex wave which depends on
the value of $A$ near $M$.
Conversely the presence of $A$ does not seem to influence the location of
the peaks of the modulus of the complex wave. Although we will not deal with this problem, 
we believe that in order to have a 
local $C^2$ convergence of the solutions, some assumption about the smallness of the 
magnetic potential $A$ should be added, as done in \cite{ku} for minima of $V$.

\bigskip
Finally we point out that
Theorem~\ref{th:nond} and Theorem~\ref{CS:th:main} hold for problems involving
more general nonlinearities. See Remark \ref{nonl} in the last section.

\vspace{10pt}

\section{The variational framework}
We work in the real Hilbert space $E$ obtained as the completion of
$C_0^\infty (\Rn,\C)$ with respect to the norm associated to the inner product
\[
\langle u\mid v\rangle = \operatorname{Re} \int_\Rn \nabla u \cdot \overline{\nabla v} + u \bar v 
..
\]
Solutions to ($S_\e$) are, under some conditions we are going to point out, critical points of the 
functional formally defined on $E$ as
\begin{multline}\label{feps}
f_\e (u)=\frac{1}{2}\int_\Rn \bigg(\bigg| \bigg(\frac{1}{\I}\nabla - A(\e x)\bigg)u\bigg|^2 
+|u|^2+V(\e x)|u|^2\bigg)\, dx \\- \frac{1}{p+1}\int_\Rn K(\e x)|u|^{p+1}\, dx.
\end{multline}
In what follows, we shall assume that the functions $V$, $K$ and $A$
satisfy the following assumptions:
\begin{description}
\item[(K1)] $K\in L^\infty (\Rn)\cap C^2(\Rn)$ is strictly positive and $K''$ is
bounded; 
\item[(V1)] $V\in L^\infty (\Rn)\cap C^2(\Rn)$ satisfies $\inf_{x\in\Rn}
(1+V(x)) > 0$, and $V''$ is bounded;
\item[(A1)] $A\in L^\infty (\Rn,\Rn)\cap C^1(\Rn,\Rn)$, and
the jacobian $J_A$ of $A$ is globally bounded in $\Rn$.
\end{description}
Indeed,
\begin{multline*}
\int_\Rn \bigg(\bigg| \bigg(\frac{1}{\I}\nabla - A(\e x)\bigg)u\bigg|^2 \bigg)\, dx
= \\ = \int_\Rn \bigg(|\nabla u|^2 + |A(\e x)u|^2 -2
\operatorname{Re} ( \textstyle\frac{\nabla u}{\I} \cdot A(\e x)\overline{u}) \bigg)\, dx,
\end{multline*}
and the last integral is finite thanks to the Cauchy--Schwartz
inequality and the boundedness of $A$.

\noindent It follows that $f_\e$ is actually well-defined on $E$.

In order to find possibly multiple critical points of \eqref{feps},
we follow the approach already presented in the previous chapter. In our context, we need to find
complex--valued solutions, and so some further remarks are due.

Let $\xi\in\Rn$, which will be fixed suitably later on: we look for solutions to \eqref{S} ``close'' to 
a particular solution of the equation
\begin{equation}\label{eqlimite}
\bigg(\frac{\nabla}{\I} - A(\e \xi)\bigg)^2u + u + V(\e\xi)u=K(\e\xi)|u|^{p-1}u.
\end{equation}
More precisely, we denote by $U_c \colon \Rn \to \C$ a least--energy solution
to the scalar problem
\begin{equation}\label{CS:eq:4}
-\varDelta U_c + U_c + V(\e\xi)U_c = K(\e\xi) |U_c|^{p-1}U_c
\mbox{\quad in }\Rn.
\end{equation}
By energy comparison (see \cite{ku}), one has that
\[
U_c (x)=e^{\I\sigma} U^\xi (x-y_0)
\]
for some choice of $\sigma \in [0,2\pi]$ and $y_0\in\Rn$, where
$U^\xi\colon\Rn\to\R$ is the unique solution of the problem
\begin{equation}
\begin{cases}
-\varDelta U^\xi + U^\xi + V(\e\xi) U^\xi = K(\e\xi)| U^\xi|^{p-1} U^\xi  \\
U^\xi (0) = \max_\Rn U^\xi\\
U^\xi > 0.
\end{cases}
\end{equation}
If $U$ denotes the unique solution of
\begin{equation}
\begin{cases}
-\varDelta U + U = U^p \mbox{\quad in }\Rn\\
U(0) = \max_\Rn U \\
U > 0,
\end{cases}
\end{equation}
then some elementary computations prove that
$U^\xi (x)=\alpha(\e\xi)U(\beta(\e\xi)x)$,
where
\begin{eqnarray*}
\alpha (\e\xi) &= \left(\frac{1+V(\e\xi)}{K(\e\xi)}\right)^\frac{1}{p-1} \\
\beta(\e\xi)&=(1+V(\e\xi))^{1/2}.
\end{eqnarray*}
It is easy to show, by direct computation, that the function
$u(x)= e^{i A(\e\xi)\cdot x} U_c (x)$ actually solves \eqref{eqlimite}.

For $\xi \in \Rn$ and $\sigma\in [0,2\pi]$, we set
\begin{equation}\label{def:zeta}
z^{\e\xi ,\sigma} \colon x\in\Rn \mapsto e^{\I\sigma + \I A(\e\xi)\cdot 
x}\alpha(\e\xi)U(\beta(\e\xi)(x-\xi)).
\end{equation}

Sometimes, for convenience, we shall identify $[0,2\pi]$ and $S^1\subset \C$, through $\eta = 
e^{i\sigma}$.

Introduce now the functional $F^{\e\xi,\sigma}\colon E \to \R$ defined by
\begin{multline*}
F^{\e\xi,\sigma}(u)=\frac{1}{2}\int_\Rn  \bigg(\bigg|\bigg(\frac{\nabla u}{\I} -
A(\e\xi)u \bigg)\bigg|^2 + |u|^2 + V(\e\xi)|u|^2 \bigg) \, dx \\ - \frac{1}{p+1}
\int_\Rn K(\e\xi)|u|^{p+1}\, dx,
\end{multline*}
whose critical point correspond to solutions of \eqref{eqlimite}.

The set
\[
Z^\e =\{{z^{\e\xi,\sigma}}\mid \xi \in \Rn \, \land \,
\sigma \in [0,2\pi]\} \simeq S^1 \times \Rn
\]
is a regular manifolds of critical points for the functional $F^{\e\xi,\sigma}$.

It follows from elementary differential geometry that
\begin{multline*}
T_{z^{\e\xi,\eta}}Z^\e = \operatorname{span}_{\R}
\{\tfrac{\de}{\de \sigma}
{{{{{z^{\e\xi,\sigma}}}}}},\tfrac{\de}{\de \xi_1}
{{z^{\e\xi,\sigma}}},\dots,\tfrac{\de}{\de \xi_n}
{{{z^{\e\xi,\sigma}}}}\}=\\ =\operatorname{span}_{\R}
\{\I{{{z^{\e\xi,\sigma}}}},\tfrac{\de}{\de \xi_1}
{{{z^{\e\xi,\sigma}}}},\dots,\tfrac{\de}{\de \xi_n}{{z^{\e\xi,\sigma}}}\},
\end{multline*}
where we mean by the symbol $\operatorname{span}_{\R}$ that all
the linear combinations must have real coefficients.

We remark that, for $j=1,\dots,n$,
\begin{multline*}
\frac{\de}{\de \xi_j}{{z^{\e\xi,\sigma}}}=-\frac{\de}{\de x_j}{{z^{\e\xi,\sigma}}} + O(\e|\nabla 
V(\e\xi)|)+\\
+\I\alpha (\e\xi)e^{i A(\e\xi)\cdot x + \I\sigma} U(\beta(\e\xi)(x-\xi))\left(\frac{\de}{\de 
\xi_j}\left(A(\e\xi)\cdot x\right)
+A_j(\e\xi)\right)=\\
=-\frac{\de}{\de x_j} {{z^{\e\xi,\sigma}}} + O(\e|\nabla V(\e\xi)|+\e|J_A 
(\e\xi)|)+\I{{{z^{\e\xi,\sigma}}}}
A_j(\e\xi),
\end{multline*}
so that
\[
\frac{\de}{\de \xi_j}{{z^{\e\xi,\sigma}}} =
- \frac{\de}{\de x_j}{{z^{\e\xi,\sigma}}} +
\I {z^{\e\xi,\sigma}}A_j(\e\xi)+ O(\e).
\]
Collecting these remarks,
we get that any
$\zeta \in T_{{{{{z^{\e\xi,\sigma}}}}}}Z^\e$ can be written as
$$
\zeta = \I\ell_1 {{{z^{\e\xi,\sigma}}}}+ \sum_{j=2}^{n+1} \ell_j
\frac{\de}{\de x_{j-1}}{{{{z^{\e\xi,\sigma}}}}}+O(\e)$$
for some real coefficients $\ell_1,\ell_2,\dots,\ell_{n+1}$.

The next lemma shows that $\nabla f_\e ({{{z^{\e\xi,\sigma}}}})$ gets small when
$\e\to 0$.

\begin{lemma}\label{CS:lem:1}
For all $\xi\in \Rn$, all $\eta \in S^1$ and all $\e>0$ small, one has that 
$$ 
\|\nabla f_\e ({z^{\e\xi,\sigma}})\|\leq C\left(\e |\nabla 
V(\e\xi)|+\e |\nabla K(\e\xi)|+\e |J_A (\e\xi)|+\e 
|\operatorname{div} A(\e\xi)|+\e^{2}\right),  $$ 
for some constant $C>0$.
 \end{lemma} 
\begin{proof} 
From 
\begin{multline}\label{CS:eq:7}
f_{\e}(u)  =  F^{\e\xi,\eta}(u)+\frac{1}{2}\int_\Rn \bigg( \bigg| \frac{\nabla
u}{\I} - A(\e x)u \bigg|^2 - \bigg| \frac{\nabla u}{\I}-A(\e\xi)u
\bigg|^2 \biggr)  +\\
+ \frac{1}{2}\int_{\Rn}\left[V(\e x)-V(\e\xi)\right]u^{2}
-\frac{1}{p+1}\int_\Rn [K(\e x)-K(\e\xi)]|u|^{p+1}
\end{multline}
 and since $z^{\e\xi,\eta}$ is a
critical point of $F^{\e\xi,\eta}$, one has (with $z=z^{\e\xi,\eta}$)
\begin{multline*} 
\langle \nabla f_\e (z)\mid v \rangle = \re \int_\Rn \bigg(
\frac{1}{\I}\nabla - A(\e\xi)\bigg)z \cdot (A(\e\xi)-A(\e x))\bar v 
\\+ \re\int_\Rn (A(\e\xi)-A(\e x))z \cdot \overline{\bigg( \frac{1}{\I}\nabla -
A(\e\xi)\bigg)v}  + \\ \re\int_\Rn (A(\e\xi)-A(\e x))z \cdot
(A(\e\xi)-A(\e x))\bar v \\
+ \re\int_\Rn (V(\e x)-V(\e\xi))z\bar v  - \re\int_\Rn (K(\e
x)-K(\e\xi))|z|^{p-2}z\bar v \\ = \e \operatorname{Re} \int_\Rn \frac{1}{\I}
(\operatorname{div}\, A(\e x)) z \bar v  + 2 \operatorname{Re} \int_\Rn
(A(\e\xi)-A(\e x))z \cdot \overline{\bigg( \frac{1}{\I} \nabla -
A(\e\xi)\bigg)v} \\+\re\int_\Rn (V(\e x)-V(\e\xi))z\bar v  -
\re\int_\Rn (K(\e x)-K(\e\xi))|z|^{p-2}z\bar v 
\end{multline*}

\noindent
From the assumption that $|D^{2}V(x)|\leq {\rm const.}$ one infers 
\[ 
|V(\e x)-V(\e\xi)|\leq \e |\nabla V(\e\xi)|\cdot
|x-\xi|+c_{1}\e^{2}|x-\xi|^{2}.
\]
This implies 
\begin{eqnarray} 
 \int_{\Rn}|V(\e x)-V(\e\xi)|^{2}{{{z^{\e\xi,\sigma}}}}^{2}
 &\leq&  c_1\e^{2}|\nabla V(\e\xi)|^{2}\int_{\Rn}|x-\xi|^{2}z^{2}(x-\xi) + \nonumber \\ 
 & & c_{2}\e^{4}\int_{\Rn}|x-\xi|^{4}z^{2}(x-\xi). 
\label{CS:eq:1.3} 
\end{eqnarray} 
A direct calculation yields 
\begin{eqnarray*} 
\int_{\Rn}|x-\xi|^{2}z^{2}(x-\xi) & = & 
\a^{2}(\e\xi)\int_{\Rn}|y|^{2}U^{2}(\b(\e\xi) y)dy \\ 
& = & \alpha (\e\xi)^{2} \beta (\e\xi)^{-n-2}\int_{\Rn}|y'|^{2}U^{2}(y')dy'\leq c_{3}. 
\end{eqnarray*} 
{From} this (and a similar calculation for the last integral in the 
above formula) one infers 
\begin{equation} 
\int_{\Rn}|V(\e x)-V(\e\xi)|^{2}|z^{\e\xi,\sigma}|^{2}
\leq c_{4}\e^{2}|\nabla
V(\e\xi)|^{2} + c_{5}\e^{4}. 
\label{CS:eq:1.4} 
\end{equation} 
Of course, similar estimates hold for the terms involving $K$.
It then follows that
\[
\|\nabla f_\e (z^{\e\xi,\eta})\| \leq C(\e | \operatorname{div}A(\e\xi)| + \e
| \nabla V (\e\xi)| + \e | J_A (\e\xi)| + \e^2),
\]
and the lemma is proved.
  \end{proof} 

\section{The invertibility of $D^2f_\e$ on $(TZ^\e)^\bot$}

To apply the perturbative method, we need to exploit some non--degeneracy properties of
the solution ${z^{\e\xi,\sigma}}$ as a critical point of $F^{\e\xi,\sigma}$.

Let $L_{\e,\sigma,\xi}\colon (T_{{{z^{\e\xi,\sigma}}}}Z^\e)^\bot
\longrightarrow
(T_{{z^{\e\xi,\sigma}}}Z)^\bot$ be the operator defined by \[\langle
L_{\e,\sigma,\xi} v \mid w \rangle =
D^2 f_\e ({{{z^{\e\xi,\sigma}}}})(v,w)\] for all
$v,w \in (T_{{z^{\e\xi,\sigma}}}Z^\e)^\bot$.

The following elementary result will play a fundamental r\^{o}le in the present
section.

\begin{lemma}
Let $M \subset \Rn$ be a bounded set. 
Then there exists a constant $C>0$ such that for all $\xi\in M$ one has
\begin{equation}\label{norme:equiv}
\int_\Rn \bigg| \bigg( \frac{\nabla}{\I} - A(\xi) \bigg) u \bigg|^2 + |u|^2 \geq C \int_\Rn (|\nabla 
u|^2 + |u|^2) \qquad\forall u \in E.
\end{equation}
\end{lemma}

\begin{proof}
To get a contradiction, we assume on the contrary the existence of a sequence
$\{ \xi_n\}$ in $M$ and a sequence $\{u_n\}$ in $E$ such that $\|u_n\|_E = 1$
for all $n\in\mathbb{N}$ and 
\begin{equation}\label{CS:eq:contradiction}
\lim_{n\to + \infty} \bigg[ \int_\Rn \bigg| \bigg( \frac{\n}{\I} - A(\xi)
\bigg)u_n \bigg|^2 + \int_\Rn |u_n|^2 \bigg] = 0.
\end{equation}
In particular, $u_n \to 0$ strongly in $L^2 (\Rn,\C)$. Moreover, since $M$
is bounded, we can assume also $\xi_n \to \xi^* \in \overline{M}$ as
$n\to\infty$. From 

\begin{multline*}
\int_\Rn \bigg| \bigg( \frac{\n}{\I} - A(\xi_n)
\bigg)u_n \bigg|^2 = \\
 \int_\Rn \left( |\n u_n|^2 + |A(\xi_n)|^2 |u_n|^2 - 2 \re
\frac{1}{\I} \n u_n \cdot A(\xi_n) \overline{u_n} \right)
\end{multline*}
we get 
\begin{eqnarray*}
\lim_{n \to + \infty} \int_\Rn | \n u_n |^2 &=& 1 \\
\lim_{n \to + \infty} \re \int_\Rn \tfrac{1}{\I} \n u_n \cdot A(\xi_n) \overline{u_n} &=& \frac{1}{2} \,.
\end{eqnarray*}
Therefore,
\begin{eqnarray*}
\limsup_{n\to\infty} \int_\Rn |\n u_n| \, |A(\xi_n)| \, |u_n| &\geq&
\limsup_{n\to\infty} \left| \int_\Rn \textstyle\frac{1}{\I} \n u_n \cdot
A(\xi_n) \overline{u_n} \right| \\
&\geq& \limsup_{n\to\infty} \re \int_\Rn \tfrac{1}{\I} \n u_n \cdot A(\xi_n)
\overline{u_n} = \frac{1}{2} \,.
\end{eqnarray*}
From this we conclude that
\begin{eqnarray*}
\frac{1}{2} &\leq& |A(\xi^*)| \limsup_{n\to\infty} \| \n u_n \|_{L^2} \|u_n
\|_{L^2} \leq \\
&\leq& |A(\xi^*)| \limsup_{n\to\infty} \|u_n \|_{L^2} = 0,
\end{eqnarray*}
which is clearly absurd. This completes the proof of the lemma.
\end{proof}

At this point we shall prove the following result:

\begin{lemma}\label{lemma:3.1}
Given $\bar\xi >0$, there exists $C > 0$ such that for $\e$ small enough one
has
\begin{equation}\label{CS:eq:posdef}
|\langle L_{\e,\sigma,\xi}v\mid v \rangle | \geq C \|v\|^2, \quad \forall |\xi|
\leq \bar\xi, \;\forall\sigma\in [0,2\pi],\;\forall v\in
(T_{{z^{\e\xi,\sigma}}}Z^\e)^\bot.
\end{equation}
\end{lemma}

\begin{proof}
We follow the arguments in \cite{ams}, with some minor modifications due to the presence of 
$A$. Recall that $$T_{{z^{\e\xi,\sigma}}}Z=\operatorname{span}_{\R} \{ \frac{\de}{\de
\xi_1}{z^{\e\xi,\sigma}},\dots ,\frac{\de}{\de
\xi_n}{z^{\e\xi,\sigma}},{i z^{\e\xi,\sigma}} \}.$$
Define 
$$
\mathcal{V}=\operatorname{span}_{\R} \{ \tfrac{\de}{\de
x_1}{z^{\e\xi,\sigma}},\dots ,\tfrac{\de}{\de
x_n}{z^{\e\xi,\sigma}},{z^{\e\xi,\sigma}} ,i{z^{\e\xi,\sigma}}\}.
$$
As in the previous chapter, it suffices to prove \eqref{CS:eq:posdef} for all
$v\in\operatorname{span}_\R \{{{z^{\e\xi,\sigma}}},\phi\}$,
where $\phi \perp \mathcal{V}$.
More precisely, we shall prove that for some constants $C_1 >0$, $C_2 > 0$,
for all $\e$ small enough and all $|\xi|\leq \bar\xi$ the following hold:
\begin{equation}\label{CS:eq:13}
\langle L_{\e,\sigma,\xi}{{z^{\e\xi,\sigma}}}\mid {z^{\e\xi,\sigma}} \rangle \leq
-C_1 < 0,
\end{equation}
\begin{equation}\label{CS:eq:14}
\langle L_{\e,\sigma,\xi}\phi\mid \phi \rangle
\geq C_2 \|\phi\|^2 \quad\forall\phi\perp\mathcal{V}.
\end{equation}
For the reader's convenience, we reproduce here the expression for the second
derivative of $F^{\e\xi,\sigma}$:
\begin{multline*}
D^2 F^{\e\xi,\sigma}(u)(v,v)=\int_\Rn \bigg| \bigg(
\frac{\nabla}{\I}-A(\e\xi)\bigg)v\bigg|^2 + |v|^2 + V(\e\xi)|v|^2 \\
- K(\e\xi)\bigg[(p-1) \re \int_\Rn |u|^{p-3} \re (u \bar v) u \bar v  +
\int_\Rn |u|^{p-1}|v|^2 \bigg].
\end{multline*}
Moreover, since ${{z^{\e\xi,\sigma}}}$ is a solution of \eqref{eqlimite}, we
immediately get
\[
\int_\Rn
\bigg(  \bigg| \bigg( \frac{\nabla}{\I}-A(\e\xi)\bigg){{{{z^{\e\xi,\sigma}}}}}\bigg|^2
+ V(\e\xi)|{z^{\e\xi,\sigma}}|^2 +|{{z^{\e\xi,\sigma}}}|^2
 \bigg)  = K(\e\xi)\int_\Rn
|{z^{\e\xi,\sigma}}|^{p+1}.
\]
From this it follows readily that we can find some $c_0 > 0$ such that for all
$\e > 0$ small, all $|\xi| \leq \bar \xi$ and all $\sigma\in [0,2\pi]$ it
results
\begin{equation} D^2
F^{\e\xi,\sigma}({{z^{\e\xi,\sigma}}}) ({z^{\e\xi,\sigma}},{z^{\e\xi,\sigma}})<c_0 <0.
\end{equation}
Recalling \eqref{CS:eq:7}, we find
\begin{multline*}
\langle L_{\e,\sigma,\xi}{{z^{\e\xi,\sigma}}}\mid {z^{\e\xi,\sigma}} \rangle =D^2
F^{\e\xi,\sigma}({{z^{\e\xi,\sigma}}})({z^{\e\xi,\sigma}},{z^{\e\xi,\sigma}})+ \\ +\int_\Rn [V(\e 
x)-V(\e\xi)]|{{z^{\e\xi,\sigma}}}|^2
-\int_\Rn [K(\e x)-K(\e\xi)]|{{z^{\e\xi,\sigma}}}|^{p+1}
+\\ \int_\Rn
\bigg( \bigg| \bigg(\frac{\nabla}{\I}-A(\e x)\bigg){{{z^{\e\xi,\sigma}}}}\bigg|^2 - \bigg|
\bigg(\frac{\nabla}{\I}-A(\e \xi)\bigg){{{z^{\e\xi,\sigma}}}}\bigg|^2
\bigg) \,.
\end{multline*}
It follows that
\begin{multline}
\langle L_{\e,\sigma,\xi}{{z^{\e\xi,\sigma}}}\mid {z^{\e\xi,\sigma}} \rangle \leq D^2
F^{\e\xi,\sigma}({{z^{\e\xi,\sigma}}})({z^{\e\xi,\sigma}},{z^{\e\xi,\sigma}}) + \\ +c_1
\e |\nabla V(\e\xi)|+c_2\e |\nabla K (\e\xi)| + c_3 \e |J_A (\e\xi)| + c_4
\e^2. \end{multline}
Hence \eqref{CS:eq:13} follows. The proof
of \eqref{CS:eq:14} is more involved. We first prove the following
claim.

\vspace{5pt}

{\em Claim.} There results \begin{equation}\label{CS:eq:17}
D^2 F^{\e\xi}({z^{\e\xi,\sigma}})(\phi,\phi) \geq c_1 \|\phi\|^2 \quad
\forall\phi\perp\mathcal{V}.
\end{equation}

\vspace{3pt}

Recall that the complex ground state $U_c$ introduced in \eqref{CS:eq:4} is a
critical point of mountain-pass type for the corresponding energy functional
$J\colon E \longrightarrow \R$ defined by
\begin{equation}
J(u)=\frac{1}{2}\int_\Rn ( |\nabla u|^2 + |u|^2 + V(\e\xi) |u|^2 ) -
\frac{1}{p+1}\int_\Rn K(\e\xi)|u|^{p+1}  .
\end{equation}
Let
\[
\mathcal{M} = \left\{u\in E \colon
\textstyle\int_\Rn (|\n u|^2 + |u|^2 + V(\e\xi)|u|^2 ) 
=\textstyle \int_\Rn |u|^{p+1}  \right\}
\]
be the Nehary manifold of $J$, which has codimension one. Let
\[
\mathcal{N}=\left\{u\in E \colon \textstyle\int_\Rn \left( \left| \left( \frac{\n}{\I} - A(\e\xi)
\right) u\right|^2 + |u|^2 + V(\e\xi)|u|^2 \right) 
= \textstyle\int_\Rn |u|^{p+1} \ \right\}
\]
be the Nehari manifold of $F^{\e\xi,\sigma}$.
One checks readily that $\operatorname{codim} \mathcal{N}=1$. 
Recall (\cite{ku}) that $U_c$ is, up to multiplication by a constant phase,
the {\em unique} minimum of $J$ restricted to $\mathcal{M}$. Now, for every
$u\in\mathcal{M}$, the function $x \mapsto e^{i A(\e\xi)\cdot x} u(x)$ lies in
$\mathcal{N}$, and viceversa. Moreover
\[
J(u)=F^{\e\xi,\sigma} (e^{i A(\e\xi)\cdot x} u).
\]
This immediately implies that $\min_{\mathcal{N}} F^{\e\xi,\sigma}$ is
achieved at a point which differs from $e^{i A(\e\xi)\cdot x} U_c (x)$ at most
for a constant phase. In other words, $z^{\e\xi,\sigma}$ is a critical point
for $F^{\e\xi,\sigma}$ of mountain-pass type, and the claim follows by standard
results (see \cite{Chang}).

\medskip

Let $R\gg 1$ and consider a radial smooth function 
$\chi_{1}:\R^{n}\longrightarrow \R$ such that 
\begin{equation}\label{CS:eq:c1}  \chi_{1}(x) = 1, \quad \hbox{ for
} |x| \leq R; \qquad  \chi_{1}(x) = 0, \quad \hbox{ for } |x| \geq 2 R; 
\end{equation} 
\begin{equation}\label{CS:eq:c2} 
|\n \chi_{1}(x)| \leq \frac{2}{R}, \quad \hbox{ for } R \leq |x| \leq 2 R. 
\end{equation} 
We also set $ \chi_{2}(x)=1-\chi_{1}(x)$. 
Given $\phi$ let us consider the functions 
$$ 
\phi_{i}(x)=\chi_{i}(x-\xi)\phi(x),\quad i=1,2. 
$$ 
A straightforward computation yields: 
$$ 
\irn |\phi|^2 = \irn |\phi_1|^2 + \irn |\phi_2|^2 + 2\re\irn \phi_{1} \, \bar 
\phi_{2}, 
$$ 
$$ 
\irn |\n \phi|^2 = \irn |\n \phi_1|^2 + \irn |\n \phi_2|^2 + 2\re\irn 
\n\phi_{1} \cdot \overline{\n \phi_{2}}, 
$$ 
and hence 
$$ 
\| \phi \|^2 = \| \phi_1 \|^2 + \| \phi_2 \|^2+ 
2 \re\irn\left[ \phi_{1} \, \overline \phi_{2}
+ \n\phi_{1} \cdot \overline{\n \phi_{2}}\right]. 
$$ 
Letting $I$ denote the last integral, 
one immediately finds: 
$$ 
I=\underbrace{\irn \chi_{1}\chi_{2}(\phi^{2}+|\n \phi|^{2}) }_{I_{\phi}} + 
\underbrace{\irn\phi^{2}\n\chi_{1}\cdot \n\chi_{2}  }_{I'}+ 
\underbrace{\irn (\phi_{2}\n\chi_{1}\cdot\overline{\n\phi}+\overline{\phi_{1}}\,\n 
\phi\cdot\n\chi_{2})}_{I''}. 
$$ 
Due to the definition of $\chi$, the two integrals $I'$ and $I''$ 
reduce to integrals from $R$ and $2R$, and thus they are $o_{R}(1)\|\phi\|^{2}$,
where $o_{R}(1)$ is a function which tends to $0$, as $R \to +\infty$.
As a consequence we have that 
\begin{equation}\label{CS:eq:d} 
\| \phi \|^2 = \| \phi_1 \|^2 + \| \phi_2 \|^2 + 2I_\phi +
o_R(1)\| \phi \|^2 \,.
\end{equation} 
After these preliminaries, let us evaluate the three terms in the 
equation below: 
$$ 
(L_{\e,\sigma,\xi}\phi|\phi)= 
\underbrace{(L_{\e,\sigma,\xi}\phi_{1}|\phi_{1})}_{\a_{1}}+ 
\underbrace{(L_{\e,\sigma,\xi}\phi_{2}|\phi_{2})}_{\a_{2}}+ 
2\underbrace{(L_{\e,\sigma,\xi}\phi_{1}|\phi_{2})}_{\a_{3}}. 
$$ 
One has: 
\begin{multline*}
\alpha_{1}=\langle L_{\e,\sigma,\xi}\phi_{1}\mid \phi_{1}
\rangle =D^2
F^{\e\xi,\sigma}({{z^{\e\xi,\sigma}}})(\phi_{1},\phi_{1}) +
\int_\Rn [V(\e x)-V(\e\xi)]|\phi_{1}|^2  \\ -\int_\Rn [K(\e
x)-K(\e\xi)]|\phi_{1}|^{p+1} + \int_\Rn \bigg|
\bigg( \bigg(\frac{\nabla}{\I}-A(\e x)\bigg)\phi_{1}\bigg|^2 - \bigg|
\bigg(\frac{\nabla}{\I}-A(\e \xi)\bigg)\phi_{1}\bigg|^2 \bigg) \, .
\end{multline*}

In order to use \eqref{CS:eq:17}, we introduce the function 
$\phi^*_{1}=\phi_{1}-\psi$, where 
$\psi$ is the projection of $\phi_{1}$ onto $\cal V$: 
\begin{multline*}
\psi=(\phi_{1}|z^{\e\xi,\sigma})z^{\e\xi,\sigma} \| z^{\e\xi,\sigma} \|^{-2} + 
(\phi_{1}|iz^{\e\xi,\sigma})iz^{\e\xi,\sigma}\|z^{\e\xi,\sigma}\|^{-2}\\
+\sum (\phi_{1}|\partial_{x_{i}}{{z^{\e\xi,\sigma}}}) 
\partial_{x_{i}}{z^{\e\xi,\sigma}}\| \partial_{x_{i}}{{{z^{\e\xi,\sigma}}}} \|^{-2}. 
\end{multline*}
Then we have: 
\begin{equation} 
D^{2}F^{\e\xi}[\phi_{1},\phi_{1}]= 
D^{2}F^{\e\xi}[\phi_{1}^*,\phi_{1}^*]+ 
D^{2}F^{\e\xi}[\psi,\psi]+2\re D^{2}F^{\e\xi}[\phi^*_{1},\psi] \,.
\label{CS:eq:alfa2} 
\end{equation} 
Since ${z^{\e\xi,\sigma}}$ is orthogonal to $\partial_{x_{i}}{{{z^{\e\xi,\sigma}}}}$, 
$i=1,\ldots,n$, then one 
readily 
checks that $\phi^*_{1}\perp {\cal V}$ and hence \eqref{CS:eq:17} 
implies 
\begin{equation} 
D^{2}F^{\e\xi}[\phi^*_{1},\phi^*_{1}]\geq c_{1} 
\|\phi^*_{1}\|^{2}. 
\label{CS:eq:alfa3} 
\end{equation} 
On the other side, since $(\phi|{{z^{\e\xi,\sigma}}})=0$ it follows: 
\begin{eqnarray*} 
        (\phi_{1}|{{z^{\e\xi,\sigma}}}) & = & (\phi|{z^{\e\xi,\sigma}})-(\phi_{2}|{{z^{\e\xi,\sigma}}})= 
-(\phi_{2}|{z^{\e\xi,\sigma}}) \\ 
& = & -\re\irn\phi_{2}{{z^{\e\xi,\sigma}}}-\re\irn \n {{z^{\e\xi,\sigma}}}\cdot \n \phi_{2} \\ 
&=& -\re\irn\chi_{2}(y)z(y)\phi(y+\xi)dy-\re\irn \n z(y)\cdot \n 
\chi_{2}(y)\phi(y+\xi)dy . 
\end{eqnarray*} 
Since $\chi_{2}(x)=0$ for all $|x|<R$, and since $z(x)\to 0$ as 
$|x|=R\to\infty$, we infer 
$(\phi_{1}|{{z^{\e\xi,\sigma}}})=o_{R}(1)\|\phi\|$.  Similarly one shows that 
$(\phi_{1}|\partial_{x}{{{{z^{\e\xi,\sigma}}}}})=o_{R}(1)\|\phi\|$ and it follows that 
\begin{equation}\label{CS:eq:alfa4} 
         \|\psi\|=o_{R}(1)\|\phi\|. 
\end{equation} 
We are now in position to estimate the last two terms in equation (\ref{CS:eq:alfa2}). 
Actually, using Lemma 3.1 we get 

\begin{eqnarray}\label{CS:eq:alfa5}
& & D^{2}F^{\e\xi}[\psi,\psi] \geq  C\|\psi\|^{2}
+ V(\e\xi)\irn\psi^{2} \\ &- & \notag  K(\e\xi)\bigg[ \re (p-1)
\int_\Rn |z^{\e\xi,\sigma}|^{p-3}
\re (z^{\e\xi,\sigma} \bar \psi) z^{\e\xi,\sigma} \bar \psi  
\\ & & \notag 
+\int_\Rn |z^{\e\xi,\sigma}|^{p-1}|\psi|^2 \bigg] =
o_{R}(1)\|\phi\|^{2}. 
\end{eqnarray}
The same arguments readily imply 
\begin{equation}\label{CS:eq:alfa6} 
\re D^{2}F^{\e\xi}[\phi^*_{1},\psi]=o_{R}(1)\|\phi\|^{2}.
\end{equation} 
Putting together (\ref{CS:eq:alfa3}), (\ref{CS:eq:alfa5}) and  (\ref{CS:eq:alfa6}) 
we infer 
\begin{equation} \label{CS:eq:alfa7} 
D^{2}F^{\e\xi}[\phi_{1},\phi_{1}]\geq C\|\phi_{1}\|^{2}+o_{R}(1)\|\phi\|^{2}.
\end{equation} 
Using arguments already carried out before, one has 
\begin{eqnarray*} 
\irn|V(\e x)-V(\e\xi)|\phi_{1}^{2} & \leq 
&\e c_{2}\irn|x-\xi|\chi_1^{2}(x-\xi)\phi^{2}(x) \\ 
& \leq & \e c_{3}\irn |y|\chi_1^{2}(y)\phi^{2}(y+\xi)dy \\ 
& \leq & \e c_{4} R \|\phi\|^{2},
\end{eqnarray*} 
and similarly for the terms containing $K$.
This and (\ref{CS:eq:alfa7}) yield 
\begin{equation} 
        \a_{1}=(L_{\e,\sigma,\xi}\phi_{1}|\phi_{1})\geq c_{5}
        \|\phi_{1}\|^{2}-
        \e c_{4}R \|\phi\|^{2}+o_{R}(1)\|\phi\|^{2}. 
\label{CS:eq:alfa8} 
\end{equation} 
Let us now estimate $\a_{2}$.  One finds 
\begin{equation} 
\a_{2} = \langle L_{\e,\sigma,\xi} \phi_{2} \mid \phi_{2} \rangle \geq 
c_{6} \|\phi_{2}\|^{2}+o_{R}(1)\|\phi\|^{2}.  \label{CS:eq:alfa9} 
\end{equation} 
In a quite similar way one shows that 
\begin{equation} 
\a_{3} \geq c_{7}I_{\phi}+o_{R}(1)\|\phi\|^{2}. 
\label{CS:eq:alfa10} 
\end{equation} 
Finally, (\ref{CS:eq:alfa8}), (\ref{CS:eq:alfa9}), (\ref{CS:eq:alfa10}) and the fact 
that $I_{\phi}\geq 0$, yield 
\begin{eqnarray*} 
(L_{\e,\sigma,\xi}\phi|\phi) & = & \a_{1}+\a_{2}+2\a_{3}\\ 
& \geq & c_{8}\left[\|\phi_{1}\|^{2}+\|\phi_{2}\|^{2} 
+2 I_{\phi}\right]-c_{9} R \e\|\phi\|^{2}+ o_{R}(1)\|\phi\|^{2}. 
\end{eqnarray*} 
Recalling (\ref{CS:eq:d}) we infer that 
$$ 
(L_{\e,\sigma,\xi}\phi|\phi)\geq c_{10}\|\phi\|^{2}-c_{9} R \e\|\phi\|^{2}+ 
o_{R}(1)\|\phi\|^{2}. 
$$ 
Taking $R = \e^{-1/2}$, and choosing $\e$ small,
equation  \eqref{CS:eq:14} follows. This completes the proof. 
\end{proof}

\section{The finite dimensional reduction}
\label{CS:secfdr} 
 In this Section we will show that the existence of critical points of $f_{\e}$ 
 can be reduced to the search of critical points of an auxiliary finite 
 dimensional functional.  The proof will be carried out in two 
 subsections dealing, respectively, with  a Liapunov-Schmidt reduction, and 
 with the behaviour of the auxiliary finite dimensional functional. 

\smallskip

\section{A Liapunov-Schmidt type reduction}\label{CS:subsec:LS} The main 
result of this section is the following lemma. 
 
\begin{lemma}\label{CS:lem:w} 
For $\e>0$ small, $|\xi|\leq \overline{\xi}$ and $\sigma\in [0,2\pi]$, there
exists a unique  $w=w(\e,\sigma,\xi)\in 
(T_{z^{\e\xi,\sigma}} Z^\e)^{\perp}$ such that 
$\nabla f_\e (z^{\e\xi,\sigma} + w)\in T_{z^{\e\xi,\sigma}} Z^\e$. 
Such a $w(\e,\sigma,\xi)$ 
 is of class $C^{2}$, resp.  $C^{1,p-1}$, with respect to $\xi$, provided that 
 $p\geq 2$, resp.  $1<p<2$. 
Moreover, the functional $\Phi_\e (\sigma,\xi)=f_\e
(z^{\e\xi,\sigma}+w(\e,\sigma,\xi))$ has  the same regularity as $w$ and
satisfies:   $$ 
\n \Phi_\e(\sigma_0,\xi_0)=0\quad \Longleftrightarrow\quad \n 
f_\e\left(z_{\xi_0}+w(\e,\sigma_0,\xi_0)\right)=0. 
$$ 
\end{lemma} 
\begin{proof} 
Let $P=P_{\e\xi,\sigma}$ denote the projection onto $(T_{z^{\e\xi,\sigma}} Z^\e)^\perp$. We 
want 
to find a 
solution $w\in (T_{z^{\e\xi,\sigma}} Z)^{\perp}$ of the equation 
$P\nabla f_\e(z^{\e\xi,\sigma} +w)=0$.  One has that $\n f_\e(z+w)= 
\n f_\e (z)+D^2 f_\e(z)[w]+R(z,w)$ with $\|R(z,w)\|=o(\|w\|)$, uniformly 
with respect to 
$z=z^{\e\xi,\sigma}$, for  $|\xi|\leq \overline{\xi}$. Using the notation introduced 
in the previous section, we are led to the equation: 
$$ 
L_{\e,\sigma,\xi}w + P\n f_\e (z)+PR(z,w)=0. 
$$ 
According to Lemma \ref{lemma:3.1}, this is equivalent to 
$$ 
w = N_{\e,\xi,\sigma}(w), \quad \mbox{where}\quad 
N_{\e,\xi,\sigma}(w)=-L_{\e,\sigma,\xi}^{-1}\left( P\n f_\e (z)+PR(z,w)\right). 
$$ 
{From} Lemma \ref{CS:lem:1} it follows that 
 
\begin{equation}\label{CS:eq:N} 
        \|N_{\e,\xi,\sigma}(w)\|\leq c_1 (\e|\n V(\e\xi)|++\e|\nabla
K(\e\xi)|+\e |J_A (\e\xi)|+\e^2)+ o(\|w\|).          \end{equation} 
Then one readily checks that $N_{\e,\xi,\sigma}$ is a contraction on some ball in 
$(T_{z^{\e\xi,\sigma}} Z^\e)^{\perp}$ 
provided that $\e>0$ is small enough and $|\xi|\leq \overline{\xi}$. 
Then there exists a unique $w$ such that $w=N_{\e,\xi,\sigma}(w)$.  Let us 
point out that we cannot use the Implicit Function Theorem to find 
$w(\e,\xi,\sigma)$, because the map
$(\e,u)\mapsto P\n f_\e (u)$ fails to be
$C^2$.  However, fixed $\e>0$ small, we can apply the Implicit 
Function Theorem to the map $(\xi,\sigma,w)\mapsto P\n f_\e (z^{\e\xi,\sigma} +
w)$.  Then, in particular, the function $w(\e,\xi,\sigma)$ turns out to be of
class  $C^1$ with respect to $\xi$ and $\sigma$.  Finally, it is a standard
argument, see  \cite{AmbBad1, ABC}, to check that the critical points of
$\Phi_\e (\xi,\sigma)=f_\e (z+w)$ give rise to critical points of $f_\e$. 
\end{proof} 
 
\begin{remark}\label{remark:psi}
Since $f_\e (z^{\e\xi,\sigma})$ is actually independent of $\sigma$, the implicit function $w$ is 
constant with respect to that variable. As a result, there exists a functional $\Psi_\e \colon \Rn 
\to \R$ such that
\[
\Phi_\e (\sigma,\xi)=\Psi_\e (\xi), \qquad \forall\sigma\in [0,2\pi], \quad \forall \xi\in\Rn.
\]
In the sequel, we will omit the dependence of $w$ on $\sigma$, even it is defined over 
$S^1\times\Rn$.
\end{remark}

\begin{remark}\label{CS:rem:w} 
From (\ref{CS:eq:N}) it immediately follows that: 
\begin{equation}        \label{CS:eq:w} 
        \|w\|\leq C \left(\e |\n V(\e\xi)|+\e|\n K(\e\xi)| +\e
|J_A (\e\xi)|+\e^2\right),  \end{equation} 
where $C>0$. 
 \end{remark} 
 
\noindent  The following result can be proved by adapting the same argument as
in \cite{ams}.

\begin{lemma}\label{CS:lem:Dw} 
One has that: 
\begin{equation}        \label{CS:eq:Dw} 
        \|\nabla_\xi w\|\leq c \left(\e |\n 
V(\e\xi)|+\e|\n K(\e\xi)|+\e |J_A(\e\xi)|+O(\e^2)\right)^\gamma,
\end{equation} 
where $\gamma=\min\{1,p-1\}$ and $c > 0$ is some constant.
\end{lemma}

\section{The finite dimensional functional}\label{CS:subs:3}

The purpose of this subsection is to give an explicit form to the
finite--dimensional functional $\Phi_\e (\sigma,\xi)=\Psi_\e (\xi)=f_\e
(z^{\e\xi,\sigma}+w(\e,\xi))$.

Recall the precise definition of $z^{\e\xi,\sigma}$ given in \eqref{def:zeta}.
For brevity, we set in the sequel $z=z^{\e\xi,\sigma}$ and $w =
w(\e,\xi)$.

Since $z$ satisfies \eqref{eqlimite}, we easily find the following relations:
\begin{equation}
\int_\Rn \bigg|\bigg( \frac{\n}{\I} - A(\e\xi)\bigg)z \bigg|^2 + |z|^2 +
V(\e\xi)|z|^2=\int_\Rn K(\e\xi)|z|^{p+1}
\end{equation}
\begin{multline}
\re\int_\Rn \bigg( \frac{\n}{\I} - A(\e\xi)\bigg)z \cdot \overline{\bigg(
\frac{\n}{\I} - A(\e\xi)\bigg)w}  +  \re\int_\Rn z\bar w  \\
+\re\int_\Rn V(\e\xi)z \bar w  = \re\int_\Rn K(\e\xi) |z|^{p-1}z\bar w \,. 
\end{multline}

Hence we get
\begin{multline}\label{CS:eq:Phi}
\Phi_\e (\sigma,\xi)=f_\e (z^{\e\xi,\sigma}+w(\e,\sigma,\xi))
=\\ =K(\e \xi) \left( \frac{1}{2}-\frac{1}{p+1}\right) \int_\Rn
|z|^{p+1}
+ \frac{1}{2} \int_\Rn |A(\e\xi)-A(\e x)|^2 z^2\,  
+\\ \re\int_\Rn (A(\e\xi) -A(\e x))z \cdot (A(\e\xi)-A(\e x))\bar w 
 +\e\re\int_\Rn \frac{1}{\I} z \bar w \, \operatorname{div}A(\e x)
\\ + \frac{1}{2}\int_\Rn \bigg| \bigg(\frac{\n}{\I} - A(\e x)\bigg)w \bigg|^2 
+\re\int_\Rn [V(\e x)-V(\e\xi)]z \bar w
\\ + \frac{1}{2}\int_\Rn [V(\e x)-V(\e\xi)] |w|^2
+ \frac{1}{2} \int_\Rn [V(\e x)-V(\e\xi)] z^2
\\ +\frac{1}{2} V(\e\xi) \int_\Rn |w|^2\,
- \frac{1}{p+1} \re  \int_\Rn K(\e x)(|z+w|^{p+1} -
|z|^{p+1} - (p+1) |z|^{p-1} z \bar w )\, 
\\ + \re K(\e\xi)\int_\Rn |z|^{p-1}z\bar w 
+O(\e^2) \,.
\end{multline}
Here we have used the estimate
\[
\int_\Rn \bigg(\frac{1}{2} K(\e x)-\frac{1}{p+1} K(\e\xi) \bigg) |z|^{p+1} \,
 = \bigg( \frac{1}{2} - \frac{1}{p+1} \bigg) \int_\Rn K(\e\xi) |z|^{p+1} \,
 + O(\e^2),
\]
which follows from the boundedness of $K''$. Since we know that
\begin{eqnarray*}
\alpha (\e\xi)&=&\left( \frac{1+V(\e\xi)}{K(\e\xi)}\right)^{\frac{1}{p-1}}\\
\beta (\e\xi)&=&\left( 1+V(\e\xi)\right)^\frac{1}{2},
\end{eqnarray*}
we get immediately
\begin{equation}
\int_\Rn |z^{\e\xi,\sigma}|^{p+1}  = C_0 \, \Lambda (\e\xi)
[K(\e\xi)]^{-1},
\end{equation}
where we define the auxiliary function
\begin{equation}\label{def:auxiliary}
\Lambda (x)=\frac{(1+V(x))^\theta}{K(x)^{2/(p-1)}},\qquad \theta =
\frac{p+1}{p-1}-\frac{n}{2},
\end{equation}
and $C_0 = \|U\|_{L^2}$.
Now one can estimate the various terms in \eqref{CS:eq:Phi} by  means of
\eqref{CS:eq:w} and \eqref{CS:eq:Dw}, to prove that
\begin{equation}\label{CS:eq:expansion}
\Phi_\e (\sigma,\xi) = \Psi_\e (\xi)= C_1 \, \Lambda (\e\xi) + O(\e).
\end{equation}
Similarly,
\begin{equation}\label{CS:eq:expansionD}
\nabla \Psi_\e (\xi )=C_1 \n \Lambda (\e\xi) + \e^{1+\gamma} O(1),
\end{equation}
where $C_1 = \left(
\frac{1}{2}-\frac{1}{p+1} \right)C_0$. We omit the details, which can be
deduced without effort from \cite{ams}.

\section{Statement and proof of the main results}

\bigskip

In this section we exploit the finite-dimensional reduction
performed in the previous section to find existence and multiple
solutions of (NLS).
Recalling Lemma~\ref{CS:lem:w}, we have to look for critical points
of $\Phi_\e$ as a function of the variables
$(\sigma,\xi)\in [0,2\pi ]\times \Rn$
(or, equivalently, $(\eta,\xi)\in S^1 \times \Rn$).

In what follows, we use the
following notation: given a set $\Omega \subset \Rn$ and a number $\rho > 0$,
\[
\Omega_\rho \overset{\rm def}{=} \{ x\in\Rn \mid \operatorname{dist}(x,\Omega) < \rho \}.
\]

We start with the following result, which deals with local extrema.

\begin{theorem}\label{pippo}
Suppose that (K1), (V1) and (A1) hold. Assume moreover that there is
a compact set $M \subset \Rn$ over which $\Lambda$ achieves an isolated strict local 
minimum with value $a$. By this we mean that for some $\delta > 0$,
\begin{equation}\label{CS:eq:slm}
b\; \overset{\rm def}{=} \inf_{x\in \de M_\delta} \Lambda (x) > a.
\end{equation}
Then there exists $\e_\delta > 0$ such that $(S_\e)$ has at least $\operatorname{cat} 
(M,M_\delta)$ (orbits of) solutions concentrating near $M_\delta$, for all $0 < \e < \e_\delta$.

Conversely, assume that $K$ is a compact set $\Rn$ over which $\Lambda$
achieves an isolated strict local maximum with value $b$, namely
for some $\gamma > 0$,
\begin{equation}\label{CS:eq:max}
a \; \overset{\rm def}{=} \inf_{x\in \de K_\gamma} \Lambda (x) < b.
\end{equation}
Then there exists $\e_\gamma > 0$ such that $(S_\e)$ has at
least $\operatorname{cat} (K,K_\gamma)$ (orbits of)
solutions concentrating near $K_\gamma$, for all $0 < \e < \e_\gamma$.

\end{theorem}

\begin{proof}
As in the previous theorem, one has $\Phi_\e (\eta,\xi)=\Psi_\e (\xi)$.
Now choose $\bar\xi > 0$ in such a way that $M_\delta \subset \{x\in\Rn \mid \ |x| < \bar\xi\}$. 
Define again $\overline{\Lambda}$ as in the proof of Theorem \ref{th:5.3}. Let
\begin{eqnarray*}
N^\e &=& \{\xi\in \Rn \mid \e\xi\in M\} \\
N_\delta^\e &=& \{ \xi\in \Rn \mid \e\xi\in M_\delta\}\\
\Theta^\e &=& \{ \xi\in \Rn \mid \Psi_\e (\xi) \leq C_1 \tfrac{a+b}{2}\}.
\end{eqnarray*}
From \eqref{CS:eq:expansion} we get some $\e_\delta > 0$ such that
\begin{equation}\label{CS:eq:inclusioni}
N^\e \subset \Theta^\e \subset N_\delta^\e ,
\end{equation}
for all $0 < \e < \e_\delta$. To apply standard category theory, we need to prove that 
$\Theta^\e$ is compact. To this end, as can be
readily checked, it suffices to prove that $\Theta^\e$ cannot touch $\de N_\delta^\e$. But if 
$\e\xi\in\de M$, one has $\overline{\Lambda} (\e\xi) \geq b$ by the very definition of $\delta$, 
and so
\[
\Psi_\e (\xi) \geq C_1 \overline{\Lambda} (\e\xi) + o_\e (1) \geq C_1 b + o_\e (1).
\]
On the other hand, for all $\xi\in \Theta^\e$ one has also $\Psi_\e (\xi) \leq C_1 \frac{a+b}{2}$. 
We can conclude from \eqref{CS:eq:inclusioni} and elementary properties of the 
Lusternik--Schnirel'man category that $\Psi_e$ has at least
\[
\operatorname{cat} (\Theta^\e,\Theta^\e) \geq \operatorname{cat} (N^\e , N_\delta^\e) = 
\operatorname{cat} (N,N_\delta)
\]
critical points in $\Theta^\e$, which correspond to at least $\operatorname{cat}(M,M_\delta)$ 
orbits of solutions to $(S_\e)$.
Now, let $(\eta^*,\xi^*)\in S^1 \times M_\delta$ a critical point of $\Phi_\e$. Hence this point  
$(\eta^*,\xi^*)$ localizes a solution $u_{\e,\eta^*,\xi^*} (x)=z^{\e\xi^*,\eta^*} (x) + 
w(\e,\eta^*,\xi^*)$ of $(S_\e)$.  Recalling the change of variable which allowed us to pass
from (NLS) to $(S_\e)$, we find that
\[
u_{\e,\eta^*,\xi^*} (x) \approx z^{\e\xi^*,\eta^*} (\tfrac{x-\xi^*}{\e}).
\]
solves (NLS). The concentration statement follows from standard arguments (\cite{ABC,ams}).
The Proof of the second part of Theorem \ref{th:5.3} follows with analogous
arguments.
\end{proof}

Theorem \ref{th:nond} is an immediate corollary of the previous one when $x_0$ is 
either a nondegenerate local maximum or minimum for $\Lambda$.
We remark that the case in which $\Lambda$ has a maximum cannot be handled
using a direct variational approach and the arguments in \cite{cingolani}
cannot be applied. 

To treat the general case, we need some more work.
In order to  present our main result, we need to introduce
some topological concepts.

\smallskip
Given a set $M \subset \Rn$, the {\em cup long} of $M$ is by definition
\[
\ell (M)=1+\sup \{k\in\mathbb{N} \mid (\exists \alpha_1,\dots,\alpha_n \in \check{H}^{*} 
(M)\setminus \{1\})(\alpha_1 \cup \dots \cup \alpha_k \neq 0)\}.
\]
If no such classes exists, we set $\ell (M)=1$. Here $\check{H}^{*}(M)$ is the Alexander 
cohomology of $M$ with real coefficients, and $\cup$ denotes the cup product. It is well known 
that $\ell (S^{n-1})=\cat (S^{n-1})=2$, and $\ell (T^n)=\cat (T^n) = n+1$, where $T^n$ is the 
standard $n$--dimensional torus. But in general, one has $\ell (M) \leq \cat (M)$.

The following definition dates back to Bott (\cite{bott}).

\begin{definition}
We say that $M$ is non-degenerate for a $C^2$ function
$I\colon \R^N \to \R$ if $M$ consists of Morse theoretically
non-degenerate critical points for the restriction $I_{|M^\bot}$.
\end{definition}

To prove our existence result, we use again Theorem \ref{th:chang}.

We are now ready to prove an existence and multiplicity result for (NLS).

\begin{theorem}\label{th:5.3}
Let (V1), (K1) and (A1) hold. If the auxiliary function $\Lambda$ has
a smooth, compact, non-degenerate manifold of critical points $M$,
then for $\e > 0$ small, the problem $(S_\e)$ has at least $\ell (M)$ (orbits of) solutions 
concentrating near points of $M$.
\end{theorem}
\begin{proof}
By Remark \ref{remark:psi}, we have to find critical points of
$\Psi_\e=\Psi_\e (\xi)$. Since $M$ is compact, we can choose $\bar\xi > 0$ so
that $|x| < \bar\xi$ for all points $x\in M$. From this moment, $\bar\xi$ is
kept fixed. 
form $\{\eta^*\} \times M$ is obviously  a non-degenerate critical manifold
We set now $V=\Rn$, $J=\Lambda$, $\Sigma = M$, and $I(\xi)=\Psi_\e (\eta,\xi /
\e)$. Select $\delta > 0$ so that $M_\delta \subset \{x \colon |x| < \bar\xi
\}$, and no critical points of $\Lambda$ are in $M_\delta$, except fot those
of $M$. Set $\mathcal{U}=M_\delta$. From \eqref{CS:eq:expansion} and
\eqref{CS:eq:expansionD} it follows that $I$ is close to $J$ in $C^1
(\overline{\mathcal{U}})$ when $\e$ is very small. We can apply Theorem
\ref{th:chang} to find at least $\ell (M)$ critical points
$\{\xi_1,\dots,\xi_{\ell (M)}\}$ for $\Psi_\e$, provided $\e$ is small enough.
Hence the orbits $S^1\times\{\xi_1\}, \dots, S^1\times\{\xi_{\ell (M)}\}$
consist of critical points for $\Phi_\e$ which produce solutions of $(S_\e)$.
The concentration statement follows as in \cite{ams}. \end{proof}

\noindent
\begin{remark}\label{nonl}   
We point out that
Theorem 1.1, Theorem \ref{pippo} 
and Theorem \ref{th:5.3} hold for problems involving
more general nonlinearities $g(x,u)$
satisfying the same assumptions in \cite{Gr} (see also Remark 5.4
in \cite{ams}).
For our approach, we need the uniqueness of the radial solution
$z$ of the corresponding scalar equation
\begin{equation}\label{eta}
- \varDelta u + u + V(\e \xi)u=g(\e\xi,u), \quad u>0, \ \ u \in W^{1,2}(\Rn) \,.
\end{equation}
Let us also remark that in \cite{cingolani} the class
of nonlinearities handled does not require that
equation \eqref{eta} has a unique solution.
\end{remark}

\chapter{Homoclinic solutions of Hamiltonian systems}

A different tool for analyzing non--compact problems has been studied by Rabier and Stuart (\cite{rs,rs1}). It makes use of a recent notion of topological degree, introduced by Fitzpatrick, Pejsachowicz and Rabier, for $C^1$ Fredholm maps between Banach spaces.
In this chapter we summarize some results obtained in the joint paper \cite{SeStu}.

\vspace{2cm}

\section{Basic definitions of the new topological degree}

Consider two real Banach spaces $X$ and $Y$. First of all, one defines the
\textbf{parity} of a continuous path $\lambda\in\lbrack a,b]\mapsto
A(\lambda)$ of bounded linear Fredholm operators with index zero from $X$ into
$Y$. It is always possible to find a \textbf{parametrix} for this path, namely
a continuous function $B\colon\lbrack a,b]\rightarrow GL(Y,X)$ such that the
composition $B(\lambda)A(\lambda)\colon X\rightarrow X$ is a compact
perturbation of the identity for every $\lambda\in\lbrack a,b]$. If $A(a)$ and
$A(b)$ belong to $GL(X,Y)$, then the parity of the path $A$ on $[a,b]$ is by
definition
\[
\pi(A(\lambda)\mid\lambda\in\lbrack a,b])=\deg(B(a)A(a))\deg(B(b)A(b)).
\]
This is a good definition in the sense that it is independent of the
parametrix $B$. The following criterion can be useful for evaluating the
parity of an admissible path.

\begin{proposition}\label{prop:crit}
\label{prop2.1} Let $A\colon\lbrack a,b]\rightarrow B(X,Y)$ be a continuous
path of bounded linear operators having the following properties.

\begin{enumerate}
\item [(i)]$A\in C^{1}([a,b],B(X,Y))$.

\item[(ii)] $A(\lambda)\colon X\rightarrow Y$ is a Fredholm operator of index
zero for each $\lambda\in\lbrack a,b]$.

\item[(iii)] There exists $\lambda_{0}\in(a,b)$ such that
\begin{equation}
A^{\prime}(\lambda_{0})[\ker A(\lambda_{0})]\oplus\rge A(\lambda
_{0})=Y\label{2.1}%
\end{equation}
in the sense of a topological direct sum.
\end{enumerate}

Then there exists $\varepsilon>0$ such that $[\lambda_{0}-\varepsilon
,\lambda_{0}+\varepsilon]\subset\lbrack a,b]$,
\begin{equation}
A(\lambda)\in GL(X,Y)\text{ for }\lambda\in\lbrack\lambda_{0}-\varepsilon
,\lambda_{0})\cup(\lambda_{0},\lambda_{0}+\varepsilon]\label{2.2}%
\end{equation}
and
\begin{equation}
\pi(A(\lambda)\mid\lambda\in\lbrack\lambda_{0}-\varepsilon,\lambda
_{0}+\varepsilon])=(-1)^{k}\label{2.3}%
\end{equation}
where $k=\dim\ker A(\lambda_{0})$.
\end{proposition}

\noindent The proof of this proposition is essentially contained in
\cite{fitz,fp}.

We remark that given a continuous path $A\colon\lbrack a,b]\rightarrow\Phi
_{0}(X,Y)$ and any $\lambda_{0}\in\lbrack a,b]$ such that $A(\lambda)\in
GL(X,Y)$ for all $\lambda\neq\lambda_{0}$, the parity $\pi(A(\lambda
)\mid\lambda\in\lbrack\lambda_{0}-\varepsilon,\lambda_{0}+\varepsilon])$ is
the same for all $\varepsilon>0$ sufficiently small. This number is then
called the parity of $A$ across $\lambda_{0}$.

This degree has many features of the ordinary Leray--Schauder degree, but the homotopy invariance is only up to a $\pm$ sign. We refer to \cite{SeStu} for a short review of this concept.

We recall here that a bounded linear operator $L\colon X\to Y$ acting between two Banach spaces $X$ and $Y$ is said to be \emph{Fredholm of index zero} if its range $L(X)$ is closed in $Y$, $\ker L$ is finite--dimensional and $\dim\ker L=\operatorname{codim} L(X)$. 

\begin{definition}
Let $X$ and $Y$ be real Banach spaces and consider a function $F\in C^1 (\Lambda \times X,Y)$, where $\Lambda$ is an open interval. Let $P(\lambda,x)=\lambda$ be the projection of $\R\times X$ onto $\R$. We say that $\Lambda$ is an admissible interval for $F$ provided that
\begin{description}
\item{(i)} for all $(\lambda ,x)\in \Lambda \times X$, the bounded linear operator $\de_x F(\lambda , x)\colon X\to Y$ is a Fredholm operator of index zero;
\item{(ii)} for any compact subset $K\subset Y$ and any closed bounded subset $W$ of $\R\times X$ such that
\[
\inf \Lambda < \inf PW \leq \sup PW < \sup \Lambda.
\]
it results that $F^{-1}(K)\cap W$ is a compact subset of $\R\times X$.
\end{description}
\end{definition}

In figure 6.1, the set $F^{-1}(K)$ is not compact. Nevertheless, the shaded region represents the compact set $F^{-1}(K)\cap W$.

\medskip

\begin{figure}
\centering
\includegraphics[height=70mm]{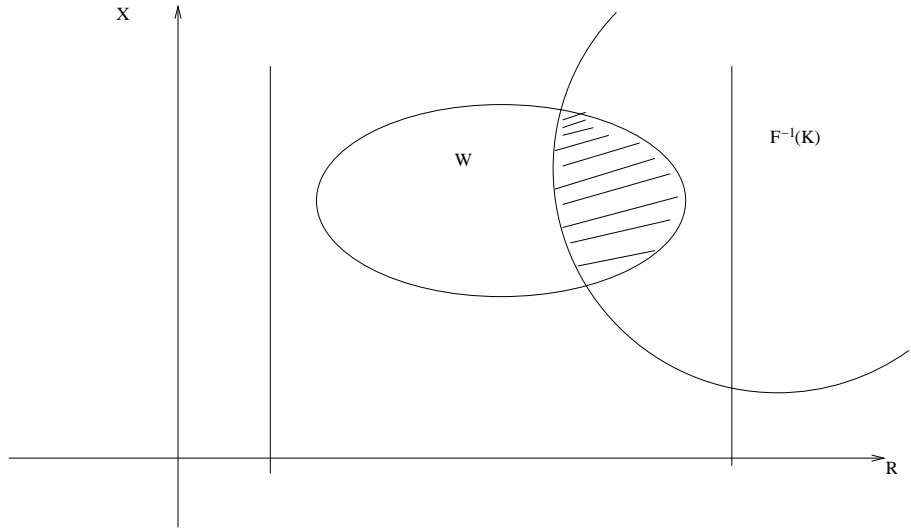}
\caption{Properness on bounded subsets}
\end{figure}

The following theorem shows that the new degree $\pi$ for Fredholm maps of index zero allows us to recast a bifurcation theorem like the one due to Rabinowitz.

\begin{theorem}\label{th:bif}
\label{th2.3} Let $X$ and $Y$ be real Banach spaces and consider a function
$F\in C^{1}(\Lambda\times X,Y)$ where $\Lambda$ is an admissible open interval
for $F$. Suppose that $\lambda_{0}\in\Lambda$ and that there exists
$\varepsilon>0$ such that $[\lambda_{0}-\varepsilon,\lambda_{0}+\varepsilon
]\subset\Lambda$,
\[
D_{x}F(\lambda,0)\in GL(X,Y)\text{ for }\lambda\in\lbrack\lambda
_{0}-\varepsilon,\lambda_{0}+\varepsilon]\setminus\{\lambda_{0}\}
\]
and
\[
\pi(D_{x}F(\lambda,0)\mid\lbrack\lambda_{0}-\varepsilon,\lambda_{0}%
+\varepsilon])=-1.
\]
Let $Z=\{(\lambda,u)\in\Lambda\times X\mid u\neq0$ and $F(\lambda,u)=0\}$ and
let $C$ denote the connected component of $Z\cup\{(\lambda_{0},0)\}$
containing $(\lambda_{0},0)$.

Then $C$ has at least one of the following properties:

\begin{description}
\item[(1)]$C$ is unbounded.

\item[(2)] The closure of $C$ contains a point $(\lambda_{1},0)$ where
$\lambda_{1}\in\Lambda\setminus\lbrack\lambda_{0}-\varepsilon,\lambda
_{0}+\varepsilon]$ and $D_{x}F(\lambda_{1},0)\notin GL(X,Y)$.

\item[(3)] The closure of $PC$ intersects the boundary of $\Lambda$.
\end{description}
\end{theorem}

\section{A concrete example}

In practice, the application of Theorem \ref{th:bif} consists of this: first, one identifies solutions of a given problem as the zeroes of a $C^1$ map $F$ depending on a parameter $\lambda$. Then one tries to single out admissible intervals for this map $F$. To this aim, one often uses the asymptotic behavior of the problem. Roughly speaking, if $F$ ``tends'' at infinity in a suitable sense to another map with a simpler structure (for example an autonomous structure), and this limit map has no non--trivial zeroes for a certain range $\Lambda$ of the parameter $\lambda$, then $\Lambda$ is admissible for $F$. The last step is to find a suitable value $\lambda_0\in\Lambda$ from which a global branch of solutions can bifurcate.

These ideas have been applied to nonlinear elliptic differential equations (\cite{rs1}), and we focus here on the case of hamiltonian systems, which present an easier technical framework.

Precisely, let us introduce the hamiltonian $H\colon \R\times \R^2 \times \R \to \R$ as
\begin{equation}
H(t,u,v,\lambda)=\frac{1}{2}\left\{ v^2 + \lambda u^2 + a(t)u^2\right\} - \frac{|u|^{\sigma+2}}{(\sigma+2)(1+\E^{-t})} + \frac{u^2 v^2}{2(\E^t +1)},
\end{equation}
where $\sigma >0$, $a\in C(\R)$ is non-negative, not identically zero, and 
\begin{equation}\label{eq:a}
\lim_{|t|\to+\infty}a(t)=0.
\end{equation}
We look for solutions of the system
\begin{equation}\label{eq:example}
\begin{cases}
Jx'=\nabla H(t,x,\lambda)\\
\lim_{|t|\to\infty} |x(t)|=0
\end{cases}
\end{equation}
where $x=(u,v)\in H^1(\R^2)$, or, more explicitely,
\begin{equation*}
\begin{cases}
-v'=\lambda u + a(t)u - \frac{1}{1+\E^{-t}}|u|^\sigma u\\
\phantom{-}u'=v+\frac{u^2 v}{1+\E^t}.
\end{cases}
\end{equation*}

Define
\begin{align*}
H^{+}(t,u,v,\lambda) &  =\frac{1}{2}\{v^{2}+\lambda u^{2}\}-\frac{\left|  u\right|  ^{\sigma+2}}{\sigma+2},\\
H^{-}(t,u,v,\lambda) &  =\frac{1}{2}\{v^{2}+\lambda u^{2}\}+\frac
{u^{2}v^{2}}{2},\\
g^{+}(u,v,\lambda) &  =(\lambda u-\left|  u\right|  ^{\sigma
}u,v)\\
g^{-}(u,v,\lambda) &  =(\lambda u+uv^{2},v+u^{2}v)\\
A_{\lambda}(t) &  =D_{(u,v)}^{2}H(t,0,0,\lambda)=\left(
\begin{array}
[c]{cc}%
\lambda+a(t) & 0\\
0 & 1
\end{array}
\right)\\
A_{\lambda}^{\pm}&  =D_{(u,v)}g^{\pm}(0,0,\lambda)=\left(
\begin{array}
[c]{cc}%
\lambda & 0\\
0 & 1
\end{array}
\right)  .
\end{align*}

The $g^\pm$ and $H^\pm$ verify

\begin{description}
\item [(1)]$g^{+}(t,0,\lambda)=g^{-}(t,0,\lambda)=0\text{ for all }%
t,\lambda\in\mathbb{R}$
\item[(2)] $\lim_{t\rightarrow\infty}\{D_{\xi}^{2}H(t,\xi,\lambda)-D_{\xi
}g^{+}(t,\xi,\lambda)\}=\lim_{t\rightarrow-\infty}\{D_{\xi}^{2}H(t,\xi
,\lambda)-\newline D_{\xi}g^{-}(t,\xi,\lambda)\}=0$, uniformly for $\xi=(u,v)$ in
bounded subsets of $\mathbb{R}^{2}$
\item[(3)] $\lim_{t\rightarrow\infty}\{H(t,\xi,\lambda)-H^{+}(t,\xi,\lambda)\}=\lim
_{t\rightarrow-\infty}\{H(t,\xi,\lambda)-H^{-}(t,\xi,\lambda)\}=0$
uniformly for $\xi$ in bounded subsets of $\mathbb{R}^{2}$
\end{description}
By \eqref{eq:a}, $H^\pm$ play the r\^{o}le of the {\it problems at infinity} we described above.

We will apply Theorem 6.1 with $X=H^1(\R,\R^2)$, $Y=L^2(\R,\R^2)$,
\[
F(\lambda,x)=Jx' - \nabla H(\cdot , x(\cdot) ,\lambda).
\]
We will also check that an admissible interval for this $F$ is $\Lambda = (-\infty,0)$.

We state and prove the main theorem concerning the system \eqref{eq:example}.

\begin{theorem}
Define
\begin{equation}\label{eq:deflambda}
\lambda_0 = \inf \left\{ \int_{\R} |\varphi'|^2 - a(t)|\varphi|^2 \mid \varphi\in H^1(\R) \land \int_{\R}|\varphi|^2 =1\right\}.
\end{equation}
Then a global branch of homoclinic solutions of \eqref{eq:example} bifurcates at $\lambda_0$ in the sense of Theorem \ref{th:bif} with $\Lambda=(-\infty,0)$.
\end{theorem}

\begin{proof}
We will apply the criterion of Proposition 6.1. 
It is immediate to prove that $F\in C^1 (\R\times H^1,L^2)$, from the very definition of the Fr\'{e}chet derivative.

We define $F^\pm\colon \R\times H^1 \to L^2$ as
\[
F^\pm(\lambda,x)=Jx' - \nabla H^\pm(\cdot , x(\cdot) ,\lambda),
\]
where, as above, $x=(u,v)$. Similarly, we define
\[
F^\pm(\lambda,x)=Jx' - \nabla H^\pm(\cdot , x(\cdot) ,\lambda).
\]
We claim that $\Lambda=(-\infty,0)$ \emph{is an admissible interval for} $F$. 

\smallskip

\textbf{First step:} $D_x F(\lambda ,0)\in \Phi_0(H^1,L^2)$ for all $\lambda <0$.

Let $L=D_x F(\lambda ,0)$, namely
\[
Lx(t)=Jx'(t) - A_\lambda (t)x(t)
\]
The linear operator $L$ is trivially bounded and also self--adjoint. From well--known results of linear functional analysis, $L$ is Fredholm of index zero provided $\rge (L)$ is closed.
With this in mind, let
$f\in L^{2}$ and suppose that there exists a sequence $\{f_{n}\}\subset \rge
L$ such that $\left\|  f-f_{n}\right\|  _{2}\rightarrow0.$ Let $\{x_{n}\}\subset H^{1}$ be such that $Lx_{n}=f_{n}.$

Let $P$ denote the orthogonal projection of $L^{2}$ onto $\ker L$ and set
$Q=I-P.$ Then $Qx_{n}=x_{n}-Px_{n}\in H^{1}$ and $Qx_{n}\in\lbrack
\ker(L)]^{\perp},$ the orthogonal complement of $\ker L$ in $L^{2}.$ Setting
\[
u_{n}=Qx_{n}\text{ we have that }u_{n}\in H^{1}\cap\lbrack\ker(L)]^{\perp
}\text{ and }Lu_{n}=f_{n}.
\]
Let us prove that the sequence $\{u_{n}\}$ is bounded in $H^{1}.$ For this we
use $S^{\pm}:H^{1}\rightarrow L^{2}$ to denote the bounded linear operators
defined by
\[
S^{\pm}x=Jx^{\prime}-A_{\lambda}^{\pm}x\text{ for all }x\in H^{1}.
\]
The two operators $S^{+}:H^{1}\rightarrow L^{2}$
and $S^{-}:H^{1}\rightarrow L^{2}$ are isomorphisms because
\[
A_\lambda^\pm =
\begin{pmatrix}
&\lambda &0\\
&0 &1
\end{pmatrix}
\]
and the spectrum of the matrix
\[
J
\begin{pmatrix}
&\lambda &0\\
&0 &1
\end{pmatrix}
\]
is
\[
\sigma\left(
J
\begin{pmatrix}
&\lambda &0\\
&0 &1
\end{pmatrix}
\right)=
\begin{cases}
\{ \pm \I \sqrt\lambda \} &\text{ if $\lambda >0$},\\
\{0\} &\text{ if $\lambda =0$},\\
\{\pm \sqrt{|\lambda|}\} &\text{ if $\lambda <0$}.
\end{cases}
\]
Hence this spectrum is, for $\lambda <0$, disjoint from $\I \R$, the imaginary axis; by a standard result (see \cite{hs,joost}), this implies that $S^\pm$ are both isomorphisms. 

So there exists a
constant $k$ such that
\begin{equation}
\left\|  S^{\pm}x\right\|  _{2}\geq k\left\|  x\right\|  \text{ for all }x\in
H^{1}.\label{5.1}%
\end{equation}
Supposing that $\left\|  u_{n}\right\|  \rightarrow\infty,$ we set
$w_{n}=\frac{u_{n}}{\left\|  u_{n}\right\|  }.$ Then $\{w_{n}\}\subset H^{1}$
with $\left\|  w_{n}\right\|  =1$ for all $n\in\mathbb{N}.$ Passing to a
subsequence, we can suppose that $w_{n}\rightharpoonup w$ weakly in $H^{1},$
and hence that $Lw_{n}\rightharpoonup Lw$ weakly in $L^{2}.$ Furthermore,
\[
Lw_{n}=\frac{Lu_{n}}{\left\|  u_{n}\right\|  }=\frac{f_{n}}{\left\|
u_{n}\right\|  }\text{ and so }\left\|  Lw_{n}\right\|  _{2}=\frac{\left\|
f_{n}\right\|  _{2}}{\left\|  u_{n}\right\|  }\rightarrow0
\]
since $\left\|  f_{n}\right\|  _{2}\rightarrow\left\|  f\right\|  _{2}$ and
$\left\|  u_{n}\right\|  \rightarrow\infty.$ Thus $Lw=0.$ But $w_{n}\in
H^{1}\cap\lbrack\ker(L)]^{\perp}$ for all $n\in\mathbb{N},$ from which it
follows that $w\in H^{1}\cap\lbrack\ker(L)]^{\perp}.$ Thus $w=0$ and
$w_{n}\rightharpoonup0$ weakly in $H^{1}.$ Consequently,
\[
w_{n}\rightarrow0\text{ uniformly on }[-R,R]\text{ for any }R\in(0,\infty).
\]
But, for all $t\in\mathbb{R},$%
\[
Jw_{n}^{\prime}(t)=A_{\lambda}(t)w_{n}(t)+\frac{f_{n}(t)}{\left\|
u_{n}\right\|  }%
\]
and so
\[
\left\|  w_{n}^{\prime}\right\|  _{L^{2}(-R,R)}\leq\sup_{t\in\mathbb{R}%
}\left\|  A_{\lambda}(t)\right\|  \left\|  w_{n}\right\|  _{L^{2}(-R,R)}%
+\frac{\left\|  f_{n}\right\|  _{2}}{\left\|  u_{n}\right\|  },
\]
showing that $\left\|  w_{n}^{\prime}\right\|  _{L^{2}(-R,R)}\rightarrow0$ as
$n\rightarrow\infty,$ for all $R\in(0,\infty).$ In particular,
\[
\left\|  w_{n}\right\|  _{H^{1}(-R,R)}\rightarrow0\text{ as }n\rightarrow
\infty.
\]

Now choose any $\varepsilon>0.$ There is a constant
$r\in(0,\infty)$ such that
\[
\left|  A_{\lambda}(t)-A_{\lambda}^{+}\right|  \leq\varepsilon\text{ for
all }t\geq r\text{ and }\left|  A_{\lambda}(t)-A_{\lambda}^{-}\right|
\leq\varepsilon\text{ for all }t\leq-r.
\]
Choose a constant $R>r+\frac{1}{\varepsilon}$ and a function $\varphi\in
C^{1}(\mathbb{R})$ such that

\medskip

\noindent$0 \leq\varphi(t)\leq1$ for all $t\in\mathbb{R}$, $\varphi(t)=0$ for
$t\leq r$, $\varphi(t)=1$ for $t\geq R$ and $\left|  \varphi^{\prime
}(t)\right| \leq\varepsilon$ for all $t\in\mathbb{R}$.

\medskip

Consider now the function $z_{n}(t)=\varphi(t)w_{n}(t).$ Clearly $z_{n}\in
H^{1}$ and
\begin{align*}
S^{+}z_{n}(t) &  =\varphi^{\prime}(t)Jw_{n}(t)+\varphi(t)Jw_{n}^{\prime
}(t)-A_{\lambda}^{+}z_{n}(t)\\
&  =\varphi^{\prime}(t)Jw_{n}(t)+\varphi(t)Lw_{n}(t)+\varphi(t)\{A_{\lambda
}(t)-A_{\lambda}^{+}\}w_{n}(t)\\
&  =\varphi^{\prime}(t)Jw_{n}(t)+\varphi(t)\frac{f_{n}(t)}{\left\|
u_{n}\right\|  }+\varphi(t)\{A_{\lambda}(t)-A_{\lambda}^{+}\}w_{n}(t).
\end{align*}
Thus,
\begin{align*}
\left\|  S^{+}z_{n}\right\|  _{2} &  \leq\varepsilon\left\|  w_{n}\right\|
_{2}+\frac{\left\|  f_{n}\right\|  _{2}}{\left\|  u_{n}\right\|  }+\sup_{t\geq
r}\left|  A_{\lambda}(t)-A_{\lambda}^{+}\right|  \left\|  w_{n}\right\|
_{2}\\
&  \leq\varepsilon+\frac{\left\|  f_{n}\right\|  _{2}}{\left\|  u_{n}\right\|
}+\varepsilon
\end{align*}
since $\left\|  w_{n}\right\|  _{2}\leq1.$ Hence, by (\ref{5.1})
\[
\left\|  w_{n}\right\|  _{H^{1}(R,\infty)}=\left\|  z_{n}\right\|
_{H^{1}(R,\infty)}\leq\left\|  z_{n}\right\|  \leq\frac{1}{k}\{2\varepsilon
+\frac{\left\|  f_{n}\right\|  _{2}}{\left\|  u_{n}\right\|  }\}.
\]
A similar argument, using $\varphi(-t)w_{n}(t)$ instead of $z_{n},$ shows
that
\[
\left\|  w_{n}\right\|  _{H^{1}(-\infty,-R)}\leq\frac{1}{k}\{2\varepsilon
+\frac{\left\|  f_{n}\right\|  _{2}}{\left\|  u_{n}\right\|  }\}.
\]
Finally, we have shown that, for all $n\in\mathbb{N},$%
\begin{align*}
\left\|  w_{n}\right\|  ^{2} &  =\left\|  w_{n}\right\|  _{H^{1}(-\infty
,-R)}^{2}+\left\|  w_{n}\right\|  _{H^{1}(-R,R)}^{2}+\left\|  w_{n}\right\|
_{H^{1}(R,\infty)}^{2}\\
&  \leq\frac{2}{k^{2}}\{2\varepsilon+\frac{\left\|  f_{n}\right\|  _{2}%
}{\left\|  u_{n}\right\|  }\}^{2}+\left\|  w_{n}\right\|  _{H^{1}(-R,R)}^{2}%
\end{align*}
and, letting $n\rightarrow\infty,$%
\[
\limsup_{n\rightarrow\infty}\left\|  w_{n}\right\|  ^{2}\leq\frac{2}{k^{2}%
}\{2\varepsilon\}^{2}%
\]
since $\left\|  w_{n}\right\|  _{H^{1}(-R,R)}^{2}\rightarrow0,\left\|
f_{n}\right\|  _{2}\rightarrow\left\|  f\right\|  _{2}$ and $\left\|
u_{n}\right\|  \rightarrow\infty.$ But $\left\|  w_{n}\right\|  \equiv1$ and
$\varepsilon>0$ can be chosen so that $\frac{2}{k^{2}}\{2\varepsilon\}^{2}<1.$
This contradiction establishes the boundedness of the sequence $\{u_{n}\}$ in
$H^{1}.$

By passing to a subsequence, we can now suppose that $u_{n}\rightharpoonup u$
weakly in $H^{1}$, and consequently that $Lu_{n}\rightharpoonup Lu$ weakly in
$L^{2}.$ However, $Lu_{n}=f_{n}$ and $\left\|  f_{n}-f\right\|  _{2}%
\rightarrow0,$ showing that $Lu=f$. This proves that $\rge L$ is a closed
subspace of $L^{2}$ and we have shown that $D_{x}F(\lambda,0)\in\Phi_{0}%
(H^{1},L^{2}).$

\smallskip

\textbf{Second step:} the properness of $F$.

In the sequel, we will use the following statement, whose proof is omitted and can be found in our paper \cite{SeStu}: \emph{if $D_x F(\lambda ,0)\in \Phi_0 (H^1,L^2)$, then $D_x F(\lambda,x)\in \Phi_0(H^1,L^2)$ for all $x\in H^1$}.

The next step is to show that $F$ is a proper map on the class of closed and bounded subsets of $\Lambda \times H^1$. Hence, let us fix a compact subset $K\subset L^2$ and a closed, bounded subset $W\subset \Lambda\times H^1$, in such a way that
\[
-\infty < \inf PW \leq \sup PW <0.
\]
Choose any sequence $\{(\lambda_n,x_n)\}$ from $F^{-1}(K)\cap W$; we may assume that
\[
\lambda_n\to\lambda\leq \sup PW <0,
\]
\[
\|F(\lambda_n,x_n)-y\|_2\to 0
\]
for some $y\in L^2$. Now, take any $x\in W$ and any $\lambda ,\mu$ in $\Lambda$. Then
\begin{equation*}
\|F(\lambda,x)-F(\mu,x)\|_2 = \|D_x H(\cdot ,x,\lambda)-D_x H(\cdot ,x,\mu)\|_2
\end{equation*}
and 
\begin{equation*}
D_x H(\cdot ,x,\lambda)-D_x H(\cdot ,x,\mu) = (\lambda -\mu) \int_0^1 D_\lambda D_{(u,v)} H(\cdot, x, s\lambda + (1-s)\mu)\, ds.
\end{equation*}
Set $L=\sup_{x\in W} \|x\|_\infty <\infty$, since $W$ is bounded in $H^1$. If $K=\sqrt{L+|\lambda|+1}$, then
\[
\|F(\lambda,x)-F(\mu,x)\|_2 \leq C(K) \|x\|_2 |\lambda-\mu|
\]
for some constant $C(K)>0$ depending only on $K$. We have thus proved that the family $\{F(\cdot,x_n)\}_{n\geq 1}$ is equicontinuous at every point $\lambda\in\Lambda$. This immediately implies that
\[
\|F(\lambda,x_n)-y\|_2\to 0.
\]
Now we ``freeze'' $\lambda$. Since $W$ is bounded, we may assume that
\[
x_n\weakto x\in H^1.
\]

\medskip

\emph{We claim that: for all $\ge >0$ there exists $R>0$ such that for all integers $n\geq 1$ and all $t\in\R$ such that $|t|\geq R$, there results
\[
|x_n(t)|\leq \ge.
\]}

\smallskip

Suppose this is not true. By a simple reasoning, the sequence $\{x_n\}$ must slip off to infinity. More precisely, one of the following conditions must be true:
\begin{description}
\item[(1)] There is a sequence $\{l_{k}\}\subset\mathbb{Z}$ with $\lim_{k\rightarrow
\infty}l_{k}=\infty$ and a subsequence $\{x_{n_{k}}\}$ of $\{x_{n}\}$ such
that $\widetilde{x_{k}}=\tau_{l_{k}}(x_{n_{k}})=x_{n_{k}}(\cdot
+l_{k})$ converges weakly in $H^{1}$ to an element $\widetilde{x}\neq0$
\item[(2)] There is a sequence $\{l_{k}\}\subset\mathbb{Z}$ with $\lim_{k\rightarrow
\infty}l_{k}=-\infty$ and a subsequence $\{x_{n_{k}}\}$ of $\{x_{n}\}$ such
that $\widetilde{x_{k}}=\tau_{l_{k}}(x_{n_{k}})=x_{n_{k}}(\cdot
+l_{k})$ converges weakly in $H^{1}$ to an element $\widetilde{x}\neq0$ 
\end{description}

To fix ideas, let us suppose that
$\{x_{n}\}$ has the property (1) . The
invariance by translation of the Lebesgue measure implies that $\left\|
F(\lambda,x_{n})-y\right\|  _{2}=\left\|  \tau_{l_{k}}(F(\lambda
,x_{n})-y)\right\|  _{2}$ so
\[
\left\|  \tau_{l_{k}}(F(\lambda,x_{n_{k}}))-\widetilde{y_{k}}\right\|
_{2}\rightarrow0\text{ where }\widetilde{y_{k}}(t)=y(t+l_{k})
\]
For any $\omega\in(0,\infty),$ it is easy to show that
\[
\left\|  \tau_{l_{k}}(F(\lambda,x_{n_{k}}))-\tau_{l_{k}}%
(F^{+}(\lambda,x_{n_{k}}))\right\|  _{L^{2}(-\omega,\omega)}\rightarrow0\text{
as }k\rightarrow\infty.
\]
Hence
\[
\left\|  \tau_{l_{k}}(F^{+}(\lambda,x_{n_{k}}))-\widetilde{y_{k}%
}\right\|  _{L^{2}(-\omega,\omega)}\rightarrow0\text{ as }k\rightarrow\infty.
\]
But
\begin{align*}
\tau_{l_{k}}(F^{+}(\lambda,x_{n_{k}}))(t)  & =Jx_{n_{k}}^{\prime}%
(t+l_{k})-g^{+}(x_{n_{k}}(t+l_{k}),\lambda)\\
& =J\widetilde{x_{k}^{\prime}}(t)-g^{+}(\widetilde{x_{k}}(t),\lambda
)=F^{+}(\lambda,\widetilde{x_{k}})(t)
\end{align*}
because $g^{+}$ is independent of $t$. Consequently,
\[
\left\|  F^{+}(\lambda,\widetilde{x_{k}})-\widetilde{y_{k}}\right\|
_{L^{2}(-\omega,\omega)}\rightarrow0\text{ as }k\rightarrow\infty,
\]
for all $\omega\in(0,\infty).$ Since the sequence $\{F^{+}(\lambda
,\widetilde{x_{k}})-\widetilde{y_{k}}\}$ is bounded in $L^{2},$ this implies
that $F^{+}(\lambda,\widetilde{x_{k}})-\widetilde{y_{k}}\rightharpoonup0$
weakly in $L^{2}.$ But $\widetilde{y_{k}}\rightharpoonup0$ weakly in $L^{2}$
by standard results in functional analysis, so we now have that $F^{+}(\lambda,\widetilde{x_{k}%
})\rightharpoonup0$ weakly in $L^{2}.$ However, the weak sequential continuity
of $F^{+}(\lambda,\cdot):H^{1}\rightarrow L^{2}$ implies that
\[
F^{+}(\lambda,\widetilde{x_{k}})\rightharpoonup F^{+}(\lambda,\widetilde
{x})\text{ weakly in }L^{2},
\]
so we must have $F^{+}(\lambda,\widetilde{x})=0$. We now show that this implies $x\equiv 0$, the desired contradiction. Notice that $x\in C^1 (\R)$, and set
\[
C=
\begin{pmatrix}
&0 &1\\
&1 &0
\end{pmatrix}
.
\]
Then
\begin{equation*}
\frac{d}{dt} \langle Cx(t),x(t)\rangle = 2 \langle Cx(t),x' (t)\rangle =\lambda u(t)^2 - |u|^{\sigma +2} + v^2 >0
\end{equation*}
whenever $x=(u,v)$ is not identically zero. Since $\lim_{|t|\to\infty} |x(t)|=0$ because $x\in H^1$, then $\langle Cx(t),x(t)\rangle \equiv 0$. This implies that $\langle JCx(t),g^+(x(t),\lambda)\rangle =0$ for all $t\in\R$. Finally, $x(t)=0$ for all $t$.\footnote{Of course, the same reasoning would apply to case (2).}

By combining the Sobolev imbedding $H^1 \subset L^\infty$ on bounded sets and the claim we have just proved, it is immediate to obtain that
\[
\|x_n-x\|_\infty \to 0.
\]

\medskip

Hence we have proved that $F(\lambda ,\cdot)$ is a proper map on the bounded subsets of $H^1$. Hence there must be a subsequence, still denoted by $\{x_n\}$, such that
\[
x_n\to x \text{ strongly in $H^1$.}
\]
But then $(\lambda_n,x_n)\to (\lambda,x)$ strongly, and in particular $(\lambda,x)\in W$ because $W$ is closed. This finally proves that $F^{-1}(K)\subset W$ is compact in $H^1$, and also that $\Lambda=(-\infty,0)$ is admissible for $F$.

\medskip

\textbf{Third step:} the choice of $\lambda_0$.

It is well-known (see \cite{Lieb}, theorem 11.5) that $\lambda_{0}$ defined in \eqref{eq:deflambda} belongs to $(-\infty,0)$ and that there 
exists an
element $\varphi_{0}\in H^{1}(\mathbb{R})$ such that
\[
\varphi_{0}(t)>0\text{ for all }t\in\mathbb{R}\text{ and }\lambda_{0}%
\int_{-\infty}^{\infty}\varphi(t)_{0}^{2}dt=\int_{-\infty}^{\infty}\varphi
_{0}^{\prime}(t)^{2}-a(t)\varphi_{0}(t)^{2}dt.
\]
Furthermore, $\varphi_{0}\in H^{2}(\mathbb{R})\cap C^{2}(\mathbb{R})$ and
satisfies the equation
\[
\varphi^{\prime\prime}(t)+\{\lambda_{0}+a(t)\}\varphi(t).
\]
Setting $x_{0}=(\varphi_{0},\varphi_{0}^{\prime}),$ we find that $\ker
D_{x}F(\lambda_{0},0)=\operatorname{span}\{x_{0}\}.$

Finally we observe that
\[
D_{\lambda}D_{(u,v)}^{2}H(t,0,0,\lambda)=\left[
\begin{array}
[c]{cc}%
1 & 0\\
0 & 0
\end{array}
\right]
\]
We are now ready to prove that the statements of Proposition \ref{prop:crit} apply with $A(\lambda)=D_xF(\lambda,0)$ and the selected value $\lambda_0$. Indeed, since $\ker A(\lambda_0)$ is one--dimensional,
\[
A'(\lambda_0)[\ker A(\lambda_0) \cap \rge A(\lambda_0) = \emptyset.
\]
Moreover, $\operatorname{codim}\rge A(\lambda_0) = \dim\ker A(\lambda_0)$ because $A(\lambda_0)$ is Fredholm of index zero. Since
\[
A'(\lambda_0)[\ker D_x F(\lambda_0,0)]=\left\{
\left(
\begin{array}
[c]{cc}%
1 & 0\\
0 & 0
\end{array}\right)
u\mid u\in \ker D_x F(\lambda_0,0) \right\},
\]
$A'(\lambda_0)[\ker A(\lambda_0)]$ and $\rge A(\lambda_0)$ are complementary in $L^2$:
\[
A'(\lambda_0)[\ker A(\lambda_0)] \oplus \rge A(\lambda_0) = L^2.
\]
In our situation, $\pi (D_x F(\lambda_0,0)\mid [\lambda_0-\ge,\lambda_0+\ge])=(-1)^{1}=-1$, and Theorem \ref{th:bif} applies.
\qed\end{proof}

\section{The general result}

The topological degree allows us to handle more general hamiltonian systems like

\begin{equation}
Jx'(t)=\nabla H(t,x(t),\lambda)
\end{equation}
where $x\in H^1(\R,\R^{2N})$, and $H\colon\mathbb{R}%
\times\mathbb{R}^{2N}\times\mathbb{R}\longrightarrow\mathbb{R}$ satisfies

\begin{description}
\item [(H1)]$H\in C(\mathbb{R\times R}^{2N}\times\mathbb{R})$ with
$H(t,\cdot,\lambda)\in C^{2}(\mathbb{R}^{2N})$ and $D_{\xi}H(t,0,\lambda)=0$
for all $t,\lambda\in\mathbb{R}$

\item[(H2)] The partial derivatives $D_{\xi}H,D_{\xi}^{2}H,D_{\lambda}D_{\xi
}H,D_{\lambda}D_{\xi}^{2}H$ and $D_{\xi}D_{\lambda}D_{\xi}H$ exist and are
continuous on $\mathbb{R\times R}^{2N}\times\mathbb{R}$.

\item[(H3)] For each $\lambda\in\mathbb{R},$ $D_{\xi}H(\cdot,\cdot
,\lambda):\mathbb{R}\times\mathbb{R}^{2N}\rightarrow\mathbb{R}^{2N}$ is a
$C_{\xi}^{1}-$bundle map, and $D_{\lambda}D_{\xi}^{2}H:\mathbb{R\times
}(\mathbb{R}^{2N}\times\mathbb{R})\mathbb{\rightarrow R}$ is a $C_{(\xi
,\lambda)}^{0}-$bundle map.

\item[(H4)] $D_{\xi}^{2}H(\cdot,0,0)$ and $D_{\lambda}D_{\xi}^{2}%
H(\cdot,0,0)\in L^{\infty}(\mathbb{R}).$

\item[(H$^{\infty}$)] For all $\lambda\in\mathbb{R},$ there exist two $C_{\xi}^{1}%
-$bundle maps $g^{+}(\cdot,\cdot,\lambda)$ and $g^{-}(\cdot,\cdot
,\lambda):\mathbb{R}\times\mathbb{R}^{2N}\mathbb{\rightarrow R}^{2N}$ such that

\begin{description}
\item [(1)]$g^{+}(t,0,\lambda)=g^{-}(t,0,\lambda)=0\text{ for all }%
t,\lambda\in\mathbb{R}$

\item[(2)] $\lim_{t\rightarrow\infty}\{D_{\xi}^{2}H(t,\xi,\lambda)-D_{\xi
}g^{+}(t,\xi,\lambda)\}=\lim_{t\rightarrow-\infty}\{D_{\xi}^{2}H(t,\xi
,\lambda)-\newline D_{\xi}g^{-}(t,\xi,\lambda)\}=0$, uniformly for $\xi$ in
bounded subsets of $\mathbb{R}^{2N}$

\item[(3)] $g^{+}(t+T^{+},\xi,\lambda)-g^{+}(t,\xi,\lambda)=g^{-}(t+T^{-}%
,\xi,\lambda)-g^{-}(t,\xi,\lambda)=0$ for some $T^{+},T^{-}>0$ and for all
$(t,\xi)\in\mathbb{R}\times\mathbb{R}^{2N}$.
\end{description}
\end{description}

We point out that (H$^\infty$) describes again the ``problems at infinity''.

\begin{theorem}[\cite{SeStu}]
\label{th5.5}Suppose that (H1) to (H4) and (H$^{\infty})$ are satisfied. An
open interval $\Lambda$ is admissible provided that, for all $\lambda
\in\Lambda,$ the following conditions are satisfied.

(1) The periodic, linear Hamiltonian systems
\[
Jx^{\prime}-A_{\lambda}^{+}(t)x=0\text{ and }Jx^{\prime}-A_{\lambda}^{-}(t)x=0
\]
have no characteristic multipliers on the unit circle.

(2) The asymptotic limit $g^{+}$ satisfies the following condition: If there is a real, symmetric $2N\times2N-$matrix $C$ such that
\[
\left\langle g^{+}(t,\xi,\lambda),JC\xi\right\rangle >0\text{ for all }\xi
\in\mathbb{R}^{2N}\backslash\{0\},
\]
then $\{x\in H^{1}:F^{+}(\lambda,x)=0\}=\{0\}.$

(3) The asymptotic limit $g^{-}$ satisfies the following condition: If there is a real, symmetric $2N\times2N-$matrix $C$ such that
\[
\left\langle g^{-}(t,\xi,\lambda),JC\xi\right\rangle >0\text{ for all }\xi
\in\mathbb{R}^{2N}\backslash\{0\},
\]
then $\{x\in H^{1}:F^{-}(\lambda,x)=0\}=\{0\}.$

(4) There is a point $\lambda_{0}\in\Lambda$ such that

(i) $k=\dim N(\lambda_{0})$ is odd where $N(\lambda)=\{u\in C^{2}%
(\mathbb{R},\mathbb{R}^{2N}):Ju^{\prime}(t)-A_{\lambda}(t)u(t)\equiv0$ and
$\lim_{\left|  t\right|  \rightarrow\infty}u(t)=0\},$

(ii) for every $u\in N(\lambda_{0})\backslash\{0\}$ there exists $v\in
N(\lambda_{0})$ such that
\[
\int_{-\infty}^{\infty}\left\langle T_{\lambda_{0}}(t)u(t),v(t)\right\rangle
dt\neq0\text{ where }T_{\lambda}(t)=D_{\lambda}D_{\xi}^{2}H(t,0,\lambda)\text{
and}%
\]

(iii) $\dim \{T_{\lambda_{0}}(\cdot)u:u\in N(\lambda_{0})\}=k.$

\noindent Then a global branch of homoclinic solutions of (6.9) bifurcates at $\lambda_{0}$ in the sense of Theorem 6.1
with $X=H^{1}$ and $Y=L^{2}.$
\end{theorem}

The proof of this theorem follows the same lines of the one presented for the concrete example. We refer to \cite{SeStu} for the details.

\part{Appendix}

\chapter{A useful toolkit}

Most of the concepts we recall here can be found on many textbooks in
\emph{Functional analysis, Theory of Sobolev spaces and Nonlinear
(functional) analysis}. When a result is useful but not quite standard,
we have tried to supply a precise reference.

\section{Integration by parts}

In Calculus, one studies the formula for integrating ``by parts". For
functions of several variables, there is good analog we state here.

\begin{proposition}
Let $\Omega\subset\Rn$ be a bounded domain with a regular boundary
$\de\Omega$.\footnote{At least of class $C^1$.} If $u$, $v\in
C^2(\Omega)\cap C^1(\overline{\Omega})$, then
\[
\int_\Omega u \nabla v = \int_{\de\Omega} u \frac{\de v}{\de \nu} \, d
\mathcal{H}^{n-1} - \int_\Omega \nabla u \cdot \nabla v,
\]
where $\mathcal{H}^{n-1}$ is the Hausdorff measure in dimension $n-1$.
\end{proposition}

\section{The differentiation of functional on Banach spaces}

\begin{definition}
Let $X$ be a Banach space. We say that the functional $J\colon X\to\R$ is differentiable at the point $u\in X$ if there exists a linear functional $x^\star \in X^\star$ with the following property:
\[
\lim_{\|h\|_X \to 0} \frac{J(u+h)-J(u)-x^\star (h)}{\|h\|_X}=0.
\]
In this case, we write $DJ(u)$ instead of $x^\star$. 
\end{definition}
See \cite{AmbProd} for a thorough study of differentiability in Banach and Hilbert spaces.

\begin{definition}
Let $X$ be a Banach space. We say that the functional $J\colon X\to\R$ is G-differentiable at the point $u\in X$ with \emph{G\^ateaux derivative} $f\in X^\star$ if for every $h\in X$
\[
\lim_{t\to 0} \frac{J(u+th)-J(u)-\langle f,th\rangle}{t}=0.
\]
\end{definition}

\begin{proposition}
If $J$ has a continuous G-derivative on an open set $U\subset X$ then $J\in C^1 (U,\R)$.
\end{proposition}

\begin{proposition}[\cite{Willem:minimax}]
\begin{description}
\item[(i)] Let $\Omega$ be an open subset of $\Rn$ and let $2 < p < \infty$. The functionals
\[
\psi (u)=\int_\Omega |u|^p, \quad \chi (u)=\int_\Omega |u^+|^p
\]
are of class $C^2 (L^p(\Omega),\R)$ and
\[
\langle \psi'(u),h\rangle = p\int_\Omega |u|^{p-2}uh,\quad \langle \chi'(u),h\rangle = p\int_\Omega |u^+|^{p-1}h.
\]
\item[(ii)] Let $2<p<\infty$ if $N=1,2$ and $2<p\leq 2^\star$ if $N\geq 3$. The functionals $\psi$ and $\chi$ are of class $C^2(H_0^1(\Omega),\R)$.
\item[(iii)] Let $N\geq 3$ and $p=2^\star$. The functionals $\psi$ and $\chi$ are of class $C^2(D_0^{1,2}(\Omega),\R)$.
\end{description}
\end{proposition}

\section{Some results in Functional Analysis}

\subsection{Weak convergence}

\begin{definition}
Let $X$ be a Banach space. We say that the sequence $\{x_n\}$ in $X$
converges weakly to $x$ if
\[
\Lambda (x_n) \to \Lambda (x)
\]
for all $\Lambda \in X^\star$.
\end{definition}

\begin{definition}
A functional $J\colon X\to\R$ is weakly lower semicontinuous on $X$ if
one can choose from any weakly convergent sequence $\{u_n\}$ in $X$,
\[
u_n \rightharpoonup u,
\]
a subsequence $\{u_{n_k}\}$ such that
\[
J(u)\leq \lim_{k\to\infty} J(u_{n_k}).
\]
\end{definition}

\subsection{Reflexivity}

Let $X$ be a Banach space. Let $X^{\star\star}$ denote the topological
\emph{bidual} space of $X$. Then there is a natural map $i\colon X\to
X^{\star\star}$ defined by $i_x(\Lambda)=\Lambda(x)$ for all $x\in X$
and $\Lambda\in X^\star$. It is an immediate consequence of the
Hahn--Banach theorem that $\|i_x\|= \|\Lambda\|$.

\begin{definition}
We say that $X$ is reflexive if the map $i$ is onto.
\end{definition}

\begin{theorem}
The closed unit ball of a Banach space $X$ is weakly compact if and
only if $X$ is reflexive.
\end{theorem}

\section{Sobolev spaces}

We mainly follow \cite{degio}

\begin{definition}
Let $\Omega\subset \Rn$ be an open set.
Let $W_{loc}^{1,1}(\Omega)$ be the set of functions $u\in
L_{loc}^1(\Omega)$ for which there exist $w_1,\ldots,w_n\in
L_{loc}^1(\Omega)$ such that
\[
-\int_\Omega u \frac{\de f}{\de x_j}\, dx = \int_\Omega f w_j \, dx
\]
for all $f\in C_0^\infty (\Omega)$ and $j=1,\dots ,n$. The functions $w_i$
($i=1,\dots ,n$) are uniquely determined, and called \emph{weak partial
derivatives} of $u$.
\end{definition}

\begin{definition}
For all $p\in [1,\infty]$, define
\[
W_{loc}^{1,p}(\Omega)=\left\{ u\in W_{loc}^{1,1}(\Omega) \colon u,
\frac{\de u}{\de x_1},\dots,\frac{\de u}{\de x_n}\in L_{loc}^p (\Omega)
\right\}.
\]
\end{definition}

\begin{definition}
For all $p\in [1,\infty]$, define the \emph{Sobolev space}
\[
W^{1,p}(\Omega)=\left\{ u\in W_{loc}^{1,1}(\Omega) \colon u,
\frac{\de u}{\de x_1},\ldots,\frac{\de u}{\de x_n}\in L^p (\Omega)
\right\}.
\]
\end{definition}

\begin{remark}
When $p=2$, it is customary to write $H^1(\Omega)$ instead of
$W^{1,2}(\Omega)$.
\end{remark}

\begin{theorem}
With respect to the norm
\[
\|u\|_{1,p}:= \left( \|u\|_{L^p(\Omega)}^p + \sum_{j=1}^n \|\de_j
u\|_{L^p(\Omega)}^p \right)^{1/p},
\]
the set $W^{1,p}(\Omega)$ is
\begin{enumerate}
\item a Banach space for $1\leq p \leq \infty$;
\item a separable Banach space for $1\leq p < \infty$;
\item a reflexive Banach space for $1<p<\infty$;
\item a separable Hilbert space for $p=2$.
\end{enumerate}
\end{theorem}

\begin{definition}
For all $p\in [1,\infty)$, define $W_0^{1,p}(\Omega)$ as the closure in
$W^{1,p}(\Omega)$ of $C_0^\infty (\Omega)$. Moreover, let
$H_0^1(\Omega)=W_0^{1,2}(\Omega)$.
\end{definition}

\begin{theorem}
The space $W_0^{1,p}(\Omega)$ is
\begin{enumerate}
\item a separable Banach space for $1\leq p < \infty$;
\item a reflexive and separable Banach space for $1<p<\infty$;
\item a separable Hilbert space for $p=2$.
\end{enumerate}
\end{theorem}

\begin{theorem}
Let $p\in [1,\infty)$, and let $u,v\in W^{1,p}(\Omega)$. Then
\[
u^+, u^-, |u|, \min\{u,v\}, \max\{u,v\}
\]
belong to $W^{1,p}(\Omega)$.
\end{theorem}

\subsection{Embeddings}

\begin{theorem}[Sobolev]
\begin{enumerate}
\item If $p=1$ and $n > 1$, then $W_0^{1,1}(\Omega) \subset L^\infty (\Omega)$ and
for all $u\in W_0^{1,1}(\Omega)$ there holds
\[
\|u\|_\infty \leq \prod_{j=1}^n \|\de_j u\|_1^{1/n} \leq \|\nabla u\|_1.
\]
\item If $1\leq p <n$, then $W_0^{1,p}(\Omega) \subset L^{p^\star}
(\Omega)$, and
for all $u\in W_0^{1,p}(\Omega)$ there holds
\[
\|u\|_{p^\star} \leq \frac{(n-1)p}{n-p}\prod_{j=1}^n \|\de_j u\|_p^{1/n}
\leq \frac{(n-1)p}{n-p} \|\nabla u\|_p.
\]
Here
\[
p^\star =
\begin{cases}
\frac{np}{n-p} &\text{ if } \quad p<n \\
\infty &\text{ if } \quad p=n.
\end{cases}
\]
is the \emph{critical Sobolev exponent}.
\end{enumerate}
\end{theorem}

\begin{theorem}[Poincar\'e's inequality]
For all $p\in [1,\infty)$ and for all $\Omega$ with finite Lebesgue
measure $|\Omega|$, we have
\[
\|u\|_{L^p(\Omega)} = \max \left\{ 1,\frac{(n-1)p}{n}\right\}
|\Omega|^{1/n} \prod_{j=1}^n \|\de_j u\|_{L^p(\Omega)}^{1/n}
\]
for all $u\in W_0^{1,p}(\Omega)$. In particular
\[
\|u\|_{L^p(\Omega)} = \max \left\{ 1,\frac{(n-1)p}{n}\right\}
|\Omega|^{1/n} \|\nabla u\|_{L^p(\Omega)}^{1/n}
\]
for all $u\in W_0^{1,p}(\Omega)$
\end{theorem}

\begin{definition}
Let $E\subset \Rn$ and let $\alpha \in (0,1]$. We say that
$u\colon E\to\R$ is H\"{o}lder of exponent $\alpha$, if there exists $c
\geq 0$ such that
\[
|u(x)-u(y)|\leq c |x-y|^\alpha
\]
for all $x,y\in E$.
We denote by $C^{0,\alpha}(E)$ the set of such functions.
\end{definition}

\begin{theorem}[Morrey]
Let $n < p < \infty$. Then $W_0^{1,p}(\Omega) \subset L^\infty (\Omega)$
and for any $u\in W_0^{1,p}(\Omega)$ there is precisely one
$\tilde{u}\in C^{0,\alpha}(\Rn)$ such that
\[
\tilde{u}(x)=u(x) \mbox{\quad for a.e. $x\in\Omega$.}
\]
\[
\tilde{u}(x)=0 \mbox{\quad for all } x\in \Rn\setminus \Omega.
\]
\[
\|\tilde{u}\|_\infty \leq \frac{p(n+1)}{p-n} \|u\|_{1,p}
\]
\[
|\tilde{u}(x)-\tilde{u}(y)| \leq \frac{2pn}{p-n} \|\nabla
u\|_{L^p(\Omega)} |x-y|^\alpha \mbox{\quad for all } x,y\in \Rn.
\]
\end{theorem}

\subsection{Compact embedding}

\begin{theorem}[Rellich, I]
Let $|\Omega| < \infty$ and $1\leq p < n$. Then the natural embedding
\[
W_0^{1,p}(\Omega) \to L^q (\Omega)
\]
is compact for all $q\in [1,p^\star)$.
\end{theorem}

\begin{theorem}[Rellich, II]
Let $|\Omega| < \infty$ and $n<p < \infty$. Then the natural embedding
\[
W_0^{1,p}(\Omega) \to L^\infty (\Omega)
\]
is compact.
\end{theorem}

\subsection{The space $D^{1,p}$}

Let $\Omega$ be an open subset of $\Rn$, and $p\in [1,n]$.

\begin{definition}
The Sobolev space $D^{1,p}(\Omega)$ is defined by
\[
D^{1,p}(\Omega)=\left\{ u\in L^{p^\star}(\Omega) \mid \partial_i u \in L^p (\Omega) \mbox{ for all $i=1,2,...,n$}\right\}.
\]
with the norm
\[
\|u\|_{D^{1,p}} = \|u\|_{L^{p^\star}} + \|\nabla u\|_{L^{p}}.
\]
Here $1/p^\star = 1/p - 1/n$. The space $D_0^{1,p}(\Omega)$ is the closure of $C_0^\infty (\Omega)$ with respect to the norm $\|\, \|_{D^{1,p}}$.
\end{definition}

For more details, see for example \cite{Adams}.

\section{The topological degree in finite dimension}

\begin{definition}
Let $\Omega\subset\Rn$ be an open bounded set, and let $f\colon \overline{\Omega}\to\Rn$ be a continuous function such that $f\in C^1 (\Omega)$. Fix any point $b\notin f(\partial \Omega)$ and any $0 < \ge < \operatorname{dist}(b,f(\partial \Omega))$. The \emph{topological degree of} $f$ with respect to $b\in \Rn\setminus f(\partial \Omega)$ is
\[
\deg (f,\Omega,b) = \int_\Omega \varphi (|f(x)-b|) J_f (x)\, dx,
\]
where $\phi \in C((0,\infty),\R)$ has support in $(0,\ge)$, $\int_\Rn \varphi (|x|)\, dx = 1$,  and $J_f$ is the jacobian of $f$.
\end{definition}

\begin{definition}
Let $\Omega\subset\Rn$ be an open bounded set, and let $f\colon \overline{\Omega}\to\Rn$ be a continuous function. Fix any point $b\notin f(\partial \Omega)$ and any sequence $\{f_k\}$ in $C^1 (\Omega) \cap C(\overline{\Omega})$ such that
\[
\lim_{k\to\infty} \|f_k-f\|_\infty =0 \mbox{\quad in } \overline{\Omega}.
\]
Then the \emph{topological degree} of $f$ with respect to the point $b$ is defined for all $k$ sufficiently large by the formula
\[
\deg (f,\Omega,b)=\lim_{k\to\infty} \deg (f_k,\Omega,b).
\]
\end{definition}

\begin{theorem}[Properties of the degree]\label{th:top1}
Let $\Omega$ be a bounded open subset of $\Rn$, and $b\in\Rn$. 
\begin{enumerate}
\item Let $f\in C(\overline{\Omega},\Rn)$ and $b\notin f(\partial\Omega)$. Then $\deg (f,\Omega,b)\in \mathbb{Z}$.
\item Assume that $f_1$, $f_2\colon \overline{\Omega}\to\Rn$ are continuous and $b\notin f_1(\partial\Omega)\cup f_2(\partial\Omega)$. If
\[
\|f_1-f_2\|_\infty \leq \frac{1}{4} \operatorname{dist}(b,f_1(\partial\Omega)\cup f_2(\partial\Omega)),
\]
then $\deg (f_1,\Omega,b)=\deg (f_2,\Omega,b)$.
\item Let $\in C(\overline{\Omega},\Rn)$. If $b$, $b'$ are in the same connected component of $\Rn\setminus f(\partial \Omega)$, then $\deg (f,\Omega,b)=\deg (f,\Omega,b')$.
\item Let $\Omega_1$,$\Omega_2$ be two bounded open subsets of $\Rn$ such that $\Omega_1 \cap \Omega_2 = \emptyset$, and let $f\in C(\overline{\Omega_1} \cup \overline{\Omega_2})$. If $b\in \Rn\setminus f(\partial\Omega_1)\cup f(\partial\Omega_2)$, then
\[
deg (f,\Omega_1\cup\Omega_2,b) = \deg (f,\Omega_1,b)+\deg (f,\Omega_2,b).
\]
\item Let $f\in C(\overline{\Omega},\Rn)$ and $K\subset \overline{\Omega}$ be a compact set. If $b\notin f(K)\cup f(\partial\Omega)$, then
\[
\deg (f,\Omega,b)=\deg (f,\Omega\setminus K,b).
\]
\item Let $f\in C(\overline{\Omega},\Rn)$ and $b\notin f(\partial\Omega)$. If $\deg (f,\Omega,b)\neq 0$, then $f^{-1}(\{b\})\neq \emptyset$.
\item Let $b\in\Rn$. Let $I\colon\Rn\to\Rn$ be the identity map. Then
\[
\deg (I,\Omega,b)=
\begin{cases}
1 &\text{ if } b\in\Omega,\\
0 &\text{ if } b\notin \Omega.
\end{cases}
\]
\item Let $H\colon \overline{\Omega}\times [0,1] \to\Rn$ be a continuous function, and $b\notin H(\partial\Omega \times [0,1])$. Then, for all $t\in [0,1]$,
\[
\deg (H(\cdot , t),\Omega,b)=\deg (H(\cdot ,0),\Omega,b).
\]
\item Let $f$,$g\in C(\overline{\Omega},\Rn)$ such that $f=g$ on $\partial\Omega$, and let $b\notin f(\partial\Omega)$.Then
\[
\deg (f,\Omega,b)=\deg (g,\Omega,b).
\]
\item Let $\Omega$ be symmetric with respect to the origin in $\Rn$. Let $f\in C(\overline{\Omega},\Rn)$ be an odd map such that $0\notin f(\partial\Omega)$. Then $\deg (f,\Omega,0)$ is an odd number provided $0\in\Omega$, and an even number provided $0\notin\Omega$.
\end{enumerate}
\end{theorem}

\section{The Leray--Schauder topological degree}

\begin{definition}
Let $X$ be a Banach space, and $\Omega\subset X$. If $T\colon \Omega \to X$ is a continuous map, we say that $T$ is \emph{compact} whenever $T(B)$ is relatively compact in $X$ for any bounded subset $B$ of $\Omega$.
\end{definition}

\begin{lemma}
Let $\Omega$ be a bounded open subset of $X$. If $T\colon \overline{\Omega}\to X$ is compact and it has no fixed points on $\partial \Omega$, then
\begin{equation*}
(\exists \ge >0)(\forall u)(u\in\partial\Omega \Rightarrow \|u-Tu\| \geq \ge).
\end{equation*}
\end{lemma}

\begin{lemma}\label{lem:top2}
Let $\Omega$ be a bounded open subset of $X$, and $T\colon\overline{\Omega}\to X$ be compact without fixed points on $\partial\Omega$. If $\ge >0$ satisfies $4\ge \leq \operatorname{dist}(0,(I-T)(\partial\Omega))$, then there exist a linear subspace $E_\ge \subset X$ of finite dimension and an operator $T_\ge \colon \overline{\Omega}\to E_\ge$ such that
\[
\|T_\ge u-Tu\| \leq \ge \qquad \forall u\in \overline{\Omega},
\]
\[
\|u-T_\ge u\| \geq 3 \ge, \qquad \forall u\in\partial\Omega.
\]
\end{lemma}

\begin{definition}
Let $X$ be a Banach space, $\Omega$ a bounded open subset of $X$, $T\colon \overline{\Omega}\to X$ a compact operator without fixed points on $\partial\Omega$. Let $\ge >0$, $E_\ge$ and $T_\ge \colon \overline{\Omega}\to E_\ge$ be as in Lemma \ref{lem:top2}. Consider $F$, a linear subspace of finite dimension containing $E_\ge$ such that $\Omega`_F = \Omega \cap F \neq \emptyset$.
The \emph{Leray--Schauder degree of} $T$ is defined as
\[
\deg (I_T,\Omega,0)=\deg_F (I_F - T_\ge ,\Omega_F ,0_F),
\]
where $I_F$ is the identity map in $F$, and $0_F$ is the origin in $F$.
\end{definition}

The Leray--Schauder degree has properties similar to those collected in Theorem \ref{th:top1}. Of course, all maps must be of the form $I-T$, where $T$ is compact.

We refer to \cite{kav93} for the proofs.

\backmatter
%
%
%

%
%

\printindex


\end{document}